\documentclass[leqno]{article}

\usepackage{amstex,amscd,amssymb}

\long\def\comment#1\endcomment{}

\def\newopp#1{\expandafter\newcommand\csname #1\endcsname{\operatorname{#1}}}
\def\renewopp#1{\expandafter\renewcommand\csname
#1\endcsname{\operatorname{#1}}} 

\def\newnott#1{\expandafter\newcommand\csname 
#1\endcsname{{\operatorname{#1}}}}
\def\renewnott#1{\expandafter\renewcommand\csname
#1\endcsname{{\operatorname{#1}}}} 

\def\newletter#1#2{\expandafter\newcommand\csname #1\endcsname{{\csname
                   #2\endcsname #1}}}
\def\newlletter#1#2{\expandafter\newcommand\csname #1#1\endcsname{{\csname
                   #2\endcsname #1}}}

\newcommand{\Arrow}{\longrightarrow}

\newcommand{\proof}[1][Proof]{\noindent {\it #1.\ }}
\def\endproof{\hfill \ensuremath{\square}\par}

\newcommand{\defn}{\noindent {\bf Definition.\ }}

\newcommand{\ex}{\noindent {\bf Example.\ }}

\newcommand{\rem}{\noindent {\bf Remark.\ }}

\newtheorem{prop}{Proposition}[section]
\newtheorem{theorem}{Theorem}[section]

\newtheorem{thm}{Theorem}    

\newenvironment{lemma}[1][]{\trivlist
   \item[\hskip \labelsep{\bfseries Lemma#1.}]\itshape}
   {\endtrivlist}

\newenvironment{corr}[1][]{\trivlist
   \item[\hskip \labelsep{\bfseries Corollary#1.}]\itshape}
   {\endtrivlist}

\makeatletter

\def\@seccntformat#1{\csname the#1\endcsname.\hskip 5pt}

\def\punkt{\refstepcounter{subsubsection}
           \noindent{\bf \thesubsubsection.\ }}

\@addtoreset{equation}{section}

\makeatother

\def\cdot{{\:\raisebox{3pt}{\text{\circle*{1.5}}}}}

\def\eqref#1{\thetag{\ref{#1}}}
\let\usualref\ref
\renewcommand{\ref}[1]{{\normalfont\usualref{#1}}}
\newcommand{\wt}[1]{{\widetilde{#1}}}
\newcommand{\wh}[1]{{\widehat{#1}}}
\renewcommand{\lim}{\operatornamewithlimits{lim}}
\newcommand{\llim}{\displaystyle\operatornamewithlimits{Lim}_{\longleftarrow}}

\renewcommand{\phi}{\varphi}
\newcommand{\eps}{\varepsilon}
\newcommand{\half}{\frac{1}{2}}
\newcommand{\halfi}{\frac{\sqrt{-1}}{2}}
\newcommand{\PT}{{\operatorname{pt}}}
\newcommand{\Eend}{{{\cal E}\!nd\:}}
\newcommand{\prehodge}{{{\cal W}{\cal H}odge}}
\newcommand{\can}{{\sf can}}
\newcommand{\id}{{\sf id}}

\newcommand{\cp}{{\C P}}
\newcommand{\calo}{{\cal O}}
\newcommand{\loc}[1]{#1_{loc}}
\newcommand{\tv}{{\overline{T}V}}
\newcommand{\tm}{{\overline{T}M}}
\newcommand{\cc}{\overline{\C}}

\newcommand{\wB}{\wt{\B}}
\newcommand{\wD}{\wt{D}}
\newcommand{\conj}{\overline{\ }}
\newcommand{\h}{{\Bbb H}}
\newcommand{\6}{\partial}
\newcommand{\del}{\delta}

\newletter{R}{Bbb}
\newletter{C}{Bbb}
\newlletter{C}{cal}
\newlletter{H}{cal}
\newletter{Q}{Bbb}
\newletter{Z}{Bbb}

\newletter{G}{cal}
\newletter{M}{cal}
\newletter{X}{frak}
\newletter{U}{cal}
\newlletter{R}{cal}
\newlletter{Z}{cal}
\newletter{D}{cal}
\newletter{B}{cal}
\newletter{I}{cal}
\newletter{N}{cal}
\newletter{T}{cal}
\newlletter{L}{cal}
\newletter{A}{cal}
\newletter{K}{cal}
\newletter{F}{cal}
\newletter{E}{cal}
\newletter{J}{frak}
\newlletter{J}{cal}
\newletter{W}{cal}
\newletter{V}{cal}

\newnott{SB}

\newopp{End}
\newopp{ev}
\newopp{Lin}
\newopp{Res}
\renewopp{Re}
\renewopp{Im}
\renewopp{dim}
\newopp{Spec}
\newopp{Ker}
\newopp{Bun}
\newopp{Hom}
\newopp{Ob}
\newopp{Der}
\newopp{Diff}
\newopp{Shv}
\newopp{Alt}
\newopp{Mult}
\newopp{gr}

\begin{document}

\title{Hyperk\"ahler structures on total spaces of holomorphic
cotangent bundles}

\author{D. Kaledin}

\date{{\em Independent University of Moscow}}

\maketitle

\section*{Introduction}

A hyperk\"ahler manifold is by definition a Riemannian manifold $M$
equipped with two anti-commuting almost complex structures $I$, $J$
parallel with respect to the Levi-Civita connection.  Hyperk\"ahler
manifolds were introduced by Calabi in \cite{C}. Since then they
have become the topic of much research. We refer the reader to
\cite{Bes} and to \cite{HKLR} for excellent overviews of the
subject.

Let $M$ be a hyperk\"ahler manifold. The almost complex structures
$I$ and $J$ generate an action of the quaternion algebra $\h$ in the
tangent bundle $\Theta(M)$ to the manifold $M$. This action is
parallel with respect to the Levi-Civita connection. Every
quaternion $h \in \h$ with $h^2 = -1$, in particular, the product $K
= IJ \in \h$, defines by means of the $\h$-action an almost complex
structure $M_h$ on $M$. This almost complex structure is also
parallel, hence integrable and K\"ahler. Thus every hyperk\"ahler
manifold $M$ is canonically K\"ahler, and in many different
ways. For the sake of convenience, we will consider $M$ as a
K\"ahler manifold by means of the complex structure $M_I$, unless
indicated otherwise.

One of the basic facts about hyperk\"ahler manifolds is that the
K\"ahler manifold $M_I$ underlying a hyperk\"ahler manifold $M$ is
canonically holomorphically symplectic. To see this, let $\omega_J$,
$\omega_K$ be the K\"ahler forms for the complex structures $M_J$,
$M_K$ on the manifold $M$, and consider the $2$-form $\Omega =
\omega_J + \sqrt{-1}\omega_K$ on $M$. It is easy to check that the
form $\Omega$ is of Hodge type $(2,0)$ for the complex structure
$M_I$ on $M$. Since it is obviously non-degenerate and closed, it is
holomorphic, and the K\"ahler manifold $M_I$ equipped with the form
$\Omega$ is a holomorphically symplectic manifold.

It is natural to ask whether every holomorphically symplectic
manifold $\langle M,\Omega\rangle$ underlies a hyperk\"ahler
structure on $M$, and if so, then how many such hyperk\"ahler
structures are there. Note that if such a hyperk\"ahler structure
exists, it is completely defined by the K\"ahler metric $h$ on
$M$. Indeed, the K\"ahler forms $\omega_J$ and $\omega_K$ are by
definition the real and imaginary parts of the form $\Omega$, and
the forms $\omega_J$ and $\omega_K$ together with the metric define
the complex structures $J$ and $K$ on $M$ and, consequently, the
whole $\h$-action in the tangent bundle $\Theta(M)$. For the sake of
simplicity, we will call a metric $h$ on a holomorphically
symplectic manifold $\langle M,\Omega\rangle$ hyperk\"ahler if the
Riemannian manifold $\langle M, h \rangle$ with the quaternionic
action associated to the pair $\langle \Omega, h \rangle$ is a
hyperk\"ahler manifold.

It is known (see, e.g., \cite{Beauv}) that if the holomorphically
symplectic manifold $M$ is compact, for example, if $M$ is a
$K3$-surface, then every K\"ahler class in $H^{1,1}(M)$ contains a
unique hyperk\"ahler metric. This is, in fact, a consequence of the
famous Calabi-Yau Theorem, which provides the canonical Ricci-flat
metric on $M$ with the given cohomology class. This Ricci-flat
metric turns out to be hyperk\"ahler. Thus in the compact case
holomorphically symplectic and hyperk\"ahler manifolds are
essentially the same.

The situation is completely different in the general case. For
example, all holomorphically symplectic structures on the formal
neighborhood of the origin $0 \in \C^{2n}$ in the $2n$-dimensional
complex vector space $\C^{2n}$ are isomorphic by the appropriate
version of the Darboux Theorem. On the other hand, hyperk\"ahler
structures on this formal neighborhood form an infinite-dimensional
family (see, e.g., \cite{HKLR}, where there is a construction of a
smaller, but still infinite-dimensional family of hyperk\"ahler
metrics defined on the whole $\C^{2n}$). Thus, to obtain meaningful
results, it seems necessary to restrict our attention to
holomorphically symplectic manifolds belonging to some special
class.

The simplest class of non-compact holomorphically symplectic
manifolds is formed by total spaces $T^*M$ to the cotangent bundle
to complex manifolds $M$. In fact, the first examples of
hyperk\"ahler manifolds given by Calabi in \cite{C} were of this
type, with $M$ being a K\"ahler manifold of constant holomorphic
sectional curvature (for example, a complex projective space). It
has been conjectured for some time that every total space $T^*M$ of
the cotangent bundle to a K\"ahler manifold admits a hyperk\"ahler
structure. The goal of this paper is to prove that this is indeed
the case, if one agrees to consider only an open neighborhood $U
\subset T^*M$ of the zero section $M \subset T^*M$.  Our main result
is the following.

\begin{thm}\label{th.1}
Let $M$ be a complex manifold equipped with a K\"ahler metric. The
metric on $M$ extends to a hyperk\"ahler metric $h$ defined in the
formal neighborhood of the zero section $M \subset T^*M$ in the
total space $T^*M$ to the holomorphic cotangent bundle to $M$. The
extended metric $h$ is invariant under the action of the group
$U(1)$ on $T^*M$ given by dilatations along the fibers of the
canonical projection $\rho:T^*M \to M$.  Moreover, every other
$U(1)$-invariant hyperk\"ahler metric on the holomorphically
symplectic manifold $T^*M$ becomes equal to $h$ after a holomorphic
symplectic $U(1)$-equivariant automorphism of $T^*M$ over $M$.
Finally, if the K\"ahler metric on $M$ is real-analytic, then the
formal hyperk\"ahler metric $h$ converges to a real-analytic metric
in an open neighborhood $U \subset T^*M$ of the zero section $M
\subset T^*M$.
\end{thm}

Many of the examples of hyperk\"ahler metrics obtained by
Theorem~\ref{th.1} are already known. (See, e.g., \cite{K1},
\cite{K2}, \cite{Nak}, \cite{H}, \cite{BG}, \cite{Sw}.) In these
examples $M$ is usually a generalized flag manifold or a homogeneous
space of some kind. On the other hand, very little is known for
manifolds of general type. In particular, it seems that even for
curves of genus $g \geq 2$ Theorem~\ref{th.1} is new.

We would like to stress the importance of the $U(1)$-invariance
condition on the metric in the formulation of
Theorem~\ref{th.1}. This condition for a total space $T^*M$ of a
cotangent bundle is equivalent to a more general compatibility
condition between a $U(1)$-action and a hyperk\"ahler structure on a
smooth manifold introduced by Hitchin (see, e.g., \cite{H}). Thus
Theorem~\ref{th.1} can be also regarded as answering a question of
Hitchin's in \cite{H}, namely, whether every K\"ahler manifold can
be embedded as the sub-manifold of $U(1)$-fixed points in a
$U(1)$-equivariant hyperk\"ahler manifold. On the other hand, it is
this $U(1)$-invariance that guarantees the uniqueness of the metric
$h$ claimed in Theorem~\ref{th.1}.

We also prove a version of Theorem~\ref{th.1} ``without the
metrics''. The K\"ahler metric on $M$ in this theorem is replaced
with a holomorphic connection $\nabla$ on the cotangent bundle to
$M$ without torsion and $(2,0)$-curvature. We call such connections
{\em K\"ah\-le\-ri\-an}. The total space of the cotangent bundle $T^*M$ is
replaced with the total space $\tm$ of the complex-conjugate to the
tangent bundle to $M$. (Note that {\em a priori} there is no complex
structure on $\tm$, but the $U(1)$-action by dilatations on this
space is well-defined.) The analog of the notion of a hyperk\"ahler
manifold ``without the metric'' is the notion if a hypercomplex
manifold (see, e.g., \cite{Bo}). We define a version of Hitchin's
condition on the $U(1)$-action for hypercomplex manifolds and prove
the following.

\begin{thm}\label{th.2}
Let $M$ be a complex manifold, and let $\tm$ be the total space of
the complex-conjugate to the tangent bundle to $M$ equipped with an
action of the group $U(1)$ by dilatation along the fibers of the
projection $\tm \to M$. There exists a natural bijection between the
set of all K\"ah\-le\-ri\-an connections on the cotangent bundle to $M$
and the set of all isomorphism classes of $U(1)$-equivariant
hypercomplex structures on the formal neighborhood of the zero
section $M \subset \tm$ in $\tm$ such that the projection $\rho:\tm
\to M$ is holomorphic. If the K\"ah\-le\-ri\-an connection on $M$ is
real-analytic, the corresponding hypercomplex structure is defined
in an open neighborhood $U \subset \tm$ of the zero section.
\end{thm}

Our main technical tool in this paper is the relation between
$U(1)$\--equi\-va\-ri\-ant hyperk\"ahler manifolds and the theory of
$\R$-Hodge structures discovered by Deligne and Simpson (see
\cite{De2}, \cite{Simpson}). To emphasize this relation, we use the
name {\em Hodge manifolds} for the hypercomplex manifolds equipped
with a compatible $U(1)$-action.

It must be noted that many examples of hyperk\"ahler manifolds
equipped with a compatible $U(1)$-action are already known. Such
are, for example, many of the manifolds constructred by the
so-called hyperk\"ahler redution from flat hyperk\"ahler spaces (see
\cite{H} and \cite{HKLR}). An important class of such manifolds is
formed by the so-called quiver varieties, studied by Nakajima
(\cite{Nak}). On the other hand, the moduli spaces $\M$ of flat
connections on a complex manifold $M$, studied by Hitchin
(\cite{Hcurves}) when $M$ is a curve and by Simpson (\cite{Smoduli}
in the general case, also belong to the class of Hodge manifolds, as
Simpspn has shown in \cite{Simpson}. Some parts of our theory,
especially the uniqueness statement of Theorem~\ref{th.1}, can be
applied to these known examples.

\medskip

We now give a brief outline of the paper. Sections 1-3 are
preliminary and included mostly to fix notation and
terminology. Most of the facts in these sections are well-known.
\begin{itemize}
\item In Section 1 we have collected the necessary facts from linear
algebra about quaternionic vector spaces and $\R$-Hodge
structures. Everything is standard, with an exception, perhaps, of
the notion of {\em weakly Hodge map}, which we introduce in
Subsection~\ref{w.H.sub}.

\item We begin Section 2 with introducing a technical notion of a
{\em Hodge bundle} on a smooth manifold equipped with a
$U(1)$-action. This notion will be heavily used throughout the
paper.  Then we switch to our main subject, namely, various
differential-geomteric objects related to a quaternion action in the
tangent bundle. The rest of Section 2 deals with almost quaternionic
manifolds and the compatibility conditions between an almost
quaternionic structure and a $U(1)$-action on a smooth manifold $M$.

\item In Section 3 we describe various integrability conditions on
an almost quaternionic structure. In particular, we recall the
definition of a hypercomplex manifold and introduce
$U(1)$-equivariant hypercomplex manifolds under the name of {\em
Hodge manifolds}. We then rewrite the definition of a Hodge manifold
in the more convenient language of Hodge bundles, to be used
throughout the rest of the paper. Finally, in
Subsection~\ref{polarization} we discuss metrics on hypercomplex and
Hodge manifolds. We recall the definition of a hyperk\"ahler
manifold and define the notion of a {\em polarization} of a Hodge
manifold. A polarized Hodge manifold is the same as a hyperk\"ahler
manifold equipped with a $U(1)$-action compatible with the
hyperk\"ahler structure of the sense of Hitchin, \cite{H}. 

\item The main part of the paper begins in Section 4. We start with
arbitrary Hodge manifolds and prove that in a neighborhood of the
subset $M^{U(1)}$ of ``regular'' $U(1)$-fixed points every such
manifold $M$ is canonically isomorphic to an open neighborhood of
the zero section in a total space $\overline{T}M^{U(1)}$ of the
tangent bundle to the fixed point set. A fixed point $m \in M$ is
``regular'' if the group $U(1)$ acts on the tangent space $T_mM$
with weights $0$ and $1$. We call this canonical isomorphism {\em
the linearization} of the regular Hodge manifold. 

The linearization construction can be considered as a hyperk\"ahler
analog of the Darboux-Weinstein Theorem in the symplectic
geometry. Apart from the cotangent bundles, it can be applied to the
Hitchin-Simpson moduli space $\M$ of flat connections on a K\"ahler
manifold $X$. The regular points in this case correspond to stable
flat connections such that the underlying holomorphic bundle is also
stable. The linearization construction provides a canonical
embedding of the subspace $\M^{reg} \subset \M$ of regular points
into the total space $T^*\M_0$ of the cotangent bundle to the space
$\M_0$ of stable holomorphic bundles on $X$. The unicity statement
of Theorem~\ref{th.1} guarantees that the hyperk\"ahler metric on
$\M^{reg}$ provided by the Simpson theory is the same as the
canonical metric constructed in Theorem~\ref{th.1}.

\item Starting with Section 5, we deal only with total spaces $\tm$
of the complex-conjugate to tangent bundles to complex manifolds
$M$. In Section 5 we describe Hodge manifolds structures on $\tm$ in
terms of certain connections on the locally trivial fibration $\tm
\to M$. ``Connection'' here is understood as a choice of the
subbundle of horizontal vectors, regardless of the vector bundle
structure on the fibration $\tm \to M$. We establish a
correspondence between Hodge manifold structures on $\tm$ and
connections on $\tm$ over $M$ of certain type, which we call {\em
Hodge connection}.

\item In Section 6 we restrict our attention to the formal
neighborhood of the zero section $M \subset \tm$. We introduce the
appropriate ``formal'' versions of all the definitions and then
establish a correspondence between formal Hodge connections on $\tm$
over $M$ and certain differential operators on the manifold $M$
itself, which we call {\em extended connections}. We also introduce
a certain canonical algebra bundle $\B^\cdot(M,\C)$ on the complex
manifold $M$, which we call {\em the Weil algebra} of the manifold
$M$. Extended connections give rise to natural derivations of the
Weil algebra.

\item Before we can proceed with classification of extended
connections on the manifold $M$ and therefore of regular Hodge
manifolds, we need to derive some linear-algebraic facts on the Weil
algebra $\B^\cdot(M,\C)$. This is the subject of Section 7. We begin
with introducing a certain version of the de Rham complex of a
smooth complex manifold, which we call {\em the total de Rham
complex}. Then we combine it the material of Section 6 to define the
so-called {\em total Weil algebra} of the manifold $M$ and establish
some of its properties. Section 7 is the most technical part of the
whole paper. The reader is advised to skip reading it until needed.

\item Section 8 is the main section of the paper. In this section we
prove, using the technical results of Section 7, that extended
connections on $M$ are in a natural one-to-one correspondence with
K\"ah\-le\-ri\-an connections on the cotangent bundle to $M$
(Theorem~\ref{kal=ext}). This proves the formal part of
Theorem~\ref{th.2}.

\item In Section 9 we deal with polarizations. After some
preliminaries, we use Theorem~\ref{th.2} to deduce the formal part
of Theorem~\ref{th.1} (see Theorem~\ref{metrics}).

\item Finally, in Section~\ref{convergence} we study the convergence
of our formal series and prove Theorem~\ref{th.1} and
Theorem~\ref{th.2} in the real-analytic case.

\item We have also attached a brief Appendix, where we sketch a more
conceptual approach to some of the linear algebra used in the paper,
in particular, to our notion of a weakly Hodge map. This approach
also allows us to describe a simple and conceptual proof of
Proposition~\ref{ac}, the most technical of the facts proved in
Section 7. The Appendix is mostly based on results of Deligne and
Simpson (\cite{De2}, \cite{Simpson}).
\end{itemize}

\noindent
{\bf Acknowledgment.} I would like to thank A. Beilinson,
D. Kazhdan, A. Levin, L. Posicelsky, A. Shen and A. Tyurin for
stimulating discussions. I am especially grateful to my friends
Misha Verbitsky and Tony Pantev for innumerable discussions,
constant interest and help, without which this paper most probably
would not have been written. I would also like to mention here how
much I have benefited from a course on moduli spaces and Hodge
theory given by Carlos Simpson at MIT in the Fall of 1993. On a
different note, I would like to express my dearest gratitude to
Julie Lynch, formerly at International Press in Cambridge, and also
to the George Soros's Foundation and to CRDF for providing me with a
source of income during the preparation of this paper.

\tableofcontents

\section{Preliminary facts from linear algebra}

\subsection{Quaternionic vector spaces}

\punkt Throughout the paper denote by $\h$ the $\R$-algebra of quaternions. 

\defn A {\em quaternionic vector space} is a finite-dimensional left module
over the algebra $\h$. 

Let $V$ be a quaternionic vector space.  Every algebra embedding
$I:\C \to \h$ defines by restriction an action of $\C$ on
$V$. Denote the corresponding complex vector space by $V_I$. 

Fix once and for all an algebra embedding $I:\C \to \h$ and call the
complex structure $V_I$ {\em the preferred complex structure} on
$V$. 

\punkt \label{u.acts.on.h}
Let the group $\C^*$ act on the algebra $\h$ by conjugation by
$I(\C^*)$.  Since $I(\R^*) \subset I(\C^*)$ lies in the center of the
algebra $\h$, this action factors through the map $N:\C^* \to
\C^*/\R^* \cong U(1)$ from $\C^*$ to the one-dimensional unitary group
defined by $N(a) = a^{-1}\overline{a}$. Call this action {\em the
standard action} of $U(1)$ on $\h$.

The standard action commutes with the multiplication and leaves
invariant the subalgebra $I(\C)$. Therefore it extends to an action
of the complex algebraic group $\C^* \supset U(1)$ on the algebra
$\h$ considered as a right complex vector space over $I(\C)$. Call
this action {\em the standard action} of $\C^*$ on $\h$.

\punkt \defn An {\em equivariant} quaternionic vector space is a quaternionic
vector space $V$ equipped with an action of the group $U(1)$ such that the
action map 
$$
\h \otimes_\R V \to V
$$
is $U(1)$-equivariant. 

The $U(1)$-action on $V$ extends to an action of $\C^*$ on the complex
vector space $V_I$. The action map $\h \otimes_\R V \to V$ factors through a
map 
$$
\Mult:\h \otimes_\C V_I \to V_I 
$$
of complex vector spaces. This map is $\C^*$-equivariant if and only if $V$
is an equivariant quaternionic vector space. 

\punkt The category of complex algebraic representations $V$ of the
group $\C^*$ is equivalent to the category of graded vector spaces $V
= \oplus V^\cdot$. We will say that a representation $W$ is {\em of
weight $i$} if $W = W^i$, that is, if an element $z \in \C^*$ acts
on $W$ by multiplication by $z^k$. For every representation $W$ we
will denote by $W(k)$ the representation corresponding to the
grading
$$
W(k)^i = W^{k+i}.
$$

\punkt \label{h=c+c} 
The algebra $\h$ considered as a complex vector space by means of right 
multiplication by $I(\C)$ decomposes $\h = I(\C) \oplus \cc$ with respect
to the standard $\C^*$-action. The first summand is of weight $0$, and the
second is of weight $1$. This decomposition is compatible with the left
$I(\C)$-actions as well and induces for every complex vector space $W$ a
decomposition
$$
\h \otimes_\C W = W \oplus \cc \otimes_\C W.
$$
If $W$ is equipped with an $\C^*$-action, the second summand is canonically
isomorphic to $\overline{W}(1)$, where $\overline{W}$ is the vector space
complex-conjugate to $W$. 

\punkt Let $V$ be an equivariant quaternionic vector space. The action map 
$$
\Mult:\h \otimes_\C V_I \cong V_I \oplus \cc \otimes V_I \to V_I
$$
decomposes $\Mult = \id \oplus j$ for a certain map $j:\overline{V_I}(1)
\to V_I$. The map $j$ satisfies $j \circ \overline{j} = -\id$, and we
obviously have the following. 

\begin{lemma}\label{explicit.eqvs}
The correspondence $V \longmapsto \langle V_I, j \rangle$ is an
equivalence of categories bet\-ween the category of equivariant
quaternionic vector spaces and the category of pairs $\langle W, j
\rangle$ of a graded complex vector space $W$ and a map
$j:\overline{W^\cdot} \to W^{1-\cdot}$ satisfying $j \circ
\overline{j} = -\id$.
\end{lemma}

\punkt We will also need a particular class of equivariant quaternionic
vector spaces which we will call regular. 

\defn An equivariant quaternionic vector space $V$ is called {\em regular}
if every irreducible $\C^*$-subrepresentation of $V_I$ is either trivial or 
of weight $1$. 

\begin{lemma}\label{regular.quaternionic}
Let $V$ be a regular $\C^*$-equivariant quaternionic vector space and let
$V^0_I \subset V_I$ be the subspace of $\C^*$-invariant vectors. Then the
action map
$$
V^0_I \oplus \overline{V^0_I} \cong \h \otimes_\C V^0_I \to V_I
$$
is invertible. 
\end{lemma}

\proof Let $V^1_I \subset V_I$ be the weight $1$ subspace with respect to
the $ gm$-action.  Since $V$ is regular, $V_I = V^0_I \oplus V^1_I$. On the
other hand, $j:\overline{V_I} \to V_I$ interchanges $V^0_I$ and
$V^1_I$. Therefore $V^1_I \cong \overline{V^0_I}$ and we are done.
\endproof

Thus every regular equivariant quaternionic vector space is a sum of
several copies of the algebra $\h$ itself. The corresponding Hodge structure 
has Hodge numbers $h^{1,0} = h^{0,1}$, $h^{p,q} = 0$ otherwise. 

\subsection{The complementary complex structure}\label{complementary} 
 
\punkt 
Let $J:\C \to \h$ be another algebra embedding. Say that embeddings $I$ and
$J$ are {\em complementary} if
$$
J(\sqrt{-1}) I(\sqrt{-1}) = - I(\sqrt{-1}) J(\sqrt{-1}).
$$
Let $V$ be an equivariant quaternionic vector space.  The standard
$U(1)$-action on $\h$ induces an action of $U(1)$ on the set of all algebra
embeddings $\C \to \h$.  On the subset of embeddings complementary to $I$
this action is transitive.  Therefore for every two embeddings $J_1,J_2:\C
\to \h$ complementary to $I$ the complex vector spaces $V_{J_1}$ and
$V_{J_2}$ are canonically isomorphic. We will from now on choose for
convenience an algebra embedding $J:\C \to \h$ complementary to $I$ and
speak of {\em the complementary complex structure} $V_J$ on $V$; however,
nothing depends on this choice. 

\punkt For every equivariant quaternionic vector space $V$ the
complementary embedding $J:\C \to \h$ induces an isomorphism
$$
J \otimes \id:\C \otimes_\R V \to \h \otimes_{I(\C)} V_I 
$$
of complex vector spaces. Let $\Mult:\h \otimes_{I(\C)} V_I \to V_I$, $\Mult:\C
\otimes_\R V \to V_J$ be the action maps. Then there exists a unique
isomorphism $H:V_J \to V_I$ of complex vector spaces such that the diagram 
$$
\begin{CD}
\C \otimes_\R V @>{J \otimes \id}>> \h \otimes_{I(\C)} V_I\\
@V{\Mult}VV                          @V{\Mult}VV          \\
V_J             @>{H}>>              V_I
\end{CD}
$$
is commutative. Call the map $H:V_J \to V_I$ {\em the standard
isomorphism} between the complementary and the preferred complex structures
on the equivariant quaternionic vector space $V$. 

\punkt Note that both $V_I$ and $V_J$ are canonically isomorphic to $V$ as
real vector spaces; therefore the map $H:V_J \to V_I$ is in fact an
automorphism of the real vector space $V$. Up to a constant this
automorphism is given by the action of the element $I(\sqrt{-1}) +
J(\sqrt{-1}) \in \h$ on the $\h$-module $V$. 

\subsection{$\protect\R$-Hodge structures} 

\punkt Recall that {\em a pure $\R$-Hodge structure $W$ of weight
$i$} is a pair of a graded complex vector space $W = \oplus_{p+q=i}
W^{p,q}$ and a {\em real structure} map $\conj:\overline{W^{p,q}}
\to W^{q,p}$ satisfying $\conj \circ \conj = \id$. The bigrading
$W^{p,q}$ is called {\em the Hodge type bigrading}. The dimensions
$h^{p,q} = \dim_\C W^{p,q}$ are called {\em the Hodge numbers} of
the pure $\R$-Hodge structure $W$. Maps between pure Hodge
structures are by definition maps of their underlying complex vector
spaces compatible with the bigrading and the real structure maps.

\punkt Recall also that for every $k$ the {\em Hodge-Tate} pure $\R$-Hodge
structure $\R(k)$ of weight $-2k$ is by definition the
$1$-dimensional complex vector space with complex conjugation equal
to $(-1)^k$ times the usual one, and with Hodge bigrading
$$
\R(k)^{p,q} = \begin{cases} \R(k), \quad p=q=-k,\\ 0, \quad
\text{otherwise}. \end{cases} 
$$
For a pure $\R$-Hodge structure $V$ denote, as usual, by $V(k)$ the
tensor product $V(k) = V \otimes \R(k)$. 

\punkt \label{w.1} 
We will also need special notation for another common $\R$-Hodge
structure, which we now introduce. Note that for every complex $V$
be a complex vector space the complex vector space $V \otimes_\R \C$
carries a canonical $\R$-Hodge structure of weight $1$ with Hodge
bigrading given by
$$
V^{1,0} = V \subset V \otimes_\R \C \qquad\qquad V^{0,1} =
\overline{V} \otimes_\R \C. 
$$
In particular, $\C \otimes_\R \C$ carries a natural $\R$-Hodge
structure of weight $1$ with Hodge numbers $h^{1,0} = h^{0,1} =
1$. Denote this Hodge structure by $\W_1$.  

\punkt Let $\langle W, \conj \rangle$ be a pure Hodge structure, and
denote by $W_\R \subset W$ the $\R$-vector subspace of elements
preserved by $\conj$. Define the {\em Weil operator} $C:W \to W$ by
$$
C = (\sqrt{-1})^{p-q}:W^{p,q} \to W^{q,p}.
$$
The operator $C:W \to W$ preserves the $\R$-Hodge structure, in
particular, the subspace $W_\R \subset W$.  On pure $\R$-Hodge
structures of weight $0$ the Weil operator $C$ corresponds to the
action of $-1 \in U(1) \subset \C^*$ in the corresponding
representation.

\punkt \label{Weil}
For a pure Hodge structure $W$ of weight $i$ let  
$$
j = C \circ \conj:\overline{W^{\cdot,1-\cdot}} \to W^{1-\cdot,\cdot}. 
$$
If $i$ is odd, in particular, if $i=1$, then $j \circ \overline{j} =
-\id$. Together with Lemma~\ref{explicit.eqvs} this gives the following.

\begin{lemma}\label{eqvs.hodge}
The category of equivariant quaternionic vector spaces is equivalent to the
category of pure $\R$-Hodge structures of weight $1$. 
\end{lemma} 

\punkt Let $V$ be an equivariant quaternionic vector space, and let
$\langle W, \conj \rangle$ be $\R$-Hodge structure of weight $1$ 
corresponding to $V$ under the equivalence of
Lemma~\ref{eqvs.hodge}. By definition the complex vector space $W$
is canonically isomorphic to the complex vector space $V_I$ with the
preferred complex structure on $V$. It will be more convenient for
us to identify $W$ with the complementary complex vector space $V_J$
by means of the standard isomorphism $H:V_J \to V_I$. The
multiplication map 
$$
\Mult:V_I \otimes_\R \C \cong V \otimes_\C \h \to V_J
$$
is then a map of $\R$-Hodge structures. The complex conjugation $\conj:W \to
\overline{W}$ is given by
\begin{equation}\label{i.conj}
\conj = C \circ H \circ j \circ H^{-1} = C \circ i:V_J \to \overline{V_J}, 
\end{equation}
where $i:V_J \to \overline{V_J}$ is the action of the element
$I(\sqrt{-1}) \subset \h$. 

\subsection{Weakly Hodge maps}\label{w.H.sub}

\punkt Recall that the category of pure $\R$-Hodge structures of weight $i$
is equivalent to the category of pairs $\langle V, F^\cdot\rangle$ of a
real vector space $V$ and a decreasing filtration $F^\cdot$ on the complex
vector space $V_\C = V \otimes_\R \C$ satisfying 
$$
V_\C = \bigoplus_{p+q=i}F^pV_\C \cap \overline{F^qV_\C}.
$$
The filtration $F^\cdot$ is called {\em the Hodge filtration}. 
The Hodge type grading and the Hodge filtration are related by $V^{p,q} =
F^pV_\C \cap \overline{F^qV_\C}$ and $F^pV_\C = \oplus_{k \geq
p}V^{k,i-k}$. 

\punkt \label{weakly.hodge} Let $\langle V,F^\cdot\rangle$ and
$\langle W,F^\cdot \rangle$ be pure $\R$-Hodge structures of weights
$n$ and $m$ respectively. Usually maps of pure Hodge structures are
required to preserve the weights, so that $\Hom(V,W)=0$ unless
$n=m$. In this paper we will need the following weaker notion of
maps between pure $\R$-Hodge structures.

\defn An $\R$-linear map $f:V \to W$ is said to be {\em weakly
Hodge} if it preserves the Hodge filtrations.

Equivalently, the complexified map $f:V_\C \to W_\C$ must decompose 
\begin{equation}\label{H.type}
f = \sum_{0 \leq p \leq m-n}f^{p,m-n-p},
\end{equation}
where the map $f^{p,m-n-p}:V_\C \to W_\C$ is of Hodge type
$(p,m-n-p)$. Note that this condition is indeed weaker than the
usual definition of a map of Hodge structures: a weakly Hodge map
$f:V \to W$ can be non-trivial if $m$ is strictly greater than
$n$. If $m < n$, then $f$ must be trivial, and if $m=n$, then weakly
Hodge maps from $V$ to $W$ are the same as usual maps of $\R$-Hodge
structures.

\punkt We will denote by $\prehodge$ the category of pure $\R$-Hodge
structures of arbitrary weight with weakly Hodge maps as
morphisms. For every integer $n$ let $\prehodge_n$ be the full
subcategory in $\prehodge$ consisting of pure $\R$-Hodge structures
of weight $n$, and let $\prehodge_{\geq n}$ be the full subcategory
of $\R$-Hodge structures of weight not less than $n$. Since weakly
Hodge maps between $\R$-Hodge structures of the same weight are the
same as usual maps of $\R$-Hodge structures, the category
$\prehodge_n$ is the usual category of pure $\R$-Hodge structures of
weight $n$.

\punkt \label{w.k} Let $\W_1 = \C \otimes_\R \C$ be the pure Hodge
structure of weight $1$ with Hodge numbers $h^{1,0}=h^{0,1}=1$, as
in \ref{w.1}.  The diagonal embedding $\C \to \C \otimes_\R \C$
considered as a map $w_1:\R \to \W_1$ from the trivial pure
$\R$-Hodge structure $\R$ of weight $0$ to $\W_1$ is obviously
weakly Hodge. It decomposes $w_1 = w_1^{1,0} + w_1^{0,1}$ as in
\eqref{H.type}, and the components $w_1^{1,0}:\C \to \W_1^{1,0}$ and
$w_1^{0,1}:\C \to \W_1^{0,1}$ are isomorphisms of complex vectors
spaces. Moreover, for every pure $\R$-Hodge structure $V$ of weight
$n$ the map $w_1$ induces a weakly Hodge map $w_1:V \to \W_1 \otimes
V$, and the components $w_1^{1,0}:V_\C \to V_\C \otimes \W_1^{1,0}$
and $w_1^{0,1}:V_\C \to V_\C \otimes \W_1^{0,1}$ are again
isomorphisms of complex vector spaces.

More generally, for every $k \geq 0$ let $\W_k = S^k\W_1$ be the
$k$-th symmetric power of the Hodge structure $\W_1$.  The space
$\W_k$ is a pure $\R$-Hodge structure of weight $k$, with Hodge
numbers $h^{k,0} = h^{k-1,1} =\ldots = h^{0,k} = 1$. Let $w_k:\R \to
\W_k$ be the $k$-th symmetric power of the map $w_1:\R \to
\W_1$. For every pure $\R$-Hodge structure $V$ of weight $n$ the map
$w_k$ induces a weakly Hodge map $w_k:V \to \W_k \otimes V$, and the
components 
$$
w_k^{p,q}:V_\C \to V_\C \otimes \W_k^{p,k-p}, \qquad 0 \leq p \leq k
$$
are isomorphisms of complex vector spaces. 

\punkt \label{w.k.uni} 
The map $w_k$ is a universal weakly Hodge map from a pure $\R$-Hodge
structures $V$ of weight $n$ to a pure $\R$-Hodge structure of
weight $n+k$. More precisely, every weakly Hodge map $f:V \to V'$
from $V$ to a pure $\R$-Hodge structure $V'$ of weight $n+k$ factors
uniquely through $w_k:V \to \W_k \otimes V$ by means of a map
$f':\W_k \otimes V \to V'$ preserving the pure $\R$-Hodge
structures. Indeed, $V_\C \otimes \W_k = \bigoplus_{0 \leq p \leq k}
V_\C \otimes \W_k^{p,k-p}$, and the maps $w_k^{p,k-p}:V_\C \to V_\C
\otimes \W_k^{p,k-p}$ are invertible.  Hence to obtain the desired
factorization it is necessary and sufficient to set
$$
f' = f^{p,k-p} \circ \left(w_k^{p,k-p}\right)^{-1}:V_\C \otimes
\W_k^{p,k-p} \to V_\C \to \V'_\C, 
$$
where $f = \sum_{0 \leq p \leq k}f^{p,k-p}$ is the Hodge type decomposition
\eqref{H.type}. 

\punkt It will be convenient to reformulate the universal properties
of the maps $w_k$ as follows. By definition $\W_k = S^k\W_1$,
therefore the dual $\R$-Hodge structures equal $\W_k^* = S^k\W_1^*$,
and for every $n,k \geq 0$ we have a canonical projection
$\can:\W_n^* \otimes \W_k^* \to \W_{n+k}^*$. For every pure
$\R$-Hodge structure $V$ of weight $k \geq 0$ let $\Gamma(V) = V
\otimes \W_k^*$.

\begin{lemma}\label{g.ex}
The correspondence $V \mapsto \Gamma(V)$ extends to a functor
$$
\Gamma:\prehodge_{\geq 0} \to \prehodge_0
$$ 
adjoint on the right to the canonical embedding $\prehodge_0
\hookrightarrow \prehodge_{\geq 0}$.
\end{lemma}

\proof Consider a weakly Hodge map $f:V_n
\to V_{n+k}$ from $\R$-Hodge structure $V_n$ of weight $n$ to a pure
$\R$-Hodge structure $V_{n+k}$ of weight $n+k$. By the universal
property the map $f$ factors through the canonical map $w_k:V_n \to
V_n \otimes \W_k$ by means of a map $f_k:V_n \otimes \W_k \to
V_{n+k}$. Let $f_k':V_k \to V_{n+k} \otimes \W_k^*$ be the adjoint
map, and let
\begin{align*}
\Gamma(f) = \can \circ f_k':&\Gamma(V_n) = V_n \otimes \W_k^* \to
V_{n+k} \otimes \W_n^* \otimes \W_k^* \to\\ 
\to &\Gamma(V_{n+k}) = V \otimes \W_{n+k}^*. 
\end{align*}
This defines the desired functor $\Gamma:\prehodge_{\geq 0} \to
\prehodge_0$. The adjointness is obvious. 
\endproof

\rem See Appendix for a more geometric description of the functor
$\Gamma:\prehodge_{\geq 0} \to \prehodge_0$.

\punkt \label{l.r} 
The complex vector space $\W_1 = \C \otimes_\R \C$ is equipped with
a canonical skew-symmetric trace pairing $\W_1 \otimes_\C \W_1 \to
\C$. Let $\gamma:\C \to \W_1^*$ be the map dual to $w_1:\R \to \W_1$
under this pairing. The map $\gamma$ is not weakly Hodge, but it
decomposes $\gamma = \gamma^{-1,0} + \gamma^{0,-1}$ with respect to
the Hodge type grading. Denote $\gamma_l = \gamma^{-1,0}$, $\gamma_r
= \gamma^{0,-1}$. For every $0 \leq p \leq k$ the symmetric powers
of the maps $\gamma_l$, $\gamma_r$ give canonical complex-linear
embeddings
$$
\gamma_l, \gamma_r:\W_p^* \to \W_k^*. 
$$

\punkt The map $\gamma_l$ if of Hodge type $(p-k,0)$, while
$\gamma_r$ is of Hodge type $(0, p-k)$, and the maps $\gamma_l$,
$\gamma_r$ are complex conjugate to each other. Moreover, they are
each compatible with the natural maps $\can:\W_p^* \otimes \W_q^*
\to \W_{p+q}^*$ in the sense that $\can \circ (\gamma_l \otimes
\gamma_l) = \gamma_l \circ \can$ and $\can \circ (\gamma_r \otimes
\gamma_r) = \gamma_r \circ \can$. For every $p$, $q$, $k$ such that
$p+q \geq k$ we have a short exact sequence
\begin{equation}\label{p+q}
\begin{CD}
0 @>>> \W_{p+q-k} @>{\gamma_r - \gamma_l}>> \W_p^* \oplus \W_q^*
@>{\gamma_l + \gamma_l}>> \W_k^* @>>> 0 
\end{CD}
\end{equation}
of complex vector spaces. We will need this exact sequence in
\ref{gamma.use}.

\punkt \label{gamma.tensor} 
The functor $\Gamma$ is, in general, not a tensor functor. However,
the canonical maps $\can:\W_n^* \otimes \W_k^* \to \W_{n+k}^*$
define for every two pure $\R$-Hodge structures $V_1$, $V_2$ of
non-negative weight a surjective map
$$
\Gamma(V_1) \otimes \Gamma(V_2) \to \Gamma(V_1 \otimes V_2). 
$$
These maps are functorial in $V_1$ and $V_2$ and commute with the
associativity and the commutativity morphisms. Moreover, for every
algebra $\A$ in the tensor category $\prehodge_{\geq 0}$ they turn
$\Gamma(\A)$ into an algebra in $\prehodge_0$.

\subsection{Polarizations}

\punkt \label{hyperherm}
Consider a quaternionic vector space $V$, and let $h$ be a Euclidean metric
on $V$.

\defn The metric $h$ is called {\em Qua\-ter\-nionic\--Her\-mi\-tian} if
for any algebra embedding $I:\C \to \h$ the metric $h$ is the real part of
an Hermitian metric on the complex vector space $V_I$.

Equivalently, a metric is Quaternionic-Hermitian if it is invariant under the
action of the group $SU(2) \subset \h$ of unitary quaternions. 

\punkt \label{hermhodge}
Assume that the quaternionic vector space $V$ is equivariant. Say
that a metric $h$ on $V$ is {\em Hermitian-Hodge} if it is
Qua\-ter\-nionic-\-Her\-mi\-tian and, in addition, invariant under the
$U(1)$-action on $V$.

Let $V_I$ be the vector space $V$ with the preferred complex structure $I$,
and let 
$$
V_I = \bigoplus V^\cdot 
$$ 
and $j:\overline{V^\cdot} \to V^{1-\cdot}$ be as in
Lemma~\ref{explicit.eqvs}. The metric $h$ is Hermitian-Hodge if and only if
\begin{enumerate}
\item it is the real part of an Hermitian metric on $V_I$, 
\item $h(V^p,V^q) = 0$ whenever $p \neq q$, and 
\item $h(j(a),b) = - h(a,j(b))$ for every $a,b \in V$. 
\end{enumerate}

\punkt Recall that a {\em polarization} $S$ on a pure $\R$-Hodge
structure $W$ of weight $i$ is a bilinear form $S:W \otimes W \to
\R(-i)$ which is a morphism of pure Hodge structures and satisfies
\begin{align*}
S(a,b) &= (-1)^{i}S(b,a)\\
S(a,Ca) &> 0
\end{align*}
for every $a,b \in W$. (Here $C:W \to W$ is the Weil operator.) 

\punkt \label{pol} 
Let $V$ be an equivariant quaternionic vector space equipped with an
Euclidean metric $h$, and let $\langle W, \conj \rangle$ be the pure
$\R$-Hodge structure of weight $1$ associated to $V$ by
Lemma~\ref{eqvs.hodge}. Recall that $W = V_J$ as a complex vector
space. Assume that $h$ extends to an Hermitian metric $h_J$ on
$V_J$, and let $S:W \otimes W \to \R(-1)$ be the bilinear form
defined by
$$
S(a,b) = h(a,C\overline{b}), \qquad a,b \in W_\R \subset W.
$$
The form $S$ is a polarization if and only if the metric $h$ is
Hermitian-Hodge. This gives a one-to-one correspondence between the set of
Hermitian-Hodge metrics on $V$ and the set of polarizations on the Hodge
structure $W$. 

\punkt Let $W^*$ be the Hodge structure of weight $-1$ dual to $W$.
The sets of polarizations on $W$ and on $W^*$ are, of course, in a
natural one-to-one correspondence. It will be more convenient for us
to identify the set of Hermitian-Hodge metrics on $V$ with the set
of polarizations on $W^*$ rather then on $W$.

Assume that the metric $h$ on the equivariant quaternionic vector
space $V$ is Hermitian-Hodge, and let $S \in \Lambda^2(W) \subset
\Lambda^2(V \otimes \C)$ be the corresponding polarization. Extend
$h$ to an Hermitian metric $h_I$ on the complex vector space $V$
with the preferred complex structure $V_I$, and let
$$
\omega_I \in V_I \otimes \overline{V_I} \subset \Lambda^2(V
\otimes_\R \C) 
$$
be the imaginary part of the corresponding Hermitian metric on the
dual space $V_I^*$. Let $i:V_J \to \overline{V_J}$ the action of the
element $I(\sqrt{-1}) \in \h$. By \eqref{i.conj} we have
$$
\omega_I(a,b) = h(a,i(b)) = h(a,C\overline{b}) = S(a,b) 
$$
for every $a,b \in V_J \subset V \otimes \C$. Since $\omega_I$ is
real, and $V \otimes \C = V_J \oplus \overline{V_J}$, the $2$-forms
$\omega_I$ and $S$ are related by
\begin{equation}\label{omega.and.Omega}
\omega_I = \half(S + \nu(S)), 
\end{equation}
where $\nu:V \otimes \C \to \overline{V \otimes \C}$ is the usual
complex conjugation. 

\section{Hodge bundles and quaternionic manifolds}\label{hbqm.section}

\subsection{Hodge bundles}\label{hb.sub}

\punkt Throughout the rest of the paper, our main tool in studying
hyperk\"ahler structures on smooth manifolds will be the equivalence
between equivariant quaternionic vector spaces and pure $\R$-Hodge
structures of weight $1$, established in Lemma~\ref{explicit.eqvs}.
In order to use it, we will need to generalize this equivalence to
the case of vector bundles over a smooth manifold $M$, rather than
just vector spaces. We will also need to consider manifolds equipped
with a smooth action of the group $U(1)$, and we would like
our generalization to take this $U(1)$-action into account. 

Such a generalization requires, among other things, an appropriate
notion of a vector bundle equipped with a pure $\R$-Hodge structure.
We introduce and study one version of such a notion in this section,
under the name of a Hodge bundle (see
Definition~\ref{hodge.bundles}). 

\punkt Let $M$ be a smooth manifold equipped with a smooth
$U(1)$-action (or {\em a $U(1)$-manifold}, to simplify the
terminology), and let $\iota:M \to M$ be the action of the element
$-1 \in U(1) \subset \C^*$.

\defn \label{hodge.bundles} 
An {\em Hodge bundle of weight $k$} on $M$ is a pair
$\langle\E,\conj\rangle$ of a $U(1)$-equi\-va\-ri\-ant complex
vector bundle $\E$ on $M$ and a $U(1)$-equivariant bundle map
$\conj:\overline{\iota^*\E}(k) \to \E$ satisfying $\conj \circ
\iota^*\conj = \id$.

Hodge bundles of weight $k$ over $M$ form a tensor $\R$-linear
additive category, denoted by $\prehodge_k(M)$. 

\punkt\label{w.hodge}
Let $\E$, $\F$ be two Hodge bundles on $M$ of weights $m$ and
$n$. Say that a bundle map. or, more generally, a differential operator
$f:\E \to \F$ is {\em weakly Hodge} if
\begin{enumerate}
\item $f = \overline{\iota^*f}$, and 
\item there exists a decomposition $f = \sum_{0 \leq n-m} f_i$ with
$f_i$ being of degree $i$ with respect to he $U(1)$-equivariant
structure.  (In particular, $f=0$ unless $n \geq m$, and we always
have $f_k = \overline{\iota^* f_{n-m-k}}$.)
\end{enumerate}

Denote by $\prehodge(M)$ the category of Hodge bundles of arbitrary
weight on $M$, with weakly Hodge bundle maps as morphisms. For every
$i$ the category $\prehodge_i(M)$ is a full subcategory in
$\prehodge(M)$.  Introduce also the category $\prehodge^\D(M)$ with
the same objects as $\prehodge(M)$ but with weakly Hodge
differential operators as morphisms. Both the categories
$\prehodge(M)$ and $\prehodge^\D(M)$ are additive $\R$-linear tensor
categories.

\punkt\label{H-type} 
For a weakly Hodge map $f:\E \to \F$ call the canonical
decomposition
$$
f = \sum_{0 \leq i \leq m-n} f_i 
$$
{\em the $H$-type decomposition}. 
\label{gamma.m}
For a Hodge bundle $\E$ on $M$ of non-negative weight $k$ let
$\Gamma(\E) = \E \otimes \W_k^*$, where $\W_k$ is the canonical pure
$\R$-Hodges structure introduced in \ref{w.k}.  The universal
properties of the Hodge structures $\W_k$ and Lemma~\ref{g.ex}
generalize immediately to Hodge bundles. In particular, $\Gamma$
extends to a functor $\Gamma:\prehodge_{\geq 0}(M) \to
\prehodge_0(M)$ adjoint on the right to the canonical embedding.

\punkt \label{u.trivial}
If the $U(1)$-action on $M$ is trivial, then Hodge bundles of weight
$i$ are the same as real vector bundles $\E$ equipped with a Hodge
type grading $\E = \oplus_{p+q=i} \E^{p,q}$ on the complexification
$\E_\C = \E \otimes_\R \C$. In particular, if $M=\PT$ is a single
point, then $\prehodge(M) \cong \prehodge^\D(M)$ is the category of
pure $\R$-Hodge structures. Weakly Hodge bundle map are then the
same as weakly Hodge maps of $\R$-Hodge structures considered in
\ref{weakly.hodge}. (Thus the notion of a Hodge bundle is indeed a
generalization of the notion of a pure $\R$-Hodge structure.) 

\punkt The categories of Hodge bundles are functorial in $M$, namely, for
every smooth map $f:M_1 \to M_2$ of smooth $U(1)$-manifolds $M_1$, $M_2$
there exist pull-back functors
\begin{align*}
f^*:\prehodge(M_1)    &\to \prehodge(M_1)   \\
f^*:\prehodge^\D(M_1) &\to \prehodge^\D(M_1).
\end{align*} 
In particular, let $M$ be a smooth $U(1)$-manifold and let $\pi:M \to \PT$
be the canonical projection. Then every $\R$-Hodge structure $V$ of weight
$i$ defines a constant Hodge bundle $\pi^*V$ on $M$, which we denote for
simplicity by the same letter $V$. Thus the trivial bundle $\R =
\Lambda^0(M) = \pi^*\R(0)$ has a natural structure of a Hodge bundle of
weight $0$. 

\punkt \label{de.Rham} 
To give a first example of Hodge bundles and weakly Hodge maps,
consider a $U(1)$-manifold $M$ equipped with a $U(1)$-invariant
almost complex structure $M_I$.  Let $\Lambda^i(M,\C) = \oplus
\Lambda^{\cdot,i-\cdot}(M_I)$ be the usual Hodge type decomposition
of the bundles $\Lambda^i(M,\C)$ of complex valued differential
forms on $M$. The complex vector bundles $\Lambda^{p,q}(M_I)$ are
naturally $U(1)$-equivariant. Let 
$$
\conj:\Lambda^{p,q}(M_I) \to \iota^*\overline{\Lambda^{q,p}(M_I)} 
$$
be the usual complex conjugation, and introduce a $U(1)$-equivariant
structure on $\Lambda^\cdot(M,\C)$ by setting
$$
\Lambda^i(M,\C) = \bigoplus_{0 \leq j \leq i} \Lambda^{j,i-j}(M)(j). 
$$
The bundle $\Lambda^i(M,\C)$ with these $U(1)$-equivariant structure
and complex conjugation is a Hodge bundle of weight $i$ on $M$. The
de Rham differential $d_M$ is weakly Hodge, and the $H$-type
decomposition for $d_M$ is in this case the usual Hodge type
decomposition $d = \6 + \bar\6$.

\punkt \rem Definition~\ref{hodge.bundles} is somewhat technical. It can be
heuristically rephrased as follows. For a complex vector bundle $\E$ on
$M$ the space of smooth global section $C^\infty(M,\E)$ is a module over the
algebra $C^\infty(M,\C)$ of smooth $\C$-valued functions on $M$, and the
bundle $\E$ is completely defined by the module $C^\infty(M,\E)$. The
$U(1)$-action on $M$ induces a representation of $U(1)$ on the algebra
$C^\infty(M,\C)$. Let $\nu:C^\infty(M,\C) \to C^\infty(M,\C)$ be
composition of the complex conjugation and the map $\iota^*:C^\infty(M,\C)
\to C^\infty(M,\C)$. The map $\nu$ is an anti-complex involution; together
with the $U(1)$-action it defines a pure $\R$-Hodge structure of weight $0$
on the algebra $C^\infty(M,\C)$. Giving a weight $i$ Hodge bundle structure
on $\E$ is then equivalent to giving a weight $i$ pure $\R$-Hodge structure
on the module $C^\infty(M,\E)$ such that the multiplication map 
$$
C^\infty(M,\C) \otimes C^\infty(M,\E) \to C^\infty(M,\E)
$$
is a map of $\R$-Hodge structures.  

\subsection{Equivariant quaternionic manifolds} 

\punkt We now turn to our main subject, namely, various
dif\-feren\-ti\-al\--ge\-o\-met\-ric structures on smooth manifolds
associated to actions of the algebra $\h$ of quaternions. 

\punkt\defn A smooth manifold $M$ is called {\em quaternionic} if it
is equipped with a smooth action of the algebra $\h$ on its
cotangent bundle $\Lambda^1(M)$. 

Let $M$ be a quaternionic manifold.  Every algebra embedding $J:\C \to \h$
defines by restriction an almost complex structure on the manifold
$M$. Denote this almost complex structure by $M_J$.

\punkt Assume that the manifold $M$ is equipped with a smooth action of the
group $U(1)$, and consider the standard action of $U(1)$ on the vector
space $\h$. Call the quaternionic structure and the $U(1)$-action on $M$
{\em compatible} if the action map 
$$
\h \otimes_\R \Lambda^1(M) \to \Lambda^1(M)
$$
is $U(1)$-equivariant. 

Equivalently, the quaternionic structure and the $U(1)$-action are
compatible if the action preserves the almost complex structure $M_I$, and
the action map 
$$
\h \otimes_\C \Lambda^{1,0}(M_I) \to \Lambda^{1,0}(M_I) 
$$
is $U(1)$-equivariant. 

\punkt \defn A quaternionic manifold $M$ equipped with a compatible smooth
$U(1)$-action is called {\em an equivariant quaternionic manifold}. 

For a $U(1)$-equivariant complex vector bundle $\E$ on $M$ denote by
$\E(k)$ the bundle $\E$ with $U(1)$-equivariant structure twisted by the
$1$-dimensional representation of weight $k$.  Lemma~\ref{explicit.eqvs}
immediately gives the following.

\begin{lemma}\label{explicit.qm}
The category of quaternionic manifolds is equivalent to the category of
pairs $\langle M_I, j \rangle$ of an almost complex manifold $M_I$ and a
$\C$-linear $U(1)$-equivariant smooth map $j:\Lambda^{0,1}(M_I)(1) \to
\Lambda^{0,1}(M_I)$ satisfying $j \circ \overline{j} = -\id$. 
\end{lemma}

\subsection{Quaternionic manifolds and Hodge bundles}

\punkt Let $M$ be a smooth $U(1)$-manifold.  Recall that we have
introduced in Subsection~\ref{hb.sub} a notion of a Hodge bundle on
$M$. Hodge bundles arise naturally in the study of quaternionic
structures on $M$ in the following way. Define a {\em quaternionic
bundle} on $M$ as a real vector bundle $\E$ equipped with a left
action of the algebra $\h$, and let $\Bun(M,\h)$ be the category of
smooth quaternionic vector bundles on the manifold $M$. Let also
$\Bun^{U(1)}(M,\h)$ be the category of smooth quaternionic bundles
$\E$ on $M$ equipped with a $U(1)$-equivariant structure such that
the $\h$-action map $\h \to \Eend(\E)$ is
$U(1)$-equivariant. Lemma~\ref{eqvs.hodge} immediately generalizes
to give the following.

\begin{lemma}\label{eqb.hodge}
The category $\Bun^{U(1)}(M,\h)$ is equivalent to the category of Hodge
bundles of weight $1$ on $M$. 
\end{lemma}

\punkt Note that if the $U(1)$-manifold $M$ is equipped with an almost complex
structure, then the decomposition $\h = \overline{\C} \oplus I(\C)$ (see
\ref{h=c+c}) induces an isomorphism 
$$
\can:\h \otimes_{I(\C)} \Lambda^{0,1} \cong \Lambda^{1,0}(M) \oplus
\Lambda^{0,1}(M) \cong \Lambda^1(M,\C).
$$
The weight $1$ Hodge bundle structure on $\Lambda^1(M,\C)$ corresponding to
the natural quaternionic structure on $\h \otimes_{I(\C)} \Lambda^{0,1}(M)$
is the same as the one considered in \ref{de.Rham}.

\punkt Assume now that the smooth $U(1)$-manifold $M$ is equipped
with a compatible quaternionic structure, and let $M_I$ be the
preferred almost complex structure on $M$.  Since $M_I$ is preserved
by the $U(1)$-action on $M$, the complex vector bundle
$\Lambda^{0,1}(M_I)$ of $(0,1)$-forms on $M_I$ is naturally
$U(1)$-equivariant.

The quaternionic structure on $\Lambda^1(M)$ induces by
Lemma~\ref{eqb.hodge} a weight-$1$ Hodge bundle structure on
$\Lambda^{0,1}(M_I)$. The corresponding $U(1)$-action on
$\Lambda^{0,1}(M_I)$ is induced by the action on $M_I$, and the real
structure map
$$
\conj:\Lambda^{1,0}(M_I)(1) \to \Lambda^{0,1}(M_I)
$$
is given by $\conj = \sqrt{-1} \left( \iota^* \circ j \right)$. (Here $j$
is induced by quaternionic structure, as in Lemma~\ref{explicit.qm}). 

\punkt \label{lambda_j=lambda_i}
Let $M_J$ be the complementary almost complex structure on the
equi\-va\-ri\-ant quaternionic manifold $M$. Recall that in
\ref{complementary} we have defined for every equivariant
quaternionic vector space $V$ the standard isomorphism $H:V_J \to
V_I$. This construction can be immediately generalized to give a
complex bundle isomorphism
$$
H:\Lambda^{0,1}(M_J) \to \Lambda^{0,1}(M_I). 
$$
Let $P:\Lambda^1(M,\C) \to \Lambda^{0,1}(M_J)$ be the natural projection,
and let $\Mult:\h \otimes_{I(\C)} \Lambda^{0,1}(M_I) \to
\Lambda^{0,1}(M_I)$ be the action map. By definition the diagram 
$$
\begin{CD}
\Lambda^1(M,\C)    @>{\can}>> \h \otimes_{I(\C)} \Lambda^{0,1}(M_I) \\
 @V{P}VV                       @V{\Mult}VV                          \\
\Lambda^{0,1}(M_J) @>{H}>>    \Lambda^{0,1}(M_I)
\end{CD}
$$
is commutative. Since the map $\Mult$ is compatible with the Hodge bundle
structures, so is the projection $P$. 

\rem This may seems paradoxical, since the complex conjugation
$\conj$ on $\Lambda^1(M,\C)$ does not preserve $\Ker P =
\Lambda^{0,1}(M_J)$. However, under our definition of a Hodge bundle
the conjugation on $\Lambda^1(M,\C)$ is $\iota^* \circ \conj$ rather
than $\conj$. Both $\conj$ and $\iota^*$ interchange
$\Lambda^{1,0}(M_J)$ and $\Lambda^{0,1}(M_J)$. 

\punkt The standard isomorphism $H:\Lambda^{0,1}(M_J) \to
\Lambda^{0,1}(M_I)$ does not commute with the Dolbeault differentials. They
are, however, related by means of the Hodge bundle structure on
$\Lambda^{0,1}(M_I)$. Namely, we have the following. 

\begin{lemma}\label{d=d_0+d_0}
The Dolbeault differential $D:\Lambda^0(M,\C) \to \Lambda^{0,1}(M_I)$ for
the almost complex structure $M_J$ is weakly Hodge. The $U(1)$-invariant
component $D_0$ in the $H$-type decomposition $D = D_0 + \overline{D_0}$ of
the map $D$ coincides with the Dolbeault differential for the almost complex
structure $M_I$. 
\end{lemma}

\proof The differential $D$ is the composition $D = P \circ d_M$ of the de
Rham differential $d_M:\Lambda^0(M,\C) \to \Lambda^1(M,\C)$ with the
canonical projection $P$. Since both $P$ and $d_M$ are weakly Hodge, so is
$D$. The rest follows from the construction of the standard isomorphism
$H$. 
\endproof 

\subsection{Holonomic derivations}

\punkt Let $M$ be a smooth $U(1)$-manifold. In order to give a
Hodge-theoretic description of the set of all equivariant quaternionic
structures on $M$, it is convenient to work not with various complex
structures on $M$, but with associated Dolbeault differentials. To do this,
recall the following universal property of the cotangent bundle
$\Lambda^1(M)$. 

\begin{lemma}\label{universal}
Let $M$ be a smooth manifold, and let $\E$ be a complex vector bundle on
$M$. Every differential operator $\6:\Lambda^0(M,\C) \to \E$ which is a
derivation with respect to the algebra structure on $\Lambda^0(M,\C)$
factors uniquely through the de Rham differential $d_M:\Lambda^0(M,\C) \to
\Lambda^1(M,\C)$ by means of a bundle map $P:\Lambda^1(M,\C) \to \E$. 
\end{lemma}

\punkt \label{holonomic} 
We first use this universal property to describe almost complex structures.
Let $M$ be a smooth manifold equipped with a complex vector bundle $\E$.

\defn A derivation $D:\Lambda^0(M,\C) \to \E$ is called {\em holonomic} if
the associated bundle map $P:\Lambda^1(M,\C) \to \E$ induces an isomorphism
of the subbundle $\Lambda^1(M,\R) \subset \Lambda^1(M,\C)$ of real
$1$-forms with the real vector bundle underlying $\E$. 

By Lemma~\ref{universal} the correspondence 
$$
M_I \mapsto \left\langle \Lambda^{0,1}(M_I), \bar\6 \right\rangle
$$
identifies the set of all almost complex structures $M_I$ on $M$ 
with the set of all pairs $\langle \E, D \rangle$ of
a complex vector bundle $\E$ and a holonomic derivation $D:\Lambda^0(M,\C)
\to \E$. 

\punkt Assume now that the smooth manifold $M$ is equipped with smooth
action of the group $U(1)$. Then we have the following. 

\begin{lemma}\label{qm.hodge}
Let $\E$ be a weight $1$ Hodge bundle on the smooth $U(1)$-manifold
$M$, and let
$$
D:\Lambda^0(M,\C) \to \E
$$ 
be a weakly Hodge holonomic derivation. There exists a unique
$U(1)$\--equi\-va\-ri\-ant quaternionic structure on $M$ such that
$\E \cong \Lambda^{0,1}(M_J)$ and $D$ is the Dolbeault differential
for the complementary almost complex structures $M_J$ on $M$.
\end{lemma}

\proof Since the derivation $M$ is holonomic, it induces an almost
complex structure $M_J$ on $M$. To construct an almost complex
structure $M_I$ complementary to $M_J$, consider the $H$-type
decomposition $D = D_0 + \overline{D_0}$ of the derivation
$D:\Lambda^0(M,\C) \to \E$ (defined in \ref{H-type}).

The map $D_0$ is also a derivation. Moreover, it is holonomic. Indeed, by
dimension count it is enough to prove that the associated bundle map
$P:\Lambda^1(M,\R) \to \E$ is injective. Since the bundle $\Lambda^1(M,\R)$
is generated by exact $1$-forms, it is enough to prove that any real valued
function $f$ on $M$ with $D_0f = 0$ is constant. However, since $D$ is
weakly Hodge,
$$
Df = D_0 f + \overline{D_0} f = D_0 f + \overline{ D_0 \overline{f} } = D_0
f + \overline{ D_0 f}, 
$$
hence $D_0f = 0$ if and only if $Df=0$. Since $D$ is holonomic, $f$ is indeed
constant. 

The derivation $D_0$, being holonomic, is the Dolbeault differential
for an almost complex structure $M_I$ on $M$.  Since $D_0$ is by
definition $U(1)$-equivariant, the almost complex structure $M_I$ is
$U(1)$-invariant. Moreover, $\E \cong \Lambda^{0,1}(M_I)$ as
$U(1)$-equivariant complex vector bundles.  By Lemma~\ref{eqb.hodge}
the weight $1$ Hodge bundle structure on $\E$ induces an equivariant
quaternionic bundle structure on $\E$ and, in turn, a structure of an
equivariant quaternionic manifold on $M$. The almost complex structure
$M_I$ coincides by definition with the preferred almost complex
structure.

It remains to notice that by Lemma~\ref{d=d_0+d_0} the Dolbeault differential 
$\bar\6_J$ for the complementary almost complex structure on $M$ indeed equals 
$D = D_0 + \overline{D_0}$. 
\endproof 

Together with Lemma~\ref{d=d_0+d_0}, this shows that the set of
equivariant quaternionic structures on the $U(1)$-manifold $M$ is
naturally bijective to the set of pairs $\langle \E, D \rangle$ of a
weight $1$ Hodge bundle $\E$ on $M$ and a weakly Hodge holonomic
derivation $D:\Lambda^0(M,\C) \to \E$.

\section{Hodge manifolds}

\subsection{Integrability}

\punkt There exists a notion of integrability for quaternionic manifolds
analogous to that for the almost complex ones.  

\defn A quaternionic manifold $M$ is called {\em hypercomplex} if for two
complementary algebra embeddings $I,J:\C \to \h$ the almost complex
structures $M_I,M_J$ are integrable. 

In fact, if $M$ is hypercomplex, then $M_I$ is integrable for any algebra
embedding $I:\C \to \h$. For a proof see, e.g., \cite{K}.

\punkt When a quaternionic manifold $M$ is equipped with a compatible
$U(1)$-action, there exist a natural choice for a pair of almost complex
structures on $M$, namely, the preferred and the complementary one.

\defn An equivariant quaternionic manifold $M$ is called a {\em Hodge
manifold} if both the preferred and the equivariant almost complex
structures $M_I$, $M_J$ are integrable. 

Hodge manifolds are the main object of study in this paper. 

\punkt \label{standard} 
There exists a simple Hodge-theoretic description of
Hodge manifolds based on Lemma~\ref{qm.hodge}. To give it (see
Proposition~\ref{explicit.hodge}), consider an equivariant quaternionic
manifold $M$, and let $M_J$ and $M_I$ be the complementary and the
preferred complex structures on $M$. The weight $1$ Hodge bundle structure
on $\Lambda^{0,1}(M_J)$ induces a weight $i$ Hodge bundle structure on the
bundle $\Lambda^{0,i}(M_J)$ of $(0,i)$-forms on $M_J$. The standard
identification $H:\Lambda^{0,1}(M_J) \to \Lambda^{0,1}(M_I)$ given in
\ref{lambda_j=lambda_i} extends uniquely to an algebra isomorphism
$H:\Lambda^{0,i}(M_J) \to \Lambda^{0,i}(M_I)$. 

Let $D:\Lambda^{0\cdot}(M_J) \to \Lambda^{0,\cdot+1}(M_J)$ be the Dolbeault
differential for the almost complex manifold $M_J$.

\begin{lemma}\label{yet.another.lemma} 
The equivariant quaternionic manifold $M$ is Hodge if and only if
the following holds. 
\begin{enumerate}
\item $M_J$ is integrable, that is, $D \circ D = 0$, and 
\item the differential $D:\Lambda^{0,i}(M_J) \to \Lambda^{0,i+1}(M_J)$
is weakly Hodge for every $i \geq 0$. 
\end{enumerate}
\end{lemma}

\proof Assume first that the conditions \thetag{i}, \thetag{ii}
hold. Condition \thetag{i} means that the complementary almost complex
structure $M_J$ is integrable. By \thetag{ii} the map $D$ is weakly Hodge.

Let $D = D_0 + \overline{D_0}$ be the $H$-type decomposition. The
map $D_0$ is an algebra derivation of
$\Lambda^{0,\cdot}(M_I)$. Moreover, by Lemma~\ref{d=d_0+d_0} the map
$D_0:\Lambda^0(M,\C) \to \Lambda^{0,1}(M_J)$ is the Dolbeault
differential $\bar\6_I$ for the preferred almost complex structure
$M_I$ on $M$. (Or, more precisely, is identified with $\bar\6_I$
under the standard isomorphism $H$.) But the Dolbeault differential
admits at most one extension to a derivation of the algebra
$\Lambda^{0,\cdot}(M_J)$. Therefore $D_0 = \bar\6_I$ everywhere.

The composition $D_0 \circ D_0$ is the $(2,0)$-component in the
$H$-type decomposition of the map $D \circ D$. Since $D \circ D = 0$,
$$
D_0 \circ D_0 = \bar\6_I \circ \bar\6_I = 0.
$$
Therefore the preferred complex structure $M_I$ is also integrable, and the
manifold $M$ is indeed Hodge. 

Assume now that $M$ is Hodge. The canonical projection $P:\Lambda^1(M,\C)
\to \Lambda^{0,1}(M_J)$ extends then to a multiplicative projection 
$$
P:\Lambda^\cdot(M,\C) \to \Lambda^{0,\cdot}(M_J)
$$
from the de Rham complex of the complex manifold $M_I$ to the Dolbeault
complex of the complex manifold $M_J$. 
The map $P$ is surjective and weakly Hodge, moreover, it 
commutes with the differentials. Since the de Rham
differential preserves the pre-Hodge structures, so does the Dolbeault
differential $D$. 
\endproof  

\punkt Lemma~\ref{yet.another.lemma} and Lemma~\ref{qm.hodge} together
immediately give the following. 

\begin{prop}\label{explicit.hodge} 
The category of Hodge manifolds is equivalent to the category of triples
$\langle M, \E, D \rangle$ of a smooth $U(1)$-manifold $M$, a weight $1$
Hodge bundle $\E$ on $M$, and a weakly Hodge algebra derivation 
$$
D = D^\cdot:\Lambda^\cdot(\E) \to \Lambda^{\cdot+1}(\E)
$$
such that $D \circ D = 0$, and the first component 
$$
D^0:\Lambda^0(M,\C) = \Lambda^0(\E) \to \E = \Lambda^1(\E)
$$
is holonomic in the sense of \ref{holonomic}. 
\end{prop}

\subsection{The de Rham complex of a Hodge manifold}

\punkt Let $M$ be a Hodge manifold. In this subsection we study in some
detail the de Rham complex $\Lambda^\cdot(M,\C)$ of the manifold $M$, in
order to obtain information necessary for the discussion of metrics on $M$
given in the Subsection~\ref{polarization}. The reader is advised to skip
this subsection until needed.

\punkt Let $\Lambda^{0,\cdot}(M_J)$ be the Dolbeault complex for the
complementary complex structure $M_J$ on $M$. By
Proposition~\ref{explicit.hodge} the complex vector bundle
$\Lambda^{0,i}(M_J)$ is a Hodge bundle of weight $i$ on $M$, and the
Dolbeault differential $D:\Lambda^{0,\cdot}(M_J) \to
\Lambda^{0,\cdot+1}(M_J)$ is weakly Hodge. Therefore $D$ admits an $H$-type
decomposition $D = D_0 + \overline{D_0}$. 

\punkt Consider the de Rham complex $\Lambda^i(M,\C)$ of the smooth
manifold $M$.  Let $\Lambda^i(M,\C) = \oplus_{p+q}
\Lambda^{p,q}(M_J)$ be the Hodge type decomposition for the
complementary complex structure $M_J$ on $M$, and let
$\nu:\Lambda^{p,q}(M_J) \to \overline{\Lambda^{q,p}(M_J)}$ be the
usual complex conjugation. Denote also
$$
f^\nu = \nu \circ f \circ \nu^{-1}
$$
for any map $f:\Lambda^\cdot(M,\C) \to \Lambda^\cdot(M,\C)$. 

Let $d_M:\Lambda^\cdot(M,\C) \to \Lambda^{\cdot+1}(M,\C)$ be the de
Rham differential, and let $d_M = D + D^\nu$ be the Hodge type
decomposition for the complementary complex structure $M_J$ on $M$. 
Since the Dolbeault differential, in turn, equals $D = D_0 +
\overline{D_0}$, we have 
$$
d_M = D_0 + \overline{D_0} + D_0^\nu + \overline{D_0}^\nu.
$$

\punkt \label{6_I}
Let $\bar\6_I:\Lambda^\cdot(M,\C) \to \Lambda^{\cdot+1}(M,\C)$ be the
Dolbeault differential for the preferred complex structure $M_I$ on $M$. As
shown in the proof of Lemma~\ref{yet.another.lemma}, the $(0,1)$-component
of the differential $\bar\6_I$ with respect to the complex structure $M_J$
equals $D_0$. Therefore the $(1,0)$-component of the complex-conjugate map
$\6_I = \bar\6_I^\nu$ equals $D_0^\nu$. Since $d_M = \bar\6_I + \6_I$, we
have
\begin{align*}
\bar\6_I &= D_0 + \overline{D_0}^\nu\\
\6_I &= \overline{D_0} + D_0^\nu
\end{align*}
 
\punkt \label{6_I^H}
The standard isomorphism $H:\Lambda^{0,1}(M_J) \to \Lambda^{0,1}(M_I)$
introduced in \ref{standard} extends uniquely to a bigraded algebra
isomorphism $H:\Lambda^{\cdot,\cdot}(M_J) \to
\Lambda^{\cdot,\cdot}(M_I)$. By definition of the map $H$, on
$\Lambda^0(M,\C)$ we have
\begin{align}\label{identities}
\begin{split}
\bar\6_I  &= H \circ D_0 \circ H^{-1}\\
\6_I      &= H \circ D_0^\nu \circ H^{-1}\\
d_M       &= \6_I + \bar\6_I = H \circ (D_0 + D_0^\nu) H^{-1}. 
\end{split}
\end{align}
The right hand side of the last identity is the algebra derivation
of the de Rham complex $\Lambda^\cdot(M,\C)$. Therefore, by
Lemma~\ref{universal} it holds not only on $\Lambda^0(M,\C)$, but on
the whole $\Lambda^\cdot(M,\C)$. The Hodge type decomposition for
the preferred complex structure $M_I$ then shows that the first two
identities also hold on the whole de Rham complex
$\Lambda^\cdot(M,\C)$.

\punkt Let now $\xi = I(\sqrt{-1}):\Lambda^{0,1}(M_J) \to
\Lambda^{1,0}(M_J)$ be the operator corresponding to the preferred
almost complex structure $M_I$ on $M$. Let also $\xi = 0$ on
$\Lambda^0(M,\C)$ and $\Lambda^{1,0}(M_J)$, and extend $\xi$ to a
derivation $\xi:\Lambda^{\cdot,\cdot}(M_J) \to
\Lambda^{\cdot-1,\cdot+1}(M_J)$ by the Leibnitz rule. We finish this
subsection with the following simple fact.

\begin{lemma}\label{xi.lemma}
On $\Lambda^{\cdot,0}(M_J) \subset \Lambda^\cdot(M,\C)$ we have
\begin{align}\label{xi.eq}
\begin{split}
\xi \circ D_0 + D_0 \circ \xi &= \overline{D}_0^\nu\\
\xi \circ \overline{D}_0 + \overline{D}_0 \circ \xi &= -D_0^\nu. 
\end{split}
\end{align}
\end{lemma}

\proof It is easy to check that both identities hold on
$\Lambda^0(M,\C)$. But both sides of these identities are algebra
derivations of $\Lambda^{\cdot,0}(M_J)$, and the right hand sides
are holonomic in the sense of \ref{holonomic}. Therefore by
Lemma~\ref{universal} both identities hold on the whole
$\Lambda^{\cdot,0}(M_J)$.  
\endproof

\subsection{Polarized Hodge manifolds}\label{polarization}

\punkt Let $M$ be a quaternionic manifold. A Riemannian metric $h$ on $M$
is called {\em Qua\-ter\-ni\-onic-\-Her\-mi\-ti\-an} if for every point $m
\in M$ the induces metric $h_m$ on the tangent bundle $T_mM$ is
Qua\-ter\-ni\-onic\--Her\-mi\-ti\-an in the sense of \ref{hyperherm}. 

\defn A {\em hyperk\"ahler manifold} is a hypercomplex manifold $M$
equipped with a Quaternionic-Hermitian metric $h$ which is K\"ahler for
both integrable almost complex structures $M_I$, $M_J$ on $M$. 

\rem In the usual definition (see, e.g., \cite{Bes}) the integrability of
the almost complex structures $M_I$, $M_J$ is omitted, since it is
automatically implied by the K\"ahler condition.

\punkt Let $M$ be a Hodge manifold equipped with a Riemannian metric $h$.
The metric $h$ is called {\em Hermitian-Hodge} if it is
Quaternionic-Hermitian and, in addition, invariant under the $U(1)$-action
on $M$. 

\defn Say that the manifold $M$ is {\em polarized} by the Hermitian-Hodge
metric $h$ if $h$ is not only Quaternionic-Hermitian, but also hyperk\"ahler.

\punkt \label{positive} Let $M$ be a Hodge manifold. By
Proposition~\ref{explicit.hodge} the holomorphic cotangent bundle
$\Lambda^{1,0}(M_J)$ for the complementary complex structure $M_J$ on $M$
is a Hodge bundle of weight $1$. The holomorphic tangent bundle
$\Theta(M_J)$ is therefore a Hodge bundle of weight $-1$.  By \ref{pol} the
set of all Hermitian-Hodge metrics $h$ on $M$ is in natural one-to-one
correspondence with the set of all polarizations on the Hodge bundle
$\Theta(M_J)$. 

Since $\theta(M)$ is of odd weight, its polarizations are
skew-symmetric as bilinear forms and correspond therefore to smooth 
$(2,0)$-forms on the complex manifold $M_J$. A $(2,0)$-form $\Omega$
defines a polarization on $\Theta(M_J)$ if and only if 
$\Omega \in C^\infty(M,\Lambda^{2,0}(M_J))$
considered as a map
$$
\Omega:\R(-1) \to \Lambda^{2,0}(M_J)
$$
is a map of weight $2$ Hodge bundles, and for an arbitrary smooth
section $\chi \in C^\infty(M,\Theta(M_J))$ we have
\begin{equation}\label{P}
\Omega(\chi,\overline{\iota^*(\chi)}) > 0.
\end{equation}  

\punkt Assume that the Hodge manifold $M$ is equipped with an
Hermitian-Hodge metric $h$. Let $\Omega \in
C^\infty(M,\Lambda^{2,0}(M_J))$ be the corresponding polarization,
and let $\omega_I \in C^\infty(M,\Lambda^{1,1}(M_I))$ be the
$(1,1)$-form on the complex manifold $M_I$ associated to the
Hermitian metric $h$.  Either one of the forms $\Omega$, $\omega_I$
completely defines the metric $h$, and by \eqref{omega.and.Omega} we
have
$$
\Omega + \nu(\Omega) = \omega_I, 
$$
where $\nu:\Lambda^\cdot(M,\C) \to \Lambda^\cdot(M,\C)$ is the complex
conjugation. 

\begin{lemma}\label{pol.hm}
The Hermitian-Hodge metric $h$ polarizes $M$ if and only if the
corresponding $(2,0)$-form $\Omega$ on $M_J$ is holomorphic, that is,
$$
D\Omega = 0,
$$
where $D$ is the Dolbeault differential for complementary complex structure
$M_J$.
\end{lemma}

\proof Let $\omega_I,\omega_J \in \Lambda^2(M,\C)$ be the K\"ahler forms
for the metric $h$ and complex structures $M_I$, $M_J$ on $M$. The metric
$h$ is hyperk\"ahler, hence polarizes $M$, if and only if $d_M\omega_I =
d_M\omega_J = 0$. 

Let $D = D_0 + \overline{D_0}$ be the $H$-type decomposition and 
let $H:\Lambda^{\cdot,\cdot}(M_J) \to \Lambda^{\cdot,\cdot}(M_J)$ be the
standard algebra identification introduced in \ref{6_I^H}. 
By definition $H(\omega_I) = \omega_J$. Moreover, by
\eqref{identities} $H^{-1} \circ d_M \circ H = D_0 +
\overline{D_0}^\nu$, hence  
$$
H(d_M\omega_J) = D_0 \omega_I + \overline{D_0}^\nu\omega_I,
$$ 
and the metric $h$ is hyperk\"ahler if and only if 
\begin{equation}\label{K}
d_M\omega_I = (D_0+\overline{D_0}^\nu)\omega_I = 0
\end{equation}
But $2\omega_I = \Omega + \nu(\Omega)$. Since $\Omega$ is of Hodge
type $(2,0)$ with respect to the complementary complex structure
$M_J$, \eqref{K} is equivalent to
$$
\overline{D_0} \Omega = D_0 \Omega = \overline{D_0}^\nu \Omega =
D_0^\nu \Omega = 0.
$$
Moreover, let $\xi$ be as in Lemma~\ref{xi.lemma}. Then $\xi(\Omega)
= 0$, and by \eqref{xi.eq} $D_0\Omega = \overline{D}_0\Omega = 0$
implies that $D_0^\nu\Omega = \overline{D}_0^\nu\Omega = 0$ as well.
It remains to notice that since the metric $h$ is Hermitian-Hodge,
$\Omega$ is of $H$-type $(1,1)$ as a section of the weight $2$ Hodge
bundle $\Lambda^{2,0}(M_J)$.  Therefore $D\Omega = 0$ implies
$\overline{D_0}\Omega = D_0\Omega = 0$, and this proves the lemma.
\endproof

\rem This statement is wrong for general hyperk\"ahler manifolds
(eve\-ry\-thing in the given proof carries over, except for the
implication $D\Omega=0 \Rightarrow
D_0\Omega=\overline{D_0}\Omega=\overline{D_0}^\nu \Omega = D_0^\nu
\Omega = 0$, which depends substantially on the $U(1)$-action). To
describe general hyperk\"ahler metrics in holomorphic terms, one has
to introduce the so-called {\em twistor spaces} (see, e.g.,
\cite{HKLR}).

\section{Regular Hodge manifolds} 

\subsection{Regular stable points} 

\punkt Let $M$ be a smooth manifold equipped with a smooth
$U(1)$-action with differential $\phi_M$ (thus $\phi_M$ is a smooth
vector field on $M$).  Since the group $U(1)$ is compact, the subset
$M^{U(1)} \subset M$ of points fixed under $U(1)$ is a smooth
submanifold.

Let $m \in M^{U(1)} \subset M$ be a point fixed under $U(1)$. Consider the
representation of $U(1)$ on the tangent space $T_m$ to $M$ at $m$.  Call
the fixed point $m$ {\em regular} if every irreducible subrepresentation of
$T_m$ is either trivial or isomorphic to the representation on $\C$ given
by embedding $U(1) \subset \C^*$. (Here $\C$ is considered as a
$2$-dimensional real vector space.)  Regular fixed points form a union of
connected component of the smooth submanifold $M^{U(1)} \subset M$.

\punkt Assume that $M$ is equipped with a complex structure preserved by
the $U(1)$-action.  Call a point $m \in M$ {\em stable} if for any $t \in
\R, t \geq 0$ there exists $\exp(\sqrt{-1}t\phi_M)m$, and the limit 
$$
m_0 \in M, m_0 = \lim_{t \to +\infty} \exp(\sqrt{-1}t\phi_M)m
$$
also exists. 

\enlargethispage{10mm}

\punkt For every stable point $m \in M$ the limit $m_0$ is obviously fixed
under $U(1)$.  Call a point $m \in M$ {\em regular stable} if it is stable
and the limit $m_0 \in M^{U(1)}$ is a regular fixed point.

Denote by $M^{reg} \subset M$ the subset of all regular stable points. The
subset $M^{reg}$ is open in $M$. 

\ex Let $Y$ be a complex manifold with a holomorphic bundle $\E$ and let
$E$ be the total space of $\E$. Let $\C^*$ act on $M$ by dilatation along the
fibers. Then every point $e \in E$ is regular stable.

\punkt Let $M$ be a Hodge manifold. Recall that the $U(1)$-action on $M$
preserves the preferred complex structure $M_I$.

\defn A Hodge manifold $M$ is called {\em regular} if $M_I^{reg} = M_I$.

\subsection{Linearization of regular Hodge manifolds} 

\punkt Consider a regular Hodge manifold $M$. Let $\Delta \subset
\C$ be the unit disk equipped with the multiplicative semigroup
structure. The group $U(1) \subset \Delta$ is embedded into $\Delta$
as the boundary circle.

\begin{lemma}\label{regular}
The action $a:U(1) \times M \to M$ extends uniquely to a holomorphic
action $\tilde{a}:\Delta \times M_I \to M_I$. Moreover, for every $x
\in \Delta \setminus \{0\}$ the action map $\wt{a}(x):M_I \to M_I$
is an open embedding.
\end{lemma}

\proof Since $M$ is regular, the exponential flow $\exp(it\phi_M)$
of the differential $\phi_M$ of the action is defined for all
positive $t \in \R$. Therefore $a:U(1) \times M \to M$ extends
uniquely to a holomorphic action
$$
\tilde{a}:\Delta^* \times M_I \to M_I, 
$$
where $\Delta^* = \Delta \backslash \{0\}$ is the punctured
disk. Moreover, the exponential flow converges as $t \to +\infty$,
therefore $\tilde{a}$ extends to $\Delta \times M_I$
continuously. Since this extension is holomorphic on a dense open
subset, it is holomorphic everywhere. This proves the first claim.

To prove the second claim, consider the subset $\wt{\Delta} \subset
\Delta^*$ of points $x \in \Delta$ such that $\wt{a}(x)$ is
injective and \'etale. The subset $\wt{\Delta}$ is closed under
multiplication and contains the unit circle $U(1) \subset
\Delta^*$. Therefore to prove that $\wt{\Delta} = \Delta^*$, it
suffices to prove that $\wt{\Delta}$ contains the interval $]0,1]
\subset \Delta^*$.

By definition we have $\wt{a}(h) = \exp(-\sqrt{-1}\log h \phi_M)$
for every $h \in ]0,1] \subset \Delta^*$. Thus we have to prove that
if for some $t \in \R, t \geq 0$ and for two points $m_1,m_2 \in M$
we have
$$
\exp(\sqrt{-1}t\phi_M)(m_1) = \exp(\sqrt{-1}t\phi_M)(m_2), 
$$
then $m_1 = m_2$. Let $m_1$, $m_2$ be such two points and let 
$$
t = \inf\{t \in \R, t \geq 0, \exp(\sqrt{-1}t\phi_M)(m_1) =
\exp(\sqrt{-1}t\phi_M)(m_2)\}. 
$$
If the point $m_0 = \exp(\sqrt{-1}t\phi_M)(m_1) =
\exp(\sqrt{-1}t\phi_M)(m_2) \in M$ is not $U(1)$-invariant, then it
is a regular point for the vector field $\sqrt{-1}\phi_M$, and by
the theory of ordinary differential equations we have $t = 0$ and
$m_1 = m_2 = m_0$.

Assume therefore that $m_0 \in M^{U(1)}$ is $U(1)$-invariant. Since
the group $U(1)$ is compact, the vector field $\sqrt{-1}\phi_M$ has
only a simple zero at $m_0 \subset M^{U(1)} \subset M$. Therefore
$m_0 = \exp(\sqrt{-1}t\phi_M)m_1$ implies that the point $m_1 \in M$
also is $U(1)$-invariant, and the same is true for the point $m_2
\in M$. But $\wt{a}(\exp t)$ acts by identity on $M^{U(1)} \subset
M$. Therefore in this case we also have $m_1=m_2=m_0$.
\endproof

\punkt Let $V = M_I^{U(1)} \subset M_I$ be the submanifold of fixed points
of the $U(1)$ action. Since the action preserves the complex structure on
$M_I$, the submanifold $V$ is complex.

\begin{lemma}
There exists a unique $U(1)$-invariant holomorphic map 
$$
\rho_M:M_I \to V
$$
such that $\rho_M|_{V} = \id$. 
\end{lemma}

\proof For every point $m \in M$ we must have $\rho_M(m) =
\displaystyle\lim_{t \to +\infty} \exp(i t \phi_M)$, which proves uniqueness.

To prove that $\rho_M$ thus defined is indeed holomorphic, notice that the
diagram 
$$
\begin{CD}
M_I           @>{0 \times \id}>>   \Delta \times M\\
 @V{\rho_M}VV                         @VV{\tilde{a}}V\\
V             @>>>                 M_I
\end{CD}
$$
is commutative. Since the action $\tilde{a}:\Delta \times M_I \to M_I$ is
holomorphic, so is the map $\rho_M$. 
\endproof 

\punkt Call the canonical map $\rho_M:M_I \to M_I^{U(1)}$ the {\em canonical
projection} of the regular Hodge manifold $M$ onto the submanifold $V \subset
M$ of fixed points. 

\begin{lemma}
The canonical projection $\rho_M:M \to M^{U(1)}$ is submersive, that
is, for every point $m \in M$ the differential $d\rho_M:T_mM \to
T_{\rho(m)}M^{U(1)}$ of the map $\rho_M$ at $m$ is surjective.
\end{lemma}

\proof Since $\rho_M|_{M^{U(1)}} = \id$, the differential $d\rho_M$
is surjective at points $m \in V \subset M$. Therefore it is
surjective on an open neighborhood $U \supset V$ of $V$ in $M$. For
any point $m \in M$ there exists a point $x \in \Delta$ such that $x
\cdot m \in U$. Since $\rho_M$ is $\Delta$-invariant, this implies that
$d\rho_M$ is surjective everywhere on $M$.  
\endproof

\punkt Let $\Theta(M/V)$ be the relative tangent bundle of the holomorphic
map $\rho:M \to V$. Let $\Theta(M)$ and $\Theta(V)$ be the tangent bundles
of $M$ and $V$ and consider the canonical exact sequence of complex bundles
$$
0 \Arrow \Theta(M/V) \Arrow \Theta(M) \overset{d\rho_M}{\Arrow}
\rho^*\Theta(V) \Arrow 0,
$$
where $d\rho_M$ is the differential of the projection $\rho_M:M \to
V$.

The quaternionic structure on $M$ defines a $\C$-linear map
$j:\Theta(M) \to \overline{\Theta}(M)$. Restricting to $\Theta(M/V)$ and
composing with $d\rho_M$, we obtain a $\C$-linear map $j:\Theta(M/V)
\to \rho^*\overline{\Theta}(V)$. 

\punkt \label{overline.T} 
Let $\tv$ be the total space of the bundle $\overline{\Theta}(V)$
complex-conjugate to the tangent bundle $\Theta(V)$, and let
$\rho:\tv \to V$ be the projection. Let the group $U(1)$ act on
$\tv$ by dilatation along the fibers of the projection $\rho$.

Since the canonical projection $\rho_M:M \to V$ is $U(1)$-invariant,
the differential $\phi_M$ of the $U(1)$-action defines a section
$$
\phi_M \in C^\infty(M,\Theta(V/M)).
$$
The section $j(\phi_M) \in C^\infty(M,\rho_M^*\overline{\Theta}(V))$
defines a map $\Lin_M:M \to \tv$ such that $\Lin_M \circ \rho =
\rho_M:M \to V$. Call the map $\Lin_M$ {\em the linearization} of the
regular Hodge manifold $M$.

\begin{prop}\label{linearization}
The linearization map $\Lin_M$ is a $U(1)$-equivariant open embedding.  
\end{prop}

\proof The map $j:\Theta(M/V) \to \rho^*\overline{\Theta}(V)$ is of
degree $1$ with respect to the $U(1)$-action, while the section
$\phi_M \in C^\infty(M,\Theta(M/V))$ is $U(1)$-invariant. Therefore
the map $\Lin_M$ is $U(1)$-equivariant.

Consider the differential $d\Lin_M:T_m(M) \to T_m(\tv)$ at a point $m \in V
\subset M$. We have
$$
d\Lin_M = d\rho_M \oplus d\rho_M \circ j:T_m(M) \to T_m(V) \oplus
\overline{T}_m(V)
$$
with respect to the decomposition $T_m(\tv) = T_m(V) \oplus
\overline{T}(V)$.  The tangent space $T_m$ is a regular quaternionic vector
space. Therefore the map $d\Lin_M$ is bijective at $m$ by
Lemma~\ref{regular.quaternionic}. Since $\Lin_M$ is bijective on $V$, this
implies that $\Lin_M$ is an open embedding on an open neighborhood $U
\subset M$ of the submanifold $V \subset M$.

To finish the proof of proposition, it suffices prove that the
linearization map $\Lin_M:M_I \to \tv$ is injective and \'etale on
the whole $M_I$. To prove injectivity, consider arbitrary two points
$m_1,m_2 \in M_I$ such that $\Lin_M(m_1)=\Lin_M(m_2)$.  There exists
a point $x \in \Delta \setminus \{0\}$ such that $x \cdot m_1, x
\cdot m_2 \in U$. The map $\Lin_M$ is $U(1)$-equivariant and
holomorphic, therefore it is $\Delta$-equivariant, and we have
$$
\Lin_M(x \cdot m_1) = x \cdot \Lin_M(m_1) = x \cdot \Lin_M(m_2) =
\Lin_M(x \cdot m_2).
$$ 
Since the map $\Lin_M:U \to \tv$ is injective, this implies that $x
\cdot m_1 = x \cdot m_2$. By Lemma~\ref{regular} the action map
$x:M_I \to M_I$ is injective. Therefore this is possible only if
$m_1 = m_2$, which proves injectivity.

To prove that the linearization map is \'etale, note that by
Lemma~\ref{regular} the action map $x:M_I \to M_I$ is not only
injective, but also \'etale. Since $\Lin_M$ is \"etale on $U$, so is
the composition $\Lin_M \circ x:M_I \to U \to \tv$ is \'etale. Since
$\Lin_M \circ x = x \circ \Lin_M$, the map $\Lin_M:M_I \to \tv$ is
\'etale at the point $m_1 \in M_I$.

Thus the linearization map is also injective and \'etale on the
whole of $M_I$. Hence it is indeed an open embedding, which proves
the propostion.  
\endproof

\subsection{Linear Hodge manifold structures}

\punkt By Proposition~\ref{linearization} every regular Hodge manifold $M$
admits a canonical open embedding $\Lin_M:M \to \tv$ into the total space
$\tv$ of the (complex-conjugate) tangent bundle to its fixed points
submanifold $V \subset M$. This embedding induces a Hodge manifold
structure on a neighborhood of the zero section $V \subset \tv$. 

In order to use the linearization construction, we will need a
characterization of all Hodge manifold structures on neighborhoods of $V
\subset \tv$ obtained in this way (see \ref{lin.def}). It is convenient to
begin with an invariant characterization of the linearization map $\Lin_M:M
\to \tv$.

\punkt \label{tau} 
Let $V$ be an arbitrary complex manifold, let $\tv$ be the total space
of the complex-conjugate to the tangent bundle $\Theta(V)$ to
$V$, and let $\rho:\tv \to V$ be the canonical projection. Contraction with
the tautological section of the bundle $\rho^*\overline{\Theta(V)}$ defines
for every $p$ a bundle map 
$$
\tau:\rho^*\Lambda^{p+1}(V,\C) \to \rho^*\Lambda^{p}(V,\C),  
$$
which we call {\em the tautological map}. In particular, the
induced map
$$
\tau:C^\infty(V,\Lambda^{0,1}(V)) \to C^\infty(\tv,\C)
$$
identifies the space $C^\infty(V,\Lambda^{0,1}(V))$ of smooth $(0,1)$-forms
on $V$ with the subspace in $C^\infty(\tv,\C)$ of function linear along the
fibers of the projection $\tv \to V$.

\punkt Let now $M$ be a Hodge manifold. Let $V \subset M_I$ be the
complex submanifold of $U(1)$-fixed points, and let $\rho_M:M \to V$
be the canonical projection. Assume that $M$ is equipped with a
smooth $U(1)$-equivariant map $f:M \to \tv$ such that $\rho_M = \rho
\circ f$. Let $\bar\6_I$ be the Dolbeault differential for the
preferred complex structure $M_I$ on $M$, and let $\phi \in
\Theta(M/V)$ be the differential of the $U(1)$-action on $M$. Let
also $j:\Lambda^{0,1}(M_I) \to \Lambda^{1,0}(M_I)$ be the map
induced by the quaternionic structure on $M$.

\begin{lemma}\label{lin.char}
The map $f:M \to \tv$ coincides with the linearization map if and only if
for every $(0,1)$-form $\alpha \in C^\infty(V,\Lambda^{0,1}(V))$ we have 
\begin{equation}\label{L}
f^*\tau(\alpha) = \langle \phi, j(\rho_M^*\alpha) \rangle.
\end{equation}
Moreover, if $f = \Lin_M$, then we have 
\begin{equation}\label{LL}
f^*\tau(\beta) = \langle \phi, j(f^*\beta) \rangle
\end{equation}
for every smooth section $\beta \in C^\infty(\tv,\rho^*\Lambda^1(V,\C))$. 
\end{lemma}

\proof Since functions on $\tv$ linear along the fibers separate points,
the correspondence
$$
f^* \circ \tau:C^\infty(V,\Lambda^{0,1}(V)) \to C^\infty(M,\C)
$$
characterizes the map $f$ uniquely, which proves the ``only if''
part of the first claim. Since by assumption $\rho_M = \rho \circ
f$, the equality \eqref{L} is a particular case of \eqref{LL} with
$\beta = \rho^*\alpha$. Therefore the ``if'' part of the first claim
follows from the second claim, which is a rewriting of the
definition of the linearization map $\Lin_M:M \to \tv$ (see
\ref{overline.T}).  
\endproof

\punkt Let now $\Lin_M:M \to \tv$ be the linearization map for the
regular Hodge manifold $M$. Denote by $U \subset \tv$ the image of
$\Lin_M$. The subset $U \subset \tv$ is open and
$U(1)$-invariant. In addition, the isomorphism $\Lin_M:M \to U$
induces a regular Hodge manifold structure on $U$.

Donote by $\Lin_U$ the linearization map for the regular Hodge
manifold $U$. Lemma~\ref{lin.char} implies the following.

\begin{corr}\label{lin.lin}
We have $\Lin_M \circ \Lin_U = \Lin_M$, thus the linearization map
$\Lin_U:U \to \tv$ coincides with the given embedding $U \hookrightarrow \tv$. 
\end{corr}

\proof Let $\alpha \in C^\infty(V,\Lambda^{0,1}(V)$ be a
$(0,1)$-form on $V$. By Lemma~\ref{lin.char} we have
$\Lin_U^*\tau(\alpha) = \langle \phi_U, j_U(\rho_U^*\alpha)\rangle$,
and it suffices to prove that
$$
\Lin_M^*(\Lin_U^*(\tau(\alpha))) = \langle \phi_M, j_M(\rho_M^*\alpha)
\rangle.
$$
By definition we have $\rho_M = \rho_U \circ \Lin_M$. Moreover, the
map $\Lin_M$ is $U(1)$-equivariant, therefore it sends $\phi_M$ to
$\phi_U$. Finally, by definition it commutes with the quaternionic
structure map $j$. Therefore 
$$
\Lin_M^*(\Lin_U^*(\tau(\alpha))) = \Lin_M^*(\langle \phi_U,
j_U(\rho_U^*\alpha) \rangle)= \langle \phi_M, j_M(\rho_M^*\alpha)
\rangle,
$$
which proves the corollary. 
\endproof 

\punkt \defn \label{lin.def}
Let $U \subset \tv$ be an open $U(1)$-invariant neighborhood of the zero
section $V \subset \tv$. A Hodge manifold structure on $\tv$ is called {\em
linear} if the associated linearization map $\Lin_U:U \to \tv$ coincides
with the given embedding $U \hookrightarrow \tv$.

By Corollary~\ref{lin.lin} every Hodge manifold structure on a subset $U
\subset \tv$ obtained by the linearization construction is linear. 

\punkt We finish this section with the following simple observation, which
we will need in the next section.

\begin{lemma}\label{aux2}
Keep the notations of Lemma~\ref{lin.char}. Moreover, assume given a
subspace $\A \subset C^\infty(\tv,\rho^*\Lambda^1(V,\C))$ such that
the image of $\A$ under the restriction map
$$
\Res:C^\infty(\tv,\rho^*\Lambda^1(V,\C)) \to
C^\infty(V,\Lambda^1(V,\C)) 
$$
onto the zero section $V \subset \tv$ is the whole space
$C^\infty(V,\Lambda^1(V,\C))$. If \eqref{LL} holds for every section
$\beta \in \A$, then it holds for every smooh section 
$$
\beta \in C^\infty(\tv,\rho^*\Lambda^1(V,\C)).
$$
\end{lemma}

\proof By assumptions sections $\beta \in \A$ generate the
restriction of the bundle $\rho^*\Lambda^1(V,\C)$ onto the zero
section $V \subset \tv$. Therefore there exists an open neighborhood
$U \subset \tv$ of the zero section $V \subset \tv$ such that the
$C^\infty(U,\C)$-submodule 
$$
C^\infty(U,\C) \cdot \A \subset C^\infty(U,\rho^*\Lambda^1(V,\C)) 
$$
is dense in the space $C^\infty(U,\rho^*\Lambda^1(V,\C))$ of smooth
sections of the pull-back bundle $\rho^*\Lambda^1(V,\C)$. Since
\eqref{LL} is continuous and linear with respect to multiplication
by smooth functions, it holds for all sections $\beta \in
C^\infty(U,\rho^*\Lambda^1(V,\C))$. Since it is also compatible with
the natural unit disc action on $M$ and $\tv$, it holds for all
sections $\beta \in C^\infty(\tv,\rho^*\Lambda^1(V,\C))$ as well.
\endproof

\section{Tangent bundles as Hodge manifolds}\label{section.5}

\subsection{Hodge connections}

\punkt The linearization construction reduces the study of arbitrary
regular Hodge manifolds to the study of linear Hodge manifold
structures on a neighborhood $U \subset \tv$ of the zero section $V
\subset \tv$ of the total space of the complex conjugate to the
tangent bundle of a complex manifold $V$. In this section we use the
theory of Hodge bundles developed in Subsection~\ref{hb.sub} in
order to describe Hodge manifold structures on $U$ in terms of
connections on the locally trivial fibration $U \to V$ of a certain
type, which we call {\em Hodge connections} (see
\ref{hodge.con}). It is this description, given in
Proposition~\ref{equiv}, which we will use in the latter part of the
paper to classify all such Hodge manifold structures.

\punkt \label{conn}
We begin with some preliminary facts about connections on locally
trivial fibrations. Let $f:X \to Y$ be an arbitrary smooth map of
smooth manifolds $X$ and $Y$. Assume that the map $f$ is submersive,
that is, the codifferential $\del_f:f^*\Lambda^1(Y) \to
\Lambda^1(X)$ is an injective bundle map.  Recall that a {\em
connection} on $f$ is by definition a splitting $\Theta:\Lambda^1(X)
\to f^*\Lambda^1(Y)$ of the canonical embedding $\del_f$.

Let $d_X$ be the de Rham differential on the smooth manifold
$X$. Every connection $\Theta$ on $f:X \to Y$ defines an algebra
derivation
$$
D = \Theta \circ d_X:\Lambda^0(X) \to f^*\Lambda^1(Y), 
$$
satisfying 
\begin{equation}\label{conn.eq}
D \rho^* h = \rho^* d_Y h
\end{equation}
for every smooth function $h \in C^\infty(Y,\R)$. Vice versa, by the
universal property of the cotangent bundle (Lemma~\ref{universal})
every algebra derivation $D:\Lambda^0(X) \to \Lambda^1(Y)$
satisfying \eqref{conn.eq} comes from a unique connection $\Theta$
on $f$.

\punkt Recall also that a connection $\Theta$ is called {\em flat}
if the associated derivation $D$ extends to an algebra
derivation
$$
D:f^*\Lambda^\cdot(Y) \to f^*\Lambda^{\cdot+1}(Y)
$$
so that $D \circ D = 0$. The splitting $\Theta:\Lambda^1(X) \to
f^*\Lambda^1(Y)$ extends in this case to an algebra map 
$$
\Theta:\Lambda^\cdot(X) \to f^*\Lambda^\cdot(Y)
$$
compatible with the de Rham differential $d_X:\Lambda^\cdot(X) \to
\Lambda^{\cdot+1}(X)$. 

\punkt We will need a slight generalization of the notion of connection. 

\defn Let $f:X \to Y$ be a smooth submersive morphism of complex
manifolds. A {\em $\C$-valued connection} $\Theta$ on $f$ is a
splitting $\Theta:\Lambda^1(Y,\C) \to f^*\Lambda^1(X,\C)$ of the
codifferential map $\del f:f^*\Lambda^1(Y,\C) \to \Lambda^1(X,\C)$
of complex vector bundles. A $\C$-valued connection $\Theta$ is
called {\em flat} if the associated algebra derivation 
$$
D = \Theta \circ d_X:\Lambda^0(X,\C) \to f^*\Lambda^1(Y,\C)
$$
extends to an algebra derivation 
$$
D:f^*\Lambda^\cdot(Y,\C) \to f^*\Lambda^{\cdot+1}(Y,\C)
$$
satisfying $D \circ D = 0$. 

As in \ref{conn}, every derivation $D:\Lambda^0(X,\C) \to
f^*\Lambda^1(Y,\C)$ satisfying \eqref{conn.eq} comes from a unique
$\C$-valued connection $\Theta$ on $f:X \to Y$. 

\rem By definition for every flat connection on $f:X \to Y$ the
subbundle of horizontal vectors in the the tangent bundle
$\Theta(X)$ is an involutive distribution. By Frobenius Theorem this
implies that the connection defines locally a trivialization of the
fibration $f$.

This is no longer true for flat $\C$-valued connections: the
subbundle of horizontal vectors in $\Theta(X) \otimes \C$ is only
defined over $\C$, and the Frobenius Theorem does not apply. One can
try to correct this by replacing the splitting
$\Theta:\Lambda^1(X,\C) \to f^*\Lambda^1(Y,\C)$ with its real part
$\Re\Theta:\Lambda^1(X) \to \Lambda^1(Y)$, but this real part is, in
general, no longer flat.

\punkt For every $\C$-valued connection $\Theta:\Lambda^1(X,\C) \to
f^*\Lambda^1(Y,\C)$ on a fibration $f:X \to Y$ the kernel
$\Ker\Theta \subset \Lambda^1(X,\C)$ is canonically isomorphic to
the quotient $\Lambda^1(X,\C)/\del_f(f^*\Lambda^1(Y,\C))$, and the
composition
$$
R = \Theta \circ d_X:\Lambda^1(X,\C)/\del_f(f^*\Lambda^1(Y,\C) \cong
\Ker\Theta \to f^*\Lambda^2(Y,\C) 
$$
is in fact a bundle map. This map is called {\em the curvature} of
the $\C$-valued connection $\Theta$. The connection $\Theta$ is flat
if and only if its curvature $R$ vanishes.

\punkt \label{conn.setup}
Let now $M$ be a complex manifold, and let $U \subset \tm$ be an
open neighborhood of the zero section $M \subset \tm$ in the total
space $\tm$ of the complex-conjugate to the tangent bundle to
$M$. Let $\rho:U \to M$ be the natural projection. Assume that $U$
is invariant with respect to the natural action of the unit disc
$\Delta \subset \C$ on $\tm$.

\punkt Since $M$ is complex, by \ref{de.Rham} the bundle
$\Lambda^1(M,\C)$ is equipped with a Hodge bundle structure of
weight $1$. The pullback bundle $\rho^*\Lambda^1(M,\C)$ is then also
equipped with a weight $1$ Hodge bundle structure. 

Our description of the Hodge manifold structures on the subset $U
\in \tm$ is based on the following notion.

\defn \label{hodge.con} 
A {\em Hodge connection} on the pair $\langle M, U\rangle$ is
a $\C$-valued connection on $\rho:U \to M$ such that the associated
derivation 
$$
D:\Lambda^0(U,\C) \to \rho^*\Lambda^1(M,\C)
$$
is weakly Hodge in the sense of \ref{w.hodge}. A Hodge connection is
called {\em flat} if it extends to a weakly Hodge derivation 
$$
D:\rho^*\Lambda^\cdot(M,\C) \to \rho^*\Lambda^{\cdot+1}(M,\C)
$$
satisfying $D \circ D = 0$. 

\punkt Assume given a flat Hodge connection $D:\Lambda^0(U,\C) \to
\rho^*\Lambda^1(M,\C)$ on the pair $\langle U,M \rangle$, and assume
in addition that the derivation $D$ is holonomic in the sense of
\ref{holonomic}. Then the pair $\langle D,
\rho^*\Lambda^1(M,\C)\rangle$ defines by
Proposition~\ref{explicit.hodge} a Hodge manifold structure on $U$.

It turns out that every Hodge manifold structure on $U$ can be
obtained in this way. Namely, we have the following. 

\begin{prop}\label{equiv}
There correspondence $D \mapsto \langle \rho^*\Lambda^1(M,\C), D
\rangle$ is a bijection between the set of all flat Hodge
connections $D$ on the pair $\langle U, M\rangle$ such that
$D:\Lambda^0(U,\C) \to \rho^*\Lambda^1(M,\C)$ is holonomic in the
sense of \ref{holonomic}, and the set of all Hodge manifold
structures on the $U(1)$-manifold $U$ such that the projection
$\rho:U_I \to M$ is holomorphic for the preferred complex structure
$U_I$ on $U$.
\end{prop}

\punkt The crucial part of the proof of Proposition~\ref{equiv} is
the following observation. 

\begin{lemma}\label{ident}
Assume given a Hodge manifold structure on the $U(1)$-manifold $U
\subset \tm$. Let $\del_\rho:\rho^*\Lambda^1(M,\C) \to
\Lambda^1(U,\C)$ be the codifferential of the projection $\rho:U \to
M$, and let $P:\Lambda^1(U,\C) \to \Lambda^{0,1}(U_J)$ be the
canonical projection. The bundle mao given by the composition
$$
P \circ \del_\rho:\rho^*\Lambda^1(M,\C) \to \Lambda^{0,1}(U_J)
$$
is an isomorphism of complex vector bundles. 
\end{lemma}

\proof Since the bundles $\rho^*\Lambda^1(M,\C)$ and
$\Lambda^{0,1}(U_J)$ are of the same rank, and the maps
$\del_\rho:\rho^*\Lambda^1(M,\C) \to \Lambda^1(U,\C)$ and $P \circ
\del_\rho$ are equivariant with respect to the action of the unit
disc on $U$, it suffices to prove the claim on $M \subset U$. Let $m
\in M$ be an arbitrary point, and let $V = T^*_m\tm$ be the
cotangent bundle at $m$ to the Hodge manifold $U \subset \tm$. Let
also $V^0 \subset V$ be the subspace of $U(1)$-invariant vectors in
$V$.

The space $V$ is an equivariant quaternionic vector space. Moreover,
the fibers of the bundles $\rho^*\Lambda^1(M,\C)$ and
$\Lambda^{0,1}(U_J)$ at the point $m$ are complex vector spaces, and
we have canonical identifications
\begin{align*}
\rho^*\Lambda^1(M,\C)|_m &\cong V^0_I \oplus \overline{V^0_I},\\
\Lambda^{0,1}(U_J)|_m &\cong V_J. 
\end{align*}
Under these identifications the map $P \circ \del_\rho$ at the point
$m$ coincides with the action map $V^0_I \oplus \overline{V^0_I} \to
V_J$, which is invertible by Lemma~\ref{regular.quaternionic}.
\endproof

\punkt By Proposition~\ref{explicit.hodge} every Hodge manifold
structure on $U$ is given by a pair $\langle \E, D\rangle$ 
of a Hodge bundle $\E$ on $U$ of weight $1$ and a holonomic
derivation $D:\Lambda^0(U,\C) \to \E$. Lemma~\ref{ident} gives an
isomorphism $\E \cong \rho^*\Lambda^1(M,\C)$, so that $D$ becomes a
flat $\C$-valued connection on $U$ over $M$. To prove
Proposition~\ref{equiv} it suffices now to prove the following. 

\begin{lemma}
The complex vector bundle isomorphism 
$$
P \circ \del_\rho:\rho^*\Lambda^1(M,\C) \to \Lambda^{0,1}(U_J)
$$
associated to a Hodge manifold structure on $U$ is compatible with
the Hodge bundle structures if and only if the projection $\rho:U_I
\to M$ is holomorphic for the preferred complex structure $U_I$ on
$U$.
\end{lemma}

\proof The preferred complex structure $U_I$ induces a Hodge bundle
structure of weight $1$ on $\Lambda^1(U,\C)$ by \ref{de.Rham}, and
the canonical projection $P:\Lambda^1(U,\C) \to \Lambda^{0,1}(U_J)$
is compatible with the Hodge bundle structures by
\ref{lambda_j=lambda_i}. If the projection $\rho:U_I \to M$ is
holomorphic, then the codifferential
$\del_\rho:\rho^*\Lambda^1(M,\C) \to \Lambda^1(U,\C)$ sends the
subbundles $\rho^*\Lambda^{1,0}(M),\rho^*\Lambda^{0,1}(M) \subset
\rho^*\Lambda^1(M,\C)$ into, respectively, the subbundles
$\Lambda^{1,0}(U_I),\Lambda^{0,1}(U_I) \subset
\Lambda^1(U,\C)$. Therefore the map $\del_\rho:\rho^*\Lambda^1(M,\C)
\to \Lambda^1(U,\C)$ is compatible with the Hodge bundle structures,
which implies the ``if'' part of the lemma.

To prove the ``only if'' part, assume that $P \circ \del_\rho$ is a
Hodge bundle isomorphism. Since the complex conjugation
$\nu:\overline{\Lambda^{0,1}(U_I)} \to \Lambda^{1,0}(U_J)$ is
compatible with the Hodge bundle structures, the projection
$\overline{P}:\Lambda^1(U,\C) \to \Lambda^{1,0}(U_J)$ and the
composition $\overline{P} \circ \del_\rho:\rho^*\Lambda^1(M,\C) \to
\Lambda^{1,0}(U_J)$ are also compatible with the Hodge bundle
structures. Therefore the map
$$
P \oplus \overline{P}:\Lambda^1(U,\C) \to \Lambda^{1,0}(U_J) \oplus
\Lambda^{0,1}(U_J) 
$$
is a Hodge bundle isomorphism, and the composition 
$$
\del_\rho \circ (P \oplus \overline{P}):\rho^*\Lambda^1(M,\C) \to
\Lambda^1(U,\C)
$$ 
is a Hodge bundle map. Therefore the codifferential
$\del_\rho:\rho^*\Lambda^1(M,\C) \to \Lambda^1(U,\C)$ is compatible
with the Hodge bundle structures. This means precisely that the
projection $\rho:U_I \to M$ is holomorphic, which finishes the proof
of the lemma and of Proposition~\ref{equiv}.  
\endproof

\subsection{The relative de Rham complex of $U$ over
$M$}\label{relative.de.rham.sub} 

\punkt Keep the notation of the last subsection. To use
Proposition~\ref{equiv} in the study of Hodge manifold structures on
the open subset $U \subset \tm$, we will need a way to check whether
a given Hodge connection on the pair $\langle U,M \rangle$ is
holonomic in the sense of \ref{holonomic}. We will also need to
rewrite the linearity condition \ref{lin.def} for a Hodge manifold
structure on $U$ in terms of the associated Hodge connection $D$. To
do this, we will use the so-called {\em relative de Rham complex} of
$U$ over $M$. For the convenience of the reader, and to fix
notation, we recall here its definition and main properties.

\punkt Since the projection $\rho:U \to M$ is submersive, the
codifferential 
$$
\del_\rho:\rho^*\Lambda^1(M,\C) \to \Lambda^1(U,\C)
$$ 
is injective. The {\em relative cotangent bundle} $\Lambda^1(U/M,C)$
is by definition the quotient bundle
$$
\Lambda^1(U/M,\C) = \Lambda^1(U,\C)/\del_\rho(\rho^*\Lambda^1(M,\C)). 
$$
Let $\pi:\Lambda^1(U,\C) \to \Lambda^1(U/M,\C)$ be the natural
projection. We have by definition the short exact sequence 
\begin{equation}\label{ex.seq}
\begin{CD}
0 @>>> \rho^*\Lambda^1(M,\C) @>\del_\rho>> \Lambda^1(U,\C) @>\pi>>
\Lambda^1(U/M,\C) @>>> 0
\end{CD}
\end{equation}
of complex vector bundles on $U$. 

\punkt The composition $d^r = \pi \circ d_U$ of the de Rham differential
$d_U$ with the projection $\pi$ is an algebra derivation 
$$
d^r:\Lambda^0(U,\C) \to \Lambda^1(U/M,\C),
$$
called {\em the relative de Rham differential}. It is a first order
differential operator, and $d^rf = 0$ if and only if the smooth
function $f:U \to \C$ factors through the projection $\rho:U \to
M$. 

Let $\Lambda^\cdot(U/M,\C)$ be the exterior algebra of the bundle
$\Lambda^1(U/M,\C)$. The projection $\pi$ extends to an algebra
map 
$$
\pi:\Lambda^\cdot(U,\C) \to \Lambda^\cdot(U/M,\C).
$$
The differential $d^r$ extends to an algebra derivation 
$$
d^r:\Lambda^\cdot(U/M,\C) \to \Lambda^{\cdot+1}(U/M,\C)
$$
satisfying $d^r \circ d^r = 0$, and we have $\pi \circ d_U = d^r
\circ \pi$. The differential graded algebra $\langle
\Lambda^\cdot(U/M,\C),d^r\rangle$ is called {\em the relative de
Rham complex} of $U$ over $M$. 

\punkt Since the relative de Rham differential $d^r$ is linear with
respect to multiplication by functions of the form $\rho^*f$ with $f
\in C^\infty(M,\C)$, it extends canonically to an operator
$$
d^r:\rho^*\Lambda^i(M,\C) \otimes \Lambda^\cdot(U/M,\C) \to
\rho^*\Lambda^i(M,\C) \otimes \Lambda^{\cdot+1}(U/M,\C). 
$$
The two-step filtration $\rho^*\Lambda^1(M,\C) \subset
\Lambda^1(U,\C)$ induces a filtration on the de Rham complex
$\Lambda^\cdot(U,\C)$, and the $i$-th associated graded quotient
of this filtration is isomorphic to the complex $\langle
\rho^*\Lambda^i(M,\C) \otimes \Lambda^\cdot(U/M,\C), d^r \rangle$.

\punkt Since $U \subset \tm$ lies in the total space of the
complex-conjugate to the tangent bundle to $M$, we have a canonical
algebra isomorphism
$$
\can:\overline{\rho^*\Lambda^\cdot(M,\C)} \to \Lambda^\cdot(U/M,\C).
$$
Let $\tau:C^\infty(M,\Lambda^1(M,\C)) \to C^\infty(U,\C)$ be
the tautological map sending a $1$-form to the corresponding linear
function on $\tm$, as in \ref{tau}. Then for every smooth
$1$-form $\alpha \in C^\infty(M,\Lambda^1(M,\C))$ we have
\begin{equation}\label{can.and.tau}
\can(\rho^*\alpha) = d^r\tau(\alpha). 
\end{equation}

\punkt The complex vector bundle $\Lambda^1(U/M,\C)$ has a natural
real structure, and it is naturally $U(1)$-equivariant. Moreover,
the decomposition $\Lambda^1(M,\C) = \Lambda^{1,0}(M) \oplus
\Lambda^{0,1}(M)$ induces a decomposition
$$
\Lambda^1(U/M,\C) = \can(\Lambda^{1,0}(M)) \oplus
\can(\Lambda^{0,1}(M)).
$$
This allows to define, as in \ref{de.Rham}, a canonical Hodge bundle
structure of weight $1$ on $\Lambda^1(U/M,\C)$. It gives rise to a
Hodge bundle structure on $\Lambda^i(U/M,\C)$ of wieght $i$, and the
relative de Rham differential $d^r:\Lambda^\cdot(U/M,\C) \to
\Lambda^{\cdot+1}(U/M,\C)$ is weakly Hodge.

\punkt The canonical isomorphism
$\can:\overline{\rho^*\Lambda^1(M,\C)} \to \Lambda^1(U/M,\C)$ is not
compatible with the Hodge bundle structures. The reason for this is
that the real structure on the Hodge bundles $\Lambda^\cdot(U/M,\C)$
is, by definition \ref{de.Rham}, twisted by $\iota^*$, where
$\iota:\tm \to \tm$ is the action of $-1 \in U(1) \subset
\C$. Therefore, while $\can$ is $U(1)$-equivariant, it is not
real. To correct this, introduce an involution
$\zeta:\Lambda^1(M,\C) \to \Lambda^1(M,\C)$ by
\begin{equation}\label{zeta}
\zeta = \begin{cases} \id &\text{ on }\Lambda^{1,0}(M) \subset
\Lambda^1(M,\C) \\
-\id &\text{ on }\Lambda^{0,1}(M) \subset \Lambda^1(M,\C)
\end{cases}
\end{equation}
and set 
\begin{equation}\label{eta}
\eta = \can \circ
\rho^*\overline{\zeta}:\rho^*\overline{\Lambda^1(M,\C)} \to
\rho^*\overline{\Lambda^1(M,\C)} \to \Lambda^1(U/M,\C)
\end{equation}
Unlike $\can$, the map $\eta$ preserves the Hodge bundle structures.

It will also be convenient to twist the tautological map
$\tau:\rho^*\Lambda^1(M,\C) \to \Lambda^0(U,\C)$ by the involution
$\zeta$. Namely, define a map $\sigma:\rho^*\Lambda^1(M,\C) \to
\Lambda^0(U,\C)$ by
\begin{equation}\label{sigma1}
\sigma = \tau \circ \rho^*\overline{\zeta}:
\rho^*\overline{\Lambda^1(M,\C)} \to
\rho^*\overline{\Lambda^1(M,\C)} \to \Lambda^0(U/M,\C)  
\end{equation}
By \eqref{can.and.tau} the twisted tautological map $\sigma$ and the
canonical map $\eta$ satisfy
\begin{equation}\label{eta.and.sigma}
\eta(\rho^*\alpha) = d^r\sigma(\alpha) 
\end{equation}
for every smooth $1$-form $\alpha \in C^\infty(M,\Lambda^1(M,\C))$. 

\punkt Let $\phi \in \Theta(U)$ be the differential of the canonical
$U(1)$-action on $U \subset \tm$. The vector field $\phi$ is real
and tangent to the fibers of the projection $\rho:U \to
M$. Therefore the contraction with $\phi$ defines an algebra
derivation 
\begin{align*}
\Lambda^{\cdot+1}(U/M,\C) &\to \Lambda^\cdot(U/M,\C)\\
\alpha &\mapsto \langle \phi, \alpha \rangle
\end{align*}
The following lemma gives a relation between this derivation, the
canonical weakly Hodge map $\eta:\rho^*\Lambda^1(M,\C) \to
\Lambda^1(U/M,\C)$ given by \eqref{eta}, and the tautological map
$\tau:\rho^*\Lambda^1(M,\C) \to \Lambda^0(U,\C)$.

\begin{lemma}\label{phi.and.tau}
For every smooth section $\alpha \in
C^\infty(U,\rho^*\Lambda^1(M,\C))$ we have 
$$
\sqrt{-1}\tau(\alpha) = \langle \phi, \eta(\alpha) \rangle \in
C^\infty(U,\C). 
$$
\end{lemma}

\proof Since the equality that we are to prove is linear with
respect to multiplication by smooth functions on $U$, we may assume
that the section $\alpha$ is the pull-back of a smooth $1$-form
$\alpha \in C^\infty(M,\Lambda^1(M,\C))$. The Lie derivative
$\LL_\phi:\Lambda^\cdot(U,\C) \to \Lambda^\cdot(U,\C)$ with respect
to the vector field $\phi$ is compatible with the projection
$\pi:\Lambda^\cdot(U,\C) \to \Lambda^\cdot(U/M,\C)$ and defines
therefore an algebra derivation $\LL_\phi:\Lambda^\cdot(U/M,\C) \to
\Lambda^\cdot(U/M,\C)$. The Cartan homotopy formula gives
\begin{equation}\label{eqqq.1}
\LL_\phi \tau(\alpha) = \langle \phi, d^r\tau(\alpha) \rangle. 
\end{equation}
The function $\tau(\alpha)$ on $\tm$ is by definition $\R$-linear
along the fibers of the projection $\rho:\tm \to M$. The subspace
$\tau(C^\infty(M,\Lambda^1(M,\C))) \subset C^\infty(U,\C)$ of such
functions decomposes as 
$$
\tau(C^\infty(M,\Lambda^1(M,\C))) =
\tau(C^\infty(M,\Lambda^{1,0}(M))) \oplus
\tau(C^\infty(M,\Lambda^{0,1}(M,\C))), 
$$
and the group $U(1)$ acts on the components with weight $1$ and
$-1$. Therefore the derivative $\LL_\phi$ of the $U(1)$-action acts
on the components by multiplication with $\sqrt{-1}$ and
$-\sqrt{-1}$. By definition of the involution $\zeta$ (see
\eqref{zeta}) this can be written as 
\begin{equation}\label{eqqq.2}
\LL_\phi \tau(\alpha) = \sqrt{-1}\tau(\zeta(\alpha)).
\end{equation}
On the other hand, by \eqref{can.and.tau} and the definition of the
map $\eta$ we have 
\begin{equation}\label{eqqq.3}
d^r\tau(\alpha) = \can(\alpha) = \eta(\zeta(\alpha)). 
\end{equation}
Substituting \eqref{eqqq.2} and \eqref{eqqq.3} into \eqref{eqqq.1}
gives  
$$
\sqrt{-1}\tau(\zeta(\alpha)) = \langle \phi, \eta(\zeta(\alpha)), 
$$
which is equivalent to the claim of the lemma. 
\endproof 

\subsection{Holonomic Hodge connections}

\punkt We will now describe a convenient way to check whether a
given Hodge connection $D$ on the pair $\langle U,M \rangle$ is
holonomic in the sense of \ref{holonomic}. To do this, we proceed as
follows.

Consider the restriction $\Lambda^1(U,\C)|_M$ of the bundle
$\Lambda^1(U,\C)$ to the zero section $M \subset U \subset \tm$, and
let 
$$
\Res:\Lambda^1(U,\C)|_M \to \Lambda^1(M,\C)
$$
be the restriction map. The kernel of the map $\Res$ coincides with
the conormal bundle to the zero section $M \subset U$, which we
denote by $S^1(M,\C)$. The map $\Res$ splits the restriction of
exact sequence \eqref{ex.seq} onto the zero section $M \subset U$,
and we have the direct sum decomposition
\begin{equation}\label{hor.vert}
\Lambda^1(U,\C)|_M = S^1(M,\C) \oplus \Lambda^1(M,\C).
\end{equation}

\punkt\label{S1} 
The $U(1)$-action on $U \subset \tm$ leaves the zero section $M
\subset U$ invariant and defines therefore a $U(1)$-action on the
conormal bundle $S^1(M,\C)$. Together with the usual real structure
twisted by the action of the map $\iota:\tm \to \tm$, this defines a
Hodge bundle structure of weight $0$ on the bundle $S^1(M,\C)$.

Note that the automorphism $\iota:\tm \to \tm$ acts as $-\id$ on the
Hodge bundle $S^1(M,\C)$, so that the real structure on $S^1(M,\C)$
is minus the usual one. Moreover, as a complex vector bundle the
conormal bundle $S^1(M,\C)$ to $M \subset \tm$ is canonically
isomorphic to the cotangent bundle $\Lambda^1(M,\C)$. The Hodge type
grading on $S^1(M,\C)$ is given in terms of this isomorphism by
$$
S^1(M,\C) = S^{1,-1}(M) \oplus S^{-1,1}(M) \cong \Lambda^{1,0}(M)
\oplus \Lambda^{0,1}(M) = \Lambda^1(M,\C).
$$

\punkt Let 
$$
C_{lin}^\infty(U,\C) = \tau(C^\infty(M,\Lambda^1(M,\C))) \subset
C^\infty(U,\C) 
$$
be the subspace of smooth functions linear along the fibers of the
canonical projection $\rho:U \subset \tm \to M$. The relative de
Rham differential defines an isomorphism
\begin{equation}\label{iso}
d^r:C_{lin}^\infty(U,\C) \to C^\infty(M,S^1(M,\C)).
\end{equation}
This isomorphism is compatible with the canonical Hodge structures of weight
$0$ on both spaces, and it is linear with respect to multiplication
by smooth functions $f \in C^\infty(M,\C)$. 

\punkt \label{pr.part}
Let now $D:\Lambda^0(U,\C) \to \rho^*\Lambda^1(M,\C)$ be a Hodge
connection on the pair $\langle U,M \rangle$, and let
$\Theta:\Lambda^1(U,\C) \to \rho^*\Lambda^(M,\C)$ be the
corresponding bundles map. Since $D$ is a $\C$-valued connection,
the restriction $\Theta|_M$ onto the zero section $M \subset M$
decomposes as
\begin{equation}\label{hor.vert.1}
\Theta = D_0 \oplus \id:S^1(M,\C) \oplus \Lambda^1(M,\C) \to
\Lambda^1(M,\C) 
\end{equation}
with respect to the direct sum decomposition \eqref{hor.vert} for a
certain bundle map $D_0:S^1(M,\C) \to \Lambda^1(M,\C)$.  

\defn The bundle map $D_0:S^1(M,\C) \to \Lambda^1(M,\C)$ is called
the {\em principal part} of the Hodge connection $D$.

\punkt Consider the map $D_0:C^\infty(M,S^1(M,\C)) \to
C^\infty(M,\Lambda^1(M,\C))$ on the spaces of smooth sections
induced by the principal part $D_0$ of a Hodge connection $D$.
Under the isomorphism \eqref{iso} this map coincides with the
restriction of the composition
$$
\Res \circ D:C^\infty(U,\C) \to C^\infty(U,\rho^*\Lambda^1(M,\C))
\to C^\infty(M,\Lambda^1(M,\C))
$$
onto the subspace $C^\infty_{lin}(U,\C) \subset
C^\infty(U,\C)$. Each of the maps $\Res$, $D$ is weakly Hodge, so
that this composition also is weakly Hodge. Since the isomorphism
\eqref{iso} is compatible with the Hodge bundle structures, this
implies that the principal part $D_0$ of the Hodge connection $D$ is
a weakly Hodge bundle map. In particular, it is purely imaginary
with respect to the usual real structure on the conormal bundle
$S^1(M,\C)$.

\punkt We can now formulate the main result of this subsection. 

\begin{lemma}\label{aux1}
A Hodge connection $D:\Lambda^0(U,\C) \to \rho^*\Lambda^1(M,\C)$ on
the pair $\langle U,M \rangle$ is holonomic in the sense of
\ref{holonomic} on an open neighborhood $U_0 \subset U$ of the zero
section $M \subset U$ if and only if its principal part
$D_0:S^1(M,\C) \to \Lambda^1(M,\C)$ is a complex vector bundle
isomorphism.
\end{lemma}

\proof By definition the derivation $D:\Lambda^0(U,\C) \to
\rho^*\Lambda^1(M,\C)$ is holonomic in the sense of \ref{holonomic}
if and only if the corresponding map 
$$
\Theta:\Lambda^1(U,\R) \to \rho^*\Lambda^1(M,\C)
$$
is an isomorphism of real vector bundles. This is an open condition. 
Therefore the derivation $D$ is holonomic on an open neighborhood
$U_0 \supset M$ of the zero section $M \subset U$ if and only if the
map $\Theta$ is an isomorphism on the zero section $M \subset U$
itself. 

According to \eqref{hor.vert.1}, the restriction $\Theta|_M$
decomposes as $\Theta_M = D_0 + \id$, and the principal part
$D_0:S^1(M,\C) \to \Lambda^1(M,\C)$ of the Hodge connection $D$ is
purely imaginary with respect to the usual real structure on
$\Lambda^1(U,\C)|_M$, while the identity map
$\id:\rho^*\Lambda^1(M,\C) \to \rho^*\Lambda^1(M,\C)$ is, of course,
real. Therefore $\Theta_M$ is an isomorphism if and only is $D_0$ is
an isomorphism, which proves the lemma.  
\endproof

\subsection{Hodge connections and linearity}\label{hodge.lin.subsec} 

\punkt \label{aux} 
Assume now given a Hodge manifold structure on the subset $U \subset
\tm$, and let $D:\Lambda^0(U,\C) \to \rho^*\Lambda^1(M,\C)$ be the
associated Hodge connection on the pair $\langle U, M\rangle$ given
by Proposition~\ref{equiv}. We now proceed to rewrite the linearity
condition \ref{lin.def} in terms of the Hodge connection $D$.

Let $j:\Lambda^1(U,\C) \to \overline{\Lambda^1(U,\C)}$ be the
canonical map defined by the quaternionic structure on $U$, and let
$\iota^*:\Lambda^1(U,\C) \to \iota^*\Lambda^1(U,\C)$ be the action
of the canonical involution $\iota:U \to U$. Let also
$D^\iota:\Lambda^0(U,\C) \to \rho^*\Lambda^1(M,\C)$ be the operator
$\iota^*$-conjugate to the Hodge connection $D$. 

We begin with the following identity.

\begin{lemma}\label{aux.lemma}
For every smooth function $f \in C^\infty(U,\C)$ we have 
$$
d^rf = \halfi\pi(j(\del_\rho(D-D^\iota)(f))),
$$
where $\pi:\Lambda^1(U,\C) \to \Lambda^1(U/M,\C)$ is the canonical
projection, and 
$$
\del_\rho:\rho^*\Lambda^1(M,\C) \to \Lambda^1(U,\C)
$$ 
is the codifferential of the projection $\rho:U \to M$. 
\end{lemma}

\proof By definition of the Hodge connection $D$ the Dolbeault
derivative $\bar\6_Jf$ coincides with the $(0,1)$-component of the
$1$-form $\del_\rho(Df) \in \Lambda^1(U,\C)$ with respect to the
complementary complex structure $U_J$ on $U$. Therefore
$$
\bar\6_J f = \half \del_\rho(Df) + \halfi j(\del_\rho(Df)). 
$$
Applying the complex conjugation $\nu:\Lambda^\cdot(U,\C) \to
\overline{\Lambda^\cdot(U,\C)}$ to this equation, we get
\begin{align*}
\begin{split}
\6_J f &= \nu\left(\half \del_\rho(D\nu(f)) + \halfi
j(\del_\rho(D\nu(f))))\right) = \\
&= \half \nu(\del_\rho(D\nu(f))) - \halfi j(\nu(\del_\rho(D\nu(f)))). 
\end{split}
\end{align*}
Since the map $\del_\rho \circ D:\Lambda^0(U,\C) \to
\rho^*\Lambda^1(M,\C)$ is weakly Hodge, we have
$$
\del_\rho(D(\iota^*\nu(f))) = \iota^*\nu(\del_\rho(Df)).
$$
Therefore $\nu(\del_\rho(D(f))) = \del_\rho(D^\iota(\nu(f)))$, and
we have 
$$
\6_J f = \half \del_\rho(D^\iota f) - \halfi j(\del_\rho(D^\iota
f)). 
$$
Thus the de Rham derivative $d_Uf$ equals 
$$
d_Uf = \6_Jf + \bar\6_Jf = \half \del_\rho((D+D^\iota)f) + \halfi
j(\del_\rho((D-D^\iota)f)). 
$$
Now, by definition $\del_\rho \circ \pi = 0$. Therefore 
$$
d^rf = \pi(d_Uf) = \halfi\pi(j(\del_\rho((D-D^\iota)f))),
$$
which is the claim of the lemma. 
\endproof 

\punkt We will also need the following fact. It can be derived
directly from Lemma~\ref{aux.lemma}, but it is more convenient to
use Lemma~\ref{aux1} and the fact that the Hodge connection
$D:\Lambda^0(U,\C) \to \rho^*\Lambda^1(M,\C)$ is holonomic.

\begin{lemma}\label{aux3}
In the notation of Lemma~\ref{aux.lemma}, let 
$$
\A = \del_\rho\left((D-D^\iota)\left(C_{lin}^\infty(U,\C)\right)\right)
\subset C^\infty(U, \rho^*\Lambda^1(M,\C))
$$ 
be the subspace of sections $\alpha \in
C^\infty(U,\rho^*\Lambda^1(M,\C))$ of the form $\alpha =
\del_\rho((D-D^\iota)f)$, where $f \in C^\infty(U,\C)$ lies in the
subspace $C_{lin}^\infty(U,\C) \subset C^\infty(U,\C)$ of smooth
functions on $U$ linear along the fibers of the projection $\rho:U
\to M$. The restriction $\Res(\A) \subset
C^\infty(M,\Lambda^1(M,\C))$ of the subspace $\A$ onto the zero
section $M \subset U$ is the whole space
$C^\infty(M,\Lambda^1(M,\C))$.
\end{lemma}

\proof Let $D_0 = \Res \circ D:C^\infty_{lin}(U,\C) \to
C^\infty(M,\Lambda^1(M,\C))$ be the principal part of the Hodge
connection $D$ in the sense of Definition~\ref{pr.part}. Since the
canonical automorphism $\iota:\tm \to \tm$ acts as $-\id$ on
$C_{lin}^\infty(U,\C)$, we have $D_0^\iota = - D_0$. Therefore
\begin{multline*}
\Res(\A) = \Res \circ (D - D^\iota)\left(C_{lin}^\infty(U,\C)\right) = \\
=(D_0-D^\iota_0)\left(C_{lin}^\infty(U,\C)\right) =
D_0\left(C_{lin}^\infty(U,\C)\right). 
\end{multline*}
Since the Hodge connection $D$ is holonomic, this space coincides
with the whole $C^\infty(M,\Lambda^1(M,\C))$ by Lemma~\ref{aux1}. 
\endproof 

\punkt We now apply Lemma~\ref{aux.lemma} to prove the following criterion
for the linearity of the Hodge manifold structure on $U$ defined by
the Hodge connection $D:\Lambda^0(U,\C) \to \rho^*\Lambda^1(M,\C)$. 

\begin{lemma}\label{explicit.lin}
The Hodge manifold structure on $U \subset \tm$ corresponding to a
Hodge connection $D:\Lambda^0(U,\C) \to \rho^*\Lambda^1(M,\C)$ is
linear in the sense of \ref{lin.def} if and only if for every smooth
function $f \in C^\infty(U,\C)$ linear along the fibers of the
projection $\rho:U \subset \tm \to M$ we have
\begin{equation}\label{expl.lin}
f = \half\sigma\left((D-D^\iota)f\right), 
\end{equation}
where $\sigma:\rho^*\Lambda^1(M,\C) \to \Lambda^0(U,\C)$ is the
twisted tautological map introduced in \eqref{sigma1}, and
$D^\iota:\Lambda^0(U,\C) \to \rho^*\Lambda^1(M,\C)$ is the operator
$\iota^*$-conjugate to $D$, as in \ref{aux}.
\end{lemma}

\proof By Lemma~\ref{lin.char} the Hodge manifold structure on $U$
is linear if and only if for every $\alpha \in
C^\infty(U,\rho^*\Lambda^1(M,\C))$ we have
\begin{equation}\label{eee.to.prove}
\langle \phi, j(\alpha) \rangle = \tau(\alpha),  
\end{equation}
where $\phi$ is the differential of the $U(1)$-action on $U$,
$j:\Lambda^1(U,\C) \to \overline{\Lambda^1(U,\C)}$ is the operator
given by the quaternionic structure on $U$, and
$\tau:\rho^*\Lambda^1(M,\C) \to \Lambda^0(U,\C)$ is the tautological
map sending a $1$-form on $M$ to the corresponding linear function
on $\tm$, as in \ref{tau}. Moreover, by Lemma~\ref{aux3} and
Lemma~\ref{aux2} the equality \eqref{eee.to.prove} holds for all
smooth sections $\alpha \in C^\infty(U,\rho^*\Lambda^1(M,\C))$ if
and only if it holds for sections of the form
\begin{equation}\label{special}
\alpha = \halfi\del_\rho((D-D^\iota)f), 
\end{equation}
where $f \in C_{lin}^\infty(U,\C) \subset C^\infty(U,\C)$ is linear
along the fibers of $\rho:U \to M$. 

Let now $f \in C^\infty(U,\C)$ be a smooth function on $U$ linear
along the fibers of $\rho:U \to M$, and let $\alpha$ be as in
\eqref{special}. Since $\phi$ is a vertical vector field on $U$ over
$M$, we have $\langle \phi, j(\alpha) \rangle = \langle \phi,
\pi(j(\alpha))\rangle$, where $\pi:\Lambda^1(U,\C) \to
\Lambda^1(U/M,\C)$ is the canonical projection. By
Lemma~\ref{aux.lemma}
\begin{equation}\label{thelema}
\langle \phi, j(\alpha) \rangle = \langle \phi, \pi(j(\alpha))
\rangle = \langle \phi, d^rf \rangle. 
\end{equation}
Since the function $f$ is linear along the fibers of $\rho:U \to M$,
we can assume that $f = \sigma(\beta)$ for a smooth $1$-form $\beta
\in C^\infty(M,\Lambda^1(M,\C)$. Then by \eqref{eta.and.sigma} and
by Lemma~\ref{phi.and.tau} the right hand side of \eqref{thelema} is
equal to
$$
\langle \phi, d^rf \rangle = \langle \phi, d^r(\sigma(\beta))
\rangle = \langle \phi, \eta(\beta) \rangle = \sqrt{-1}
\tau(\beta). 
$$
Therefore, \eqref{eee.to.prove} is equivalent to 
\begin{equation}\label{horus}
\sqrt{-1}\tau(\beta) =
\tau\left(\halfi\del_\rho((D-D^\iota)\sigma(\beta))\right).
\end{equation}
But we have $\tau = \sigma \circ \zeta$, where
$\zeta:\rho^*\Lambda^1(M,\C) \to \rho^*\Lambda^1(M,\C)$ is the
invulution introduced in \eqref{zeta}. In particular, the map
$\zeta$ is invertible, so that \eqref{horus} is in turn equivalent
to
$$
\sigma(\beta) = \half\sigma(\del_\rho((D-D^\iota)\sigma(\beta))),
$$
or, substituting back $f = \sigma(\beta)$, to 
$$
f = \half\sigma(\del_\rho((D-D^\iota)f)),
$$
which is exactly the condition \eqref{expl.lin}. 
\endproof 

\punkt \defn \label{hodge.conn.lin}
A Hodge connection $D$ on the pair $\langle U, M \rangle$ is called
{\em linear} if it satisfies the condition~\ref{expl.lin}.

We can now formulate and prove the following more useful version of
Proposition~\ref{equiv}.

\begin{prop}\label{equiv.bis}
Every linear Hodge connection $D\!:\!\Lambda^0(U,\C) \to
\rho^*\Lambda^1(M,\C)$ on the pair $\langle U, M\rangle$ defines a
linear Hodge manifold structure on an open neighborhood $V \subset
U$ of the zero section $M \subset U$, and the canonical projection
$\rho:V_I \to M$ is holomorphic for the preferred complex structure
$V_I$ on $V$. Vice versa, every such linear Hodge manifold structure
on $U$ comes from a unique linear Hodge connection $D$ on the pair
$\langle U,M \rangle$.
\end{prop}

\proof By Proposition~\ref{equiv} and Lemma~\ref{explicit.lin}, to
prove this proposition suffices to prove that if a Hodge connection
$D:\Lambda^0(U,\C) \to \rho^*\Lambda^1(M,\C)$ is linear, then it is
holonomic in the sense of \ref{holonomic} on a open neighborhood $V
\subset U$ of the zero section $M \subset U$. Lemma~\ref{aux1}
reduces this to proving that the principal part $D_0:S^1(M,\C) \to
\Lambda^1(M,\C)$ of a linear Hodge connection $D:\Lambda^0(U,\C) \to
\rho^*\Lambda^1(M,\C)$ is a bundle isomorphism. 

Let $D:\Lambda^0(U,\C) \to \rho^*\Lambda^1(M,\C)$ be such
connection. By \eqref{expl.lin} we have
$$
\half \sigma \circ (D_0 - D_0^\iota) = \id:S^1(M,\C) \to
\Lambda^1(M,\C) \to S^1(M,\C).
$$
Since $\sigma:\Lambda^1(M,\C) \to S^1(M,\C)$ is a bundle
isomorphism, so is the bundle map $D_0 - D_0^\iota:S^1(M,\C) \to
\Lambda^1(M,\C)$. As in the proof of Lemma~\ref{aux3}, we have $D_0
= - D_0^\iota$. Thus $D_0 = \half(D_0-D_0^\iota):S^1(M,\C) \to
\Lambda^1(M,\C)$ also is a bundle isomorphism, which proves the
proposition.  
\endproof

\section{Formal completions}\label{formal.section}

\subsection{Formal Hodge manifolds}

\punkt Proposition~\ref{equiv.bis} reduces the study of arbitrary
regular Hodge manifolds to the study of connections of a certain
type on a neighborhood $U \subset \tm$ of the zero section $M
\subset \tm$ in the total space $\tm$ of the tangent bundle to a
complex manifold $M$.  To obtain further information we will now
restrict our attention to the {\em formal} neighborhood of this zero
section. This section contains the appropriate definitions. We study
the convergence of our formal series in Section~\ref{convergence}.

\punkt Let $X$ be a smooth manifold and let $\Bun(X)$ be the category of
smooth real vector bundles over $X$. Let also $\Diff(X)$ be the
category with the same objects as $\Bun(X)$ but with differential operators
as morphisms.

Consider a closed submanifold $Z \subset X$. For every two real vector
bundles $\E$ and $\F$ on $X$ the vector space $\Hom(\E,\F)$ of bundle maps
from $\E$ to $\F$ is naturally a module over the ring $C^\infty(X)$ of
smooth functions on $X$.  Let $\J_Z \subset C^\infty(X)$ be the ideal of
functions that vanish on $Z$ and let $\Hom_Z(\E,\F)$ be the $\J_Z$-adic
completion of the $C^\infty(X)$-module $\Hom(\E,\F)$.

For any three bundles $\E,\F,\G$ the composition map 
$$
\Mult:\Hom(\E,\F) \otimes \Hom(\F,\G) \to \Hom(\E,\G)
$$ 
is $C^\infty(X)$-linear, hence extends to a map
$$
\Mult:\Hom_Z(\E,\F) \otimes \Hom_Z(\F,\G) \to \Hom_Z(\E,\G). 
$$
Let $\Bun_Z(X)$ be the category with the same objects as $\Bun(X)$
and for every two objects $\E$, $\F \in \Ob\Bun(X)$ with
$\Hom_Z(\E,\F)$ as the space of maps between $\F$ anf $\F$.  The
category $\Bun_Z(X)$, as well as $\Bun(X)$, is an additive tensor
category.

\punkt The space of differential operators $\Diff(\E,\F)$ is also a
$C^\infty(X)$ module, say, by left multiplication. Let $\Diff_Z(\E,\F)$ be
its $\J_Z$-completion. The composition maps in $\Diff(X)$ are no longer
$C^\infty(X)$-linear. However, they still are compatible with the
$\J_Z$-adic topology, hence extend to completions. Let $\Diff_Z(X)$ be the
category with the same objects as $\Bun(X)$ and with $\Diff_Z(\E,\F)$ as
the space of maps between two objects $\E,\F \in \Ob \Bun(X)$.

By construction we have canonical {\em $Z$-adic completion functors} 
$$
\Bun(X) \to \Bun_Z(X) \text{  and  } \Diff(X) \to \Diff_Z(X). 
$$
Call the categories $\Bun_Z(X)$ and $\Diff_Z(X)$ {\em the $Z$-adic
completions} of the categories $\Bun(X)$ and $\Diff(X)$.

\punkt When the manifold $X$ is equipped with a smooth action of
compact Lie group $G$ fixing the submanifold $Z$, the completion
construction extends to the categories of $G$-equivariant bundles on
$M$.  When $G = U(1)$, the categories $\prehodge(X)$ and
$\prehodge^\D(X)$ defined in \ref{w.hodge} also admit canonical
completions, denoted by $\prehodge_Z(X)$ and $\prehodge_Z^\D(X)$.

\punkt Assume now that the manifold $X$ is equipped with a smooth
$U(1)$-action fixing the smooth submanifold $Z \subset X$. 

\defn A {\em formal quaternionic structure} on $X$ along the
submanifold $Z \subset X$ is given by an algebra map
$$
\Mult: \h \to \Eend_{\Bun_Z(X)} \left(\Lambda^1(X)\right)
$$
from the algebra $\h$ to the algebra $\Eend_{\Bun_Z(X)}
\left(\Lambda^1(X)\right)$ of endomorphisms of the cotangent bundle
$\Lambda^1(X)$ in the category $\Bun_Z(X)$.
A formal quaternionic structure is called {\em equivariant} if the
map $\Mult$ is equivariant with respect to the natural $U(1)$-action
on both sides. 

\punkt Note that Lemma~\ref{universal} still holds in the situation
of formal completions. Consequently, everything in
Section~\ref{hbqm.section} carries over word-by-word to the case of
formal quaternionic structures. In particular, by
Lemma~\ref{qm.hodge} giving a formal equivariant quaternionic
structure on $X$ along $Z$ is equivalent to giving a pair $\langle
\E, D \rangle$ of a Hodge bundle $\E$ on $X$ and a holonomic algebra
derivation $D:\Lambda^0(X) \to \E$ in $\prehodge_Z^\D(X)$.

\punkt The most convenient way to define Hodge manifold structures
on $X$ in a formal neighborhood of $Z$ is by means of the Dolbeault
complex, as in Proposition~\ref{explicit.hodge}.

\defn A {\em formal Hodge manifold structure} on $X$ along $Z$ is a pair of
a pre-Hodge bundle $\E \in \Ob \prehodge_Z(X)$ of weight $1$ and an
algebra derivation $D^\cdot:\Lambda^\cdot\E \to \Lambda^{\cdot+1}\E$ in
$\prehodge^\D_Z(X)$ such that $D^0:\Lambda^0(\E) \to \E$ is
holonomic and $D^0 \circ D^1 = 0$. 

\punkt Let $U \subset X$ be an open subset containing $Z \subset
X$. For every Hodge manifold structure on $U$ the $Z$-adic
completion functor defines a formal Hodge manifold structure on $X$
along $Z$. Call it {\em the $Z$-adic completion} of the given
structure on $U$.

\rem Note that a Hodge manifold structure on $U$ is completely
defined by the preferred and the complementary complex structures
$U_I$, $U_J$, hence always real-analytic by the Newlander-Nirenberg
Theorem.  Therefore, if two Hodge manifold structures on $U$ have
the same completion, they coincide on every connected component of
$U$ intersecting $Z$.

\subsection{Formal Hodge manifold structures on tangent bundles}

\punkt Let now $M$ be a complex manifold, and let $\tm$ be the total
space of the complex-conjugate to the tangent bundle to $M$ equipped
with an action of $U(1)$by dilatation along the fibers of the
projection $\rho:\tm \to M$. All the discussion above applies to the
case $X = \tm$, $Z = M \subset \tm$. Moreover, the linearity
condition in the form given in Lemma~\ref{lin.char} generalizes
immediately to the formal case.

\defn A formal Hodge manifold structure on $\tm$ along $M$ is called
{\em linear} if for every smooth $(0,1)$-form $\alpha \in
C^\infty(M,\Lambda^{0,1}(M))$ we have 
$$
\tau(\alpha) = \langle \phi, j(\rho^*) \rangle \in
C^\infty_M(\tm,\C), 
$$
where $j$ is the map induced by the formal quaternionic structure on
$\tm$ and $\phi$ and $\tau$ are as in Lemma~\ref{lin.char}. 

\punkt As in the non-formal case, linear Hodge manifold structures
on $\tm$ along $M \subset \tm$ can be described in terms of
differential operators of certain type. 

\defn \label{formal.hodge.con}
A {\em formal Hodge connection} on $\tm$ along $M \subset \tm$ is an
algebra derivation
$$
D:\Lambda^0(\tm,\C) \to \rho^*\Lambda^1(M,\C)
$$
in $\prehodge_M^\D(\tm)$ such that for every smooth function $f \in
C^\infty(M,\C)$ we have $D\rho^*=\rho^*d_Mf$, as in
\eqref{conn.eq}. A formal Hodge connection is called {\em flat} if
it extends to an algebra derivation
$$
D:\rho^*\Lambda^\cdot(M,\C) \to \Lambda^{\cdot+1}(M,\C) 
$$
in $\prehodge_M^\D(\tm)$ such that $D \circ D = 0$. A formal Hodge
connection is called {\em linear} if it satisfies the condition
\eqref{expl.lin} of Lemma~\ref{explicit.lin}, that is, for every
function $f \in C^\infty_{lin}(\tm,\C)$ linear along the fibers of
the projection $\rho:\tm \to M$ we have 
$$
f = \half\sigma\left((D-D^\iota)f\right), 
$$
where $\sigma:\rho^*\Lambda^1(M,\C) \to \Lambda^0(\tm,\C)$ is the
twisted tautological map introduced in \eqref{sigma1}, the
automorphism $\iota:\tm \to \tm$ is the multiplication by $-1 \in
\C$ on every fiber of the projection $\rho:\tm \to M$, and
$D^\iota:\Lambda^0(\tm,\C) \to \rho^*\Lambda^1(M,\C)$ is the
operator $\iota^*$-conjugate to $D$, as in \ref{aux}.

The discussion in Section~\ref{section.5} generalizes immediately to
the formal case and gives the following. 

\begin{lemma}\label{form.hdg} 
Linear formal Hodge manifold structures on $\tm$ along the zero
section $M \subset \tm$ are in a natural one-to-one correspondence
with linear flat formal Hodge connections on $\tm$ along $M$.
\end{lemma}

\subsection{The Weil algebra}

\punkt Let, as before, $M$ be a complex manifold and let $\tm$ be
the total space of the complex conjugate to its tangent bundle, as
in \ref{overline.T}. In the remaining part of this section we give
a description of the set of all formal Hodge connections on $\tm$
along $M$ in terms of certain differential operators on $M$ rather
than on $\tm$. We call such operators {\em extended connections} on
$M$ (see \ref{ext.con} for the definition). Together with a
complete classification of extended connections given in the next
Section, this description provides a full classification of regular
Hodge manifolds ``in the formal neighborhood of the subset of
$U(1)$-fixed points''. 

\punkt Before we define extended connections in
Subsection~\ref{ext.con.subsec}), we need to introduce a certain
algebra bundle in $\prehodge(M)$ which we call {\em the Weil
algebra}.  We begin with some preliminary facts.

Recall (see, e.g., \cite{Del}) that every additive category $\A$
admits a canonical completion $\llim\A$ with respect to filtered
projective limits. The category $\llim\A$ is also additive, and it
is tensor if $\A$ was tensor.  Objects of the canonical completion
$\llim\A$ are called {\em pro-objects in $\A$}.

\punkt\label{rho_*} 
Let $\rho:\tm \to M$ be the canonical projection. Extend the
pullback functor $\rho^*:\Bun(M) \to \Bun(\tm)$ to a functor
$$
\rho^*:\Bun(M) \to \Bun_M(\tm) 
$$
to the $M$-adic completion $\Bun_M(\tm)$. The functor $\rho^*$
admits a right adjoint direct image functor 
$$
\rho_*:\Bun_M(\tm) \to \llim\Bun(M).
$$
Moreover, the functor $\rho_*$ extends to a functor
$$
\rho_*:\Diff_M(\tm) \to \llim\Diff(M). 
$$
Denote by $\B^0(M,\C) = \rho_*\Lambda^0(\tm)$ the direct image under
the projection $\rho:\tm \to M$ of the trivial bundle
$\Lambda^0(\tm)$ on $\tm$.

The compact Lie group $U(1)$ acts on $\tm$ by dilatation along the
fibers, and the functor $\rho_*:\Diff_M(\tm) \to \llim\Diff(M)$
obviously extends to a functor $\rho_*:\prehodge^\D_M(\tm) \to
\llim\prehodge(M)$. The restriction of $\rho_*$ to the subcategory
$\prehodge_M(\tm) \subset \prehodge^\D_M(\tm)$ is adjoint on the
right to the pullback functor $\rho^*:\prehodge(M) \to
\prehodge_M(\tm)$.

\punkt The constant bundle $\Lambda^0(\tm)$ is canonically a Hodge
bundle of weight $0$. Therefore $\B^0(M,\C) = \rho_*\Lambda^0(M,\C)$
is also a Hodge bundle of weight $0$.  Moreover, it is a commutative
algebra bundle in $\llim\prehodge_0(M)$. Let $S^1(M,\C)$ be the
conormal bundle to the zero section $M \subset \tm$ equipped with a
Hodge bundle structure of weight $0$ as in \ref{S1}, and denote by
$S^i(M,\C)$ the $i$-th symmetric power of the Hodge bundle
$S^1(M,\C)$. Then the algebra bundle $\B^0(M,\C)$ in
$\llim\prehodge_0(M)$ is canonically isomorphic
$$
\B^0(M,\C) \cong \widehat{S}^\cdot(M,\C)
$$
to the completion $\widehat{S}^\cdot(M,\C)$ of the symmetric algebra
$S^\cdot(M,\C)$ of the Hodge bundle $S^1(M,\C)$ with respect to the
augmentation ideal $S^{>0}(M,\C)$. 

Since the $U(1)$-action on $M$ is trivial, the category
$\prehodge(M)$ of Hodge bundles on $M$ is equivalent to the category
of pairs $\langle \E, \conj \rangle$ of a complex bundle $\E$
equipped with a Hodge type grading
$$
\E = \bigoplus_{p,q} \E^{p,q}
$$
and a real structure $\conj:\E^{p,q} \to \overline{\E^{q,p}}$. The
Hodge type grading on $\B^0(M,\C)$ is induced by the Hodge type
grading $S^1(M,\C) = S^{1,-1}(M) \oplus S^{-1,1}(M)$ on the
generators subbundle $S^1(M,\C) \subset \B^0(M,\C)$, which was
described in \ref{S1}.

\rem The complex vector bundle $S^1(M,\C)$ is canonically isomorphic
to the cotangent bundle $\Lambda^1(M,\C)$. However, the Hodge bundle
structures on these two bundles are different (in fact, they have
different weights).

\punkt \label{Weil.defn}
Consider the pro-bundles  
$$
\B^\cdot(M,\C) = \rho_*\rho^*\Lambda^\cdot(M,\C)
$$
on $M$. The direct sum $\oplus\B^\cdot(M,\C)$ is a graded algebra in
$\llim\Bun(M,\C)$. Moreover, since for every $i \geq 0$ the bundle
$\Lambda^i(M,\C)$ is a Hodge bundle of weight $i$ (see
\ref{de.Rham}), $\B^i(M,\C)$ is also a Hodge bundle of weight
$i$. Denote by
$$
\B^i(M,\C) = \bigoplus_{p+q=i} \B^{p,q}(M,\C) 
$$
the Hodge type bigrading on $\B^i(M,\C)$. 

The Hodge bundle structures on $\B^\cdot(M,\C)$ are compatible with
the multiplication. By the projection formula
$$
\B^\cdot(M,\C) \cong \B^0(M,\C) \otimes \Lambda^\cdot(M,\C), 
$$
and this isomorphism is compatible with the Hodge bundle structures on
both sides. 

\defn Call the algebra $\B^\cdot(M,\C)$ in $\llim\prehodge(M)$ {\em
the Weil algebra} of the complex manifold $M$. 

\punkt \label{iota.Weil} 
The canonical involution $\iota:\tm \to \tm$ induces an algebra
involution $\iota^*:\B^\cdot(M,\C) \to \B^\cdot(M,\C)$. It acts on
generators as follows
$$
\iota^* = -\id:S^1(M,\C) \to S^1(M,\C) \qquad 
\iota^* =  \id:\Lambda^1(M,\C) \to \Lambda^1(M,\C).  
$$
For every operator $N:\B^p(M,\C) \to \B^q(M,\C)$, $p$ and $q$
arbitrary, we will denote by
$$
N^\iota = \iota^* \circ N \circ \iota^*:\B^p(M,\C) \to
\B^q(M,\C)
$$ 
the operator $\iota^*$-conjugate to $N$. 

\punkt \label{sigma}
The twisted tautological map $\sigma:\rho^*\Lambda^1(M,\C) \to
\Lambda^0(\tm,\C)$ introduced in \ref{sigma1} induces via the
functor $\rho_*$ a map $\sigma:\B^1(M,\C) \to \B^0(M,\C)$. Extend
this map to a derivation
$$
\sigma:\B^{\cdot+1}(M,\C) \to \B^\cdot(M,\C)
$$
by setting $\sigma = 0$ on $S^1(M,\C) \subset \B^0(M,\C)$.  The
derivation $\B^{\cdot+1}(M,\C) \to \B^\cdot(M,\C)$ is not weakly
Hodge. However, it is real with respect to the real structure on the
Weil algebra $\B^\cdot(M,\C)$.

\punkt \label{C.Weil} 
By definition of the twisted tautological map (\ref{sigma1},
\ref{tau}), the derivation $\sigma:\B^{\cdot+1}(M,\C) \to
\B^\cdot(M,\C)$ maps the subbundle $\Lambda^1(M,\C) \subset
\B^1(M,\C)$ to the subbundle $S^1(M,\C) \subset \B^0(M,\C)$ and
defines a complex vector bundle isomorphism $\sigma:\Lambda^1(M,\C)
\to S^1(M,\C)$. To describe this isomorphism explicitly, recall that
sections of the bundle $\B^0(M,\C)$ are the same as formal germs
along $M \subset \tm$ of smooth functions on the manifold $\tm$. The
sections of the subbundle $S^1(M,\C) \subset \B^0(M,\C)$ form the
subspace of functions linear along the fibers of the canonical
projection $\rho:\tm \to M$. The isomorphism $\sigma:\Lambda^1(M,\C)
\to S^1(M,\C)$ induces an isomorphism between the space of smooth
$1$-forms on the manifold $M$ and the space of smooth functions on
$\tm$ linear long the fibers of $\rho:\tm \to M$. This isomorphism
coincides with the tautological one on the subbundle $\Lambda^{1,0}
\subset \Lambda^1(M,\C)$, and it is minus the tautological
isomorphism on the subbundle $\Lambda^{0,1} \subset
\Lambda^1(M,\C)$.

Denote by 
$$
C = \sigma^{-1}:S^1(M,\C) \to \Lambda^1(M,\C)
$$ 
the bundle isomorphism inverse to $\sigma$. Note that the complex
vector bundle isomorphism $\sigma:\Lambda^1(M,\C) \to S^1(M,\C)$ is
real. Moreover, it sends the subbundle $\Lambda^{1,0}(M) \subset
\Lambda^1(M,\C)$ to $S^{1,-1}(M) \subset S^1(M,\C)$, and it sends
$\Lambda^{0,1}(M)$ to $S^{-1,1}(M)$. Therefore the inverse
isomorphism $C:S^1(M,\C) \to \Lambda^1(M,\C)$ is weakly Hodge. It
coincides with the tautological isomorphism on the subbundle
$S^{1,-1} \subset S^1(M,\C)$, and it equals minus the tautological
isomorphism on the subbundle $S^{-1,1} \subset S^1(M,\C)$.

\subsection{Extended connections}\label{ext.con.subsec}

\punkt We are now ready to introduce the extended connections. Keep
the notation of the last subsection. 

\defn \label{ext.con} 
An {\em extended connection} on a complex manifold $M$ is a
differential operator $D:S^1(M,\C) \to \B^1(M,\C)$ which is weakly
Hodge in the sense of \ref{w.hodge} and satisfies
\begin{equation}\label{e.c}
D(fa) = fDa + a \otimes df
\end{equation}
for any smooth function $f$ and a local section $a$ of the pro-bundle
$\B^0(M,\C)$. 

\punkt \label{red}
Let $D$ be an extended connection on the manifold $M$. 
By \ref{Weil.defn} we have canonical bundle isomorphisms 
$$
\B^1(M,\C) \cong \B^0(M,\C) \otimes \Lambda^1(M,\C) \cong \bigoplus_{i
\geq 0} S^i(M,\C) \otimes \Lambda^1(M,\C). 
$$
Therefore the operator $D:S^1 \to \B^1$ decomposes 
\begin{equation}\label{aug.con}
D = \sum_{p \geq 0}D_p, \quad D_p:S^1(M,\C) \to S^i(M,\C) \otimes
\Lambda^1(M,\C). 
\end{equation}
By \eqref{e.c} all the components $D_p$ except for the $D_1$ are
weakly Hodge bundle maps on $M$, while 
$$
D_1:S^1(M,\C) \to S^1(M,\C) \otimes \Lambda^1(M,\C) 
$$
is a connection in the usual sense on the Hodge bundle $S^1(M,\C)$. 

\defn The weakly Hodge bundle map $D_0:S^1(M,\C) \to
\Lambda^1(M,\C)$ is called {\em the principal part} of the extended
connection $D$ on $M$. The connection $D_1$ is called {\em the
reduction} of the extended connection $D$.

\punkt Extended connection on $M$ are related to formal Hodge
connections on the total space $\tm$ of the complex-conjugate to the
tangent bundle to $M$ by means of the direct image functor
$$
\rho_*:\prehodge^\D_M(\tm) \to \llim\prehodge^\D(M).
$$ 
Namely, let $D:\Lambda^0(M,\C) \to \rho^*\Lambda^1(M,\C)$ be a
formal Hodge connection on $\tm$ along $M$ in the sense of
\ref{formal.hodge.con}. The restriction of the operator
$$
\rho_*D:\B^0(M,\C) \to \B^1(M,\C)
$$
to the generators subbundle $S^1(M,\C) \subset \B^0(M,\C)$ is then
an extended connection on $M$ in the sense of \ref{ext.con}. The
principal part $D_0:S^1(M,\C) \to \Lambda^1(M,\C)$ of the Hodge
connection $D$ in the sense of \ref{pr.part} coincides with the
principal part of the extended connection $\rho_*D$. 

\punkt We now write down the counterparts of the flatness and
linearity conditions on a Hodge connection on $\tm$ for the
associated extended connection on $M$. We begin with the linearity
condition \ref{hodge.conn.lin}. Let $D:S^1(M,\C) \to \B^1(M,\C)$ be
an extended connection on $M$, let $\sigma:\B^{\cdot+1}(M,\C) \to
\B^\cdot(M,\C)$ be the algebra derivation introduced in \ref{sigma},
and let
$$
D^\iota:\B^\cdot(M,\C) \to \B^{\cdot+1}(M,\C) 
$$
be the operator $\iota^*$-conjugate to $D$ as in \ref{iota.Weil}. 

\defn \label{lin.ext.con}
An extended connection $D$ is called {\em linear} if for
every local section $f$ of the bundle $S^1(M,\C)$ 
we have 
$$
f = \half \sigma((D-D^\iota)f).
$$
This is, of course, the literal rewriting of
Definition~\ref{hodge.conn.lin}. In particular, a formal Hodge
connection $D$ on $\tm$ is linear if and only if so is the extended
connection $\rho_*D$ on $M$.

\punkt \label{deriv} 
Next we rewrite the flatness condition \ref{hodge.con}. Again, let
$$
D:S^1(M,\C) \to \B^1(M,\C)
$$ 
be an extended connection on $M$.  Since the algebra pro-bundle
$\B^0(M,\C)$ is freely generated by the subbundle $S^1(M,\C) \subset
\B^1(M,\C)$, by \eqref{e.c} the operator $D:S^1(M,\C) \to
\B^1(M,\C)$ extends uniquely to an algebra derivation
$$
D:\B^0(M,\C) \to \B^1(M,\C). 
$$
Moreover, we can extend this derivation even further to a
derivation of the Weil algebra 
$$
D:\B^\cdot(M,\C) \to \B^{\cdot+1}(M,\C) 
$$
by setting 
\begin{equation}\label{D=d}
D(f \otimes \alpha) = Df \otimes \alpha + f \otimes d\alpha
\end{equation} 
for any smooth section $f \in C^\infty(M,\B^0(M,\C))$ and any smooth
form $\alpha \in C^\infty(M,\Lambda^\cdot(M,\C))$. We will call this
extension {\em the derivation, associated to the extended connection
$D$}.

Vice versa, the Weil algebra $\B^\cdot(M,\C)$ is generated by the
subbundles
$$
S^1(M,\C),\Lambda^1(M,\C) \subset \B^\cdot(M,\C). 
$$
Moreover, for every algebra derivation $D:\B^\cdot(M,\C) \to
\B^{\cdot+1}(M,\C)$ the condition \eqref{D=d} completely defines the
restriction of $D$ to the generator subbundle $\Lambda^1(M,\C)
\subset \B^1(M,\C)$. Therefore an algebra derivation
$D:\B^\cdot(M,\C) \to \B^\cdot(M,\C)$ satisfying \eqref{D=d} is
completely determined by its restriction to the generators subbundle
$S^1(M,\C) \to \B^1(M,\C)$. If the derivation $D$ is weakly Hodge,
then this restriction is an extended connection on $M$.
 
\punkt \defn \label{flat.ext.con}
The extended connection $D$ is called {\em flat} if the associated
derivation satisfies $D \circ D = 0$.

If a formal Hodge connection $D$ on $\tm$ is flat in the sense of
\ref{hodge.con}, then we have a derivation
$D:\rho^*\Lambda^\cdot(M,\C) \to \rho^*\Lambda^{\cdot+1}(M,\C)$ such
that $D \circ D = 0$. The associated derivation
$\rho_*D:\B^\cdot(M,\C) \to \B^{\cdot+1}(M,\C)$ satisfies
\eqref{D=d}. Therefore the extended connection $\rho_*D:S^1(M,\C)
\to \B^1(M,\C)$ is also flat.

\punkt It turns out that one can completely recover a Hodge
connection $D$ on $\tm$ from the corresponding extended connection
$\rho_*D$ on $M$. More precisely, we have the following. 

\begin{lemma}
The correspondence $D \mapsto \rho_*D$ is a bijection between the
set of formal Hodge connections on $\tm$ along $M \subset \tm$ and
the set of extended connections on $M$. A connection $D$ is flat,
resp. linear if and only if $\rho_*D$ is flat,
resp. linear. 
\end{lemma}

\proof To prove the first claim of the lemma, it suffices to prove
that every extended connection on $M$ comes from a unique formal
Hodge connection on the pair $\langle \tm,M\rangle$. In general, the
functor $\rho_*$ is not fully faithful on the category $\Diff(M)$,
in other words, it does not induce an isomorphism on the spaces of
differential operators between vector bundles on $\tm$. However, for
every complex vector bundle $\E$ on $\tm$ the functor $\rho_*$ does
induce an isomorphism
$$
\rho_*:\Der_M(\Lambda^0(M,\C),\E) \cong
\Der_{\B^0(M,\C)}(\B^0(M,\C),\rho_*\E) 
$$
between the space of {\em derivations} from $\Lambda^0(M,\C)$ to
$\F$ completed along $M \subset \tm$ and the space of derivations
from the algebra $\B^0(M,\C) = \rho_*\Lambda^0(M,\C)$ to the
$\B^0(M,\C)$-module $\rho_*\E$. Therefore every derivation
$$
D':\B^0(M,\C) \to \B^1(M,\C) = \rho_*\rho^*\Lambda^1(M,\C)
$$ 
must be of the form $D'=\rho_*D$ for some derivation
$$
D:\Lambda^0(\tm,\C) \to \rho^*\Lambda^1(M,\C)
$$ 
It is easy to check that $D$ is a Hodge connection if and only if
$D'=\rho_*D$ is weakly Hodge and satisfies \eqref{e.c}. By
\ref{deriv} the space of all such derivations $D':\B^0(M,\C) \to
\B^1(M,\C)$ coincides with the space of all extended connections on
$M$, which proves the first claim of the lemma.

Analogously, for every extended connection $D' = \rho_*D$ on $M$,
the canonical extension of the operator $D'$ to an algebra
derivation $D':\B^\cdot(M,\C) \to \B^{\cdot+1}(M,\C)$ constructed in
\ref{deriv} must be of the form $\rho_*D$ for a certain weakly Hodge
differential operator $D:\rho^*\Lambda^\cdot(M,\C) \to
\Lambda^{\cdot+1}(M,\C)$. If the extended connection $D'$ is flat,
then $D' \circ D' = 0$. Therefore $D \circ D = 0$, which means that
the Hodge connection $D$ is flat. Vice versa, if the Hodge
connection $D$ is flat, then it extends to a weakly Hodge derivation
$D:\rho^*\Lambda^\cdot(M,\C) \to \Lambda^{\cdot+1}(M,\C)$ so that $D
\circ D = 0$. The equality $D \circ D = 0$ implies, in particular,
that the operator $\rho_*D$ vanishes on the sections of the form 
$$
Df = df \in C^\infty(M,\Lambda^1(M,\C)) \subset
C^\infty(M,\B^1(M,\C)), 
$$
where $f \in C^\infty(M,\C)$ is a smooth function on $M$. Therefore
$\rho_*D$ coincides with the de Rham differential on the subbundle
$\Lambda^1(M,\C) \subset \B^1(M,\C)$. Hence by \ref{deriv} it is
equal to the canonical derivation $D':\B^\cdot(M,\C) \to
\B^{\cdot+1}(M,\C)$. Since $D \circ D = 0$, we have $D' \circ D' =
0$, which means that the extended connection $D'$ is flat.

Finally, the equivalence of the linearity conditions on the Hodge
connection $D$ and on the extended connection $D' = \rho_*D$ is
trivial and has already been noted in \ref{lin.ext.con}. 
\endproof 

This lemma together with Lemma~\ref{form.hdg} reduces the
classification of linear formal Hodge manifold structures on $\tm$
along the zero section $M \subset \tm$ to the classification of
extended connections on the manifold $M$ itself.

\section{Preliminaries on the Weil algebra}\label{Weil.section}

\subsection{The total de Rham complex}\label{de.rham.sub} 

\punkt Before we proceed further in the study of extended
connections on a complex manifold $M$, we need to establish some
linear-algebraic facts on the structure of the Weil algebra
$\B^\cdot(M,\C)$ defined in \ref{Weil.defn}. We also need to
introduce an auxiliary Hodge bundle algebra on $M$ which we call
{\em the total Weil algebra}. This is the subject of this
section. Most of the facts here are of a technical nature, and the
reader is advised to skip this section until needed.

\punkt We begin with introducing and studying a version of the de
Rham complex of a complex manifold $M$ which we call {\em the total
de Rham complex}. Let $M$ be a smooth complex
$U(1)$-manifold. Recall that by \ref{de.Rham} the de Rham complex
$\Lambda^\cdot(M,\C)$ of the complex manifold $M$ is canonically a
Hodge bundle algebra on $M$.  Let $\Lambda_{tot}^\cdot(M) =
\Gamma(\Lambda^\cdot(M,\C))$ be the weight $0$ Hodge bundle obtained
by applying the functor $\Gamma$ defined in \ref{gamma.m} to the de
Rham algebra $\Lambda^\cdot(M,\C)$. By \ref{gamma.tensor} the bundle
$\Lambda_{tot}^\cdot(M)$ carries a canonical algebra structure. By
\ref{de.Rham} the de Rham differential $d_M$ is weakly
Hodge. Therefore it induces an algebra derivation
$d_M:\Lambda^\cdot_{tot}(M) \to \Lambda^{\cdot+1}_{tot}(M)$ which is
compatible with the Hodge bundle structure and satisfies $d_M \circ
d_M = 0$.

\defn The weight $0$ Hodge bundle algebra $\Lambda_{tot}^\cdot(M)$ is
called {\em the total de Rham complex} of the complex manifold $M$. 

\punkt By definition
$$
\Lambda^i_{tot}(M) = \Gamma(\Lambda^i(M,\C)) = \Lambda^i(M,\C) \otimes \W^*_i,  
$$
where $\W^*_i = S^i\W^*_1$ is the symmetric power of the $\R$-Hodge
structure $\W^*_1$, as in \ref{w.k}. To describe the structure of
the algebra $\Lambda^\cdot_{tot}(M)$, we will use the following
well-known general fact. (For the sake of completeness, we have
included a sketch of its proof, see \ref{gen.symm.proof}.)

\begin{lemma}\label{symm}
Let $A$, $B$ be two objects in an arbitrary $\Q$-linear symmetric
tensor category $\A$, and let $\CC^\cdot = S^\cdot(A \otimes B)$ be
the sum of symmetric powers of the object $A \otimes B$. Note that
the object $\CC^\cdot$ is naturally a commutative algebra in $\A$ in
the obvious sense. Let also $\wt{\CC}^\cdot = \bigoplus_k S^kA
\otimes S^kB$ with the obvious commutative algebra structure. The
isomorphism $\CC^1 \cong \wt{\CC}^1 \cong A \otimes B$ extends to a
surjective algebra map $\CC^\cdot \to \wt{\CC}^\cdot$, and its
kernel $\J^\cdot \subset \CC^\cdot$ is the ideal generated by the
subobject $\J^2 = \Lambda^2(A) \otimes \Lambda^2(B) \subset S^2(A
\otimes B)$.
\end{lemma}

\punkt The category of complexes of Hodge bundles on $M$ is
obviously $\Q$-linear and tensor. Applying
Lemma~\ref{symm} to $A = \W^*_1$, $B = \Lambda^1(M,\C)
[1]$ immediately gives the following. 
 
\begin{lemma}\label{total.rel}
The total de Rham complex $\Lambda_{tot}^\cdot(M)$ of the complex
manifold $M$ is generated by its first component $\Lambda^1_{tot}(M)$,
and the kernel of the canonical surjective algebra map
$$
\Lambda^\cdot(\Lambda_{tot}^1(M)) \to \Lambda_{tot}^\cdot(M) 
$$
from the exterior algebra of the bundle $\Lambda^1_{tot}(M)$ to
$\Lambda^\cdot_{tot}(M)$ is the ideal generated by the subbundle
$$
\Lambda^2\W_1 \otimes S^2(\Lambda^1(M,\C)) \subset S^2(\Lambda^1_{tot}(M)).
$$
\end{lemma}

\punkt \label{S}
We can describe the Hodge bundle $\Lambda^1_{tot}(M)$ more explicitly
in the following way. By definition, as a $U(1)$-equivariant complex
vector bundle it equals
$$
\Lambda^1_{tot}(M) = \Lambda^1(M,\C) \otimes \W_1^* = \left(
\Lambda^{1,0}(M)(1) \oplus \Lambda^{0,1}(M)(0)\right) \otimes \left(\C(0)
\oplus \C(-1)\right),
$$
where $\Lambda^{p,q}(M)(i)$ is the $U(1)$-equivariant bundle
$\Lambda^{p,q}(M)$ tensored with the $1$-dimensional representation
of weight $i$, and $\C(i)$ is the constant $U(1)$-bundle
corresponding to the representation of weight $i$. If we denote
\begin{align*}
S^1(M,\C) &= \Lambda^{1,0}(M)(1) \oplus \Lambda^{0,1}(M)(-1) \subset
\Lambda^1_{tot}(M),\\ 
\Lambda^1_{ll}(M) &= \Lambda^{1,0}(M) \subset \Lambda^1_{tot}(M),\\
\Lambda^1_{rr}(M) &= \Lambda^{0,1}(M) \subset \Lambda^1_{tot}(M),
\end{align*}
then we have 
$$
\Lambda^1_{tot}(M) = S^1(M,\C) \oplus \Lambda^1_{ll}(M) \oplus
\Lambda^1_{rr}(M).
$$
The complex conjugation $\conj:\Lambda^1_{tot}(M) \to
\iota^*\overline{\Lambda^1_{tot}(M)}$ preserves the subbundle 
$$
S^1(M,\C) \subset \Lambda^1_{tot}(M,\C)
$$ 
and interchanges $\Lambda^1_{ll}(M)$ and $\Lambda^1_{rr}(M)$.

\punkt\label{S.Hodge.type}
If the $U(1)$-action on the manifold $M$ is trivial, then Hodge
bundles are the same as bigraded complex vector bundles with a real
structure. In this case the Hodge bigrading on the Hodge bundle
$\Lambda^1_{tot}(M,\C)$ is given by
\begin{align*}
\left(\Lambda^1_{tot}(M)\right)^{1,-1} &= S^{1,-1}(M,\C) =
\Lambda^{1,0}(M)(1),\\ 
\left(\Lambda^1_{tot}(M)\right)^{-1,1} &= S^{-1,1}(M,\C) =
\Lambda^{0,1}(M)(-1),\\ 
\left(\Lambda^1_{tot}(M)\right)^{0,0} &= \Lambda^1_{ll}(M) \oplus
\Lambda^1_{rr}(M) = \Lambda^1(M,\C).
\end{align*}
Under these identifications, the real structure on
$\Lambda^1_{tot}(M,\C)$ is minus the one induced by the usual real
structure on the complex vector bundle $\Lambda^1(M,\C)$.

\rem The Hodge bundle $S^1(M,\C)$ is canonically isomorphic to the
conormal bundle to the zero section $M \subset \tm$, which we have
described in \ref{S1}.

\punkt \label{gamma.use}
Recall now that we have defined in \ref{l.r} canonical embeddings
$\gamma_l,\gamma_r:\W_p^* \to \W_k^*$ for every $0 \leq p \leq
k$. Since $\W_0^* = \C$, for every $p,q \geq 0$ we have by
\eqref{p+q} a short exact sequence
\begin{equation}\label{cap.cup}
\begin{CD}
0 @>>> \C @>>> \W_p^* \oplus \W_q^* @>{\gamma_l \oplus \gamma_r}>>
\W_{p+q}^* @>>> 0
\end{CD}
\end{equation}
of complex vector spaces. Recall also that the embeddings
$\gamma_l$, $\gamma_r$ are compatible with the natural maps
$\can:\W_p^* \otimes \W_q^* \to \W_{p+q}^*$. Therefore the
subbundles defined by 
\begin{align*}
\Lambda^k_l(M) &= \bigoplus_{0 \leq p \leq k} \gamma_l(\W_p^*) \otimes
\Lambda^{p,k-p}(M) \subset \Lambda^k_{tot}(M) = \bigoplus_{0 \leq p \leq k} \W_k^*
\otimes \Lambda^{p,k-p}(M)\\ 
\Lambda^k_r(M) &= \bigoplus_{0 \leq p \leq k} \gamma_r(\W_p^*) \otimes
\Lambda^{k-p,p}(M) \subset \Lambda^k_{tot}(M) = \bigoplus_{0 \leq p \leq k} \W_k^*
\otimes \Lambda^{k-p,p}(M)
\end{align*}
are actually subalgebras in the total de Rham complex
$\Lambda^\cdot_{tot}(M)$. 

\punkt \label{l.r.rel} 
To describe the algebras $\Lambda^\cdot_l(M)$ and
$\Lambda^\cdot_r(M)$ explicitly, note that we obviously have
$\Lambda^\cdot_{tot}(M) = \Lambda^\cdot_l(M) +
\Lambda^\cdot_r(M)$. Moreover, in the notation of \ref{S} we have
\begin{align*}
\Lambda^1_l(M) &= S^1(M,\C) \oplus \Lambda^1_{ll}(M) \subset \Lambda^1_{tot}(M),\\
\Lambda^1_r(M) &= S^1(M,\C) \oplus \Lambda^1_{rr}(M) \subset \Lambda^1_{tot}(M).
\end{align*}
By Lemma~\ref{symm}, the algebra
$$
\Lambda^\cdot_l(M) = \left(\bigoplus_p \W_p^* \otimes
\Lambda^{p,0}(M)\right) \otimes \left( \bigoplus_q \Lambda^{0,q}(M) \right)
$$ 
is the subalgebra in the total de Rham complex $\Lambda^\cdot_{tot}(M)$
generated by $\Lambda^1_l(M)$, and the ideal of relations is
generated by the subbundle
$$
S^2(\Lambda^{1,0}(M)) \otimes \Lambda^2(\W_1^*) \subset
\Lambda^2(\Lambda^1_l(M)). 
$$
Analogously, the subalgebra $\Lambda^\cdot_r(M) \subset
\Lambda^1_{tot}(M)$ is generated by $\Lambda^1_r(M)$, and the relations
are generated by
$$
S^2(\Lambda^{0,1}(M)) \otimes \Lambda^2(\W_1^*) \subset
\Lambda^2(\Lambda^1_r(M)). 
$$

\punkt We will also need to consider the ideals in these algebras
defined by
\begin{align*}
\Lambda^k_{ll}(M) &= \bigoplus_{1 \leq p \leq k} \gamma_l(\W_p^*) \otimes
\Lambda^{p,k-p}(M) \subset \Lambda^k_l(M)\\
\Lambda^k_{rr}(M) &= \bigoplus_{1 \leq p \leq k} \gamma_r(\W_p^*) \otimes
\Lambda^{k-p,p}(M) \subset \Lambda^k_r(M)
\end{align*}
The ideal $\Lambda^\cdot_{ll}(M) \subset \Lambda^1_l(M)$ is
generated by the subbundle $\Lambda^1_{ll}(M) \subset
\Lambda^1_l(M)$, and the ideal $\Lambda^\cdot_{rr}(M) \subset
\Lambda^1_r(M)$ is generated by the subbundle $\Lambda^1_{rr}(M)
\subset \Lambda^1_r(M)$.

\punkt \label{left.right}
Denote by $\Lambda^\cdot_o(M) = \Lambda^\cdot_l(M) \cap
\Lambda^\cdot_r(M) \subset \Lambda^\cdot_{tot}(M)$ the intersection of
the subalgebras $\Lambda_l^\cdot(M)$ and $\Lambda^\cdot_r(M)$.
Unlike either of these subalgebras, the subalgebra
$\Lambda^\cdot_o(M) \subset \Lambda^\cdot_{tot}(M)$ is compatible with
the weight $0$ Hodge bundle structure on the total de Rham complex.
By \eqref{cap.cup} we have a short exact sequence
\begin{equation}\label{shrt}
\begin{CD}
0 @>>> \Lambda^\cdot(M,\C) @>>> \Lambda^\cdot_l(M) \oplus \Lambda^\cdot_r(M)
@>>> \Lambda^\cdot_{tot}(M) @>>> 0
\end{CD}
\end{equation}
of complex vector bundles on $M$. Therefore the algebra
$\Lambda^\cdot_o(M)$ is isomorphic, as a complex bundle algebra, to
the usual de Rham complex $\Lambda^\cdot(M,\C)$. As a Hodge bundle
algebra it is canonically isomorphic to the exterior algebra of the
Hodge bundle $S^1(M,\C)$ of weight $0$ on the manifold $M$.

Finally, note that the short exact sequence \eqref{shrt} induces a
direct sum decomposition
$$
\Lambda^\cdot_{tot}(M) \cong \Lambda^\cdot_{ll}(M) \oplus
\Lambda^\cdot_o(M) \oplus \Lambda^\cdot_{rr}(M). 
$$
  
\punkt \rem The total de Rham complex $\Lambda^\cdot_{tot}(M)$ is
related to Simpson's theory of Higgs bundles (see \cite{Shiggs}) in
the following way.  Recall that Simpson has proved that every
(sufficiently stable) complex bundle $\E$ on a compact complex
manifold $M$ equipped with a flat connection $\nabla$ admits a
unique Hermitian metric $h$ such that $\nabla$ and the $1$-form
$\theta = \nabla - \nabla^h \in C^\infty(M,\Lambda^1(\Eend\E))$
satisfy the so-called {\em harmonicity condition}. He also has shown
that this condition is equivalent to the vanishing of a certain
curvature-like tensor $R \in \Lambda^2(M,\Eend\E)$ which he
associated canonically to every pair $\langle \nabla, \theta
\rangle$.

Recall that flat bundles $\langle \E, \nabla \rangle$ on the
manifold $M$ are in one-to-one correspondence with free differential
graded modules $\E \otimes \Lambda^\cdot(M,\C)$ over the de Rham
complex $\Lambda^\cdot(M,\C)$.  It turns out that complex bundles
$\E$ equipped with a flat connection $\nabla$ and a $1$-form $\theta
\in C^\infty(M, \Lambda^1(\Eend\E))$ such that Simpson's tensor $R$
vanishes are in natural one-to-one correspondence with free
differential graded modules $\E \otimes \Lambda^\cdot_{tot}(M)$ over the
total de Rham complex $\Lambda^\cdot_{tot}(M)$. Moreover, a pair
$\langle\theta,\nabla\rangle$ comes from a variation of pure
$\R$-Hodge structure on $\E$ if and only if there exists a Hodge
bundle structure on $\E$ such that the product Hodge bundle
structure on the free module $\E \otimes \Lambda^\cdot_{tot}(M)$ is
compatible with the differential.

\punkt\label{gen.symm.proof}
\proof[Proof of Lemma~\ref{symm}]
For every $k \geq 0$ let $G = \Sigma_k \times \Sigma_k$ be the
product of two copies of the symmetric group $\Sigma_k$ on $k$
letters. Let $\V_k$ be the $\Q$-representation of $G_k$ induced from
the trivial representation of the diagonal subgroup $\Sigma_k
\subset G_k$. The representation $\V_k$ decomposes as
$$
\V_k = \bigoplus_V V \boxtimes V, 
$$
where the sum is over the set of irreducible representations $V$ of
$\Sigma_k$. We obviously have 
$$
\CC^k = \Hom_{G_k}\left(\V_k, A^{\otimes k} \otimes B^{\otimes
k}\right) = \bigoplus_V \Hom_{\Sigma_k}\left(V,A^{\otimes k}\right)
\otimes \Hom_{\Sigma_k}\left(V,A^{\otimes k}\right) .
$$
Let $\J^\cdot \subset \CC^\cdot$ be the ideal generated by
$\Lambda^2A \otimes \Lambda^2B \subset S^2(A \otimes B)$. It is easy
to see that
$$
\J^k = \sum_{1 \leq l \leq k-1} \bigoplus_V
\Hom_{\Sigma_k}\left(V,A^{\otimes k}\right) \otimes
\Hom_{\Sigma_k}\left(V,A^{\otimes k}\right) \subset \CC^k,
$$
where the first sum is taken over the set of $k-1$ subgroups
$\Sigma_2 \subset \Sigma_k$, the $l$-th one transposing the $l$-th
and the $l+1$-th letter, while the second sum is taken over all
irreducible constituents $V$ of the representation of $\Sigma_k$
induced from the sign representation of the corresponding $\Sigma_2
\subset \Sigma_k$. Now, there is obviously only one irreducible
representation of $\Sigma_k$ which is not encountered as an index in
this double sum, namely, the trivial one. Hence $\CC^k / \J^k = S^kA
\otimes S^kB$, which proves the lemma.  
\endproof

\subsection{The total Weil algebra}\label{t.W.sub}

\punkt \label{S.and.Lambda} 
Assume from now on that the $U(1)$-action on the complex manifold
$M$ is trivial. We now turn to studying the Weil algebra of the
manifold $M$. Let $S^1(M,\C) = S^{1,-1}(M,\C) \oplus S^{-1,1}(M,\C)$
be the weight $0$ Hodge bundle on $M$ introduced in \ref{S1}. To
simplify notation, denote
\begin{align*}
S^\cdot &= \widehat{S}^\cdot(S^1(M,\C))\\
\Lambda^\cdot &= \Lambda^\cdot(M,\C), 
\end{align*} 
where $\widehat{S}^\cdot$ is the completed symmetric power, and let 
$$
\B^\cdot = \B^\cdot(M,\C) = S^\cdot \otimes \Lambda^\cdot 
$$
be the Weil algebra of the complex manifold $M$ introduced in
\ref{Weil.defn}. Recall that the algebra $\B^\cdot$ carries a
natural Hodge bundle structure. In particular, it is equipped with a
Hodge type bigrading $\B^i = \sum_{p+q=i} \B^{p,q}$. 

\punkt \label{aug} 
We now introduce a different bigrading on the Weil algebra
$\B^\cdot$. The commutative algebra $\B^\cdot$ is freely generated
by the subbundles
$$
S^1 = S^{1,-1} \oplus S^{-1,1} \subset \B^0 \quad\text{ and }\quad
\Lambda^1 = \Lambda^{1,0} \oplus \Lambda^{0,1} \subset \B^1, 
$$
therefore to define a multiplicative bigrading on the algebra
$\B^\cdot$ it suffices to assign degrees to these generator
subbundles $S^{1,-1},S^{-1,1},\Lambda^{1,0},\Lambda^{0,1} \subset
\B^\cdot$.

\defn The {\em augmentation bigrading} on $\B^\cdot$ is the
multiplicative bigrading defined by setting 
\begin{align*}
\deg S^{1,-1} &= \deg \Lambda^{1,0} = ( 1, 0 )\\
\deg S^{-1,1} &= \deg \Lambda^{0,1} = ( 0, 1 )
\end{align*}
on generators $S^{1,-1},S^{-1,1},\Lambda^{1,0},\Lambda^{0,1} \subset
\B^{\cdot,\cdot}$. 

We will denote by $\B^{\cdot,\cdot}_{p,q}$ the component of the Weil
algebra of augmentation bidegree $(p,q)$. For any linear map
$a:\B^\cdot \to \B^\cdot$ we will denote by $a = \sum_{p,q}a_{p,q}$
its decomposition with respect to the augmentation bidegree. 

It will also be useful to consider a coarser {\em augmentation
grading} on $\B^\cdot$, defined by $\deg\B^\cdot_{p,q} = p + q$. We
will denote by $\B^\cdot_k = \bigoplus_{p+q=k}\B^\cdot_{p,q}$ the
component of $\B^\cdot$ of augmentation degree $k$. 

\punkt Note that the Hodge bidegree and the augmentation bidegree
are, in general, independent. Moreover, the complex conjugation
$\conj:\B^\cdot \to \overline{\B^\cdot}$ sends $\B^\cdot_{p,q}$ to
$\overline{\B^\cdot_{q,p}}$. Therefore the augmentation bidegree
components $\B^\cdot_{p,q} \subset \B^\cdot$ are not Hodge
subbundles. However, the coarser augmentation grading is compatible
with the Hodge structures, and the augmentation degree $k$-component
$\B^i_k \subset \B^i$ carries a natural Hodge bundle structure of
weight $i$. Moreover, the sum $\B^\cdot_{p,q} + \B^\cdot_{q,p}
\subset \B^\cdot$ is also a Hodge subbundle.

\punkt \label{total.Weil} 
We now introduce an auxiliary weight $0$ Hodge algebra bundle on
$M$, called the total Weil algebra.  Recall that we have defined in
\ref{gamma.m} a functor $\Gamma:\prehodge_{\geq 0}(M) \to
\prehodge_0(M)$ adjoint on the right to the canonical
embedding. Consider the Hodge bundle $\B_{tot}^\cdot = \Gamma(\B^\cdot)$
of weight $0$ on $M$.  By \ref{gamma.tensor} the multiplication on
$\B^\cdot$ induces an algebra structure on $\Gamma(\B^\cdot)$.

\defn The Hodge algebra bundle $\B_{tot}^\cdot$ of weight $0$ is called
{\em the total Weil algebra} of the complex manifold $M$. 

\rem For a more conceptual description of the functor $\Gamma$ and
the total Weil algebra, see Appendix.

\punkt By definition of the functor $\Gamma$ we have $\B_{tot}^k = \B^k
\otimes \W_k^* = S^\cdot \otimes \Lambda^k \otimes \W_k^* = S^\cdot
\otimes \Lambda^k_{tot}$, where $\Lambda^k_{tot} = \Lambda^k \otimes
\W_\cdot^* = \Gamma(\Lambda^k)$ is the total de Rham complex
introduced in Subsection~\ref{de.rham.sub}.  We have also introduced
in Subsection~\ref{de.rham.sub} Hodge bundle subalgebras
$\Lambda^\cdot_o, \Lambda^\cdot_l, \Lambda^\cdot_r \subset
\Lambda^\cdot_{tot}$ in the total de Rham complex $\Lambda^\cdot_{tot}$ and
ideals $\Lambda^\cdot_{ll} \subset \Lambda^\cdot_l$,
$\Lambda^\cdot_{rr} \subset \Lambda^\cdot_r$ in the algebras
$\Lambda^\cdot_l$, $\Lambda^\cdot_r$. Let
\begin{align*}
\B_o^k &= S^\cdot \otimes \Lambda_o^k \subset \B_{tot}^k\\
\B_l^k &= S^\cdot \otimes \Lambda_l^k \subset \B_{tot}^k\\
\B_r^k &= S^\cdot \otimes \Lambda_r^k \subset \B_{tot}^k
\end{align*}
be the associated subalgebras in the total Weil algebra $\B_{tot}^\cdot$
and let 
\begin{align*}
\B_{ll}^k &= S^\cdot \otimes \Lambda_{ll}^k \subset \B_l^k\\
\B_{rr}^k &= S^\cdot \otimes \Lambda_{rr}^k \subset \B_r^k
\end{align*}
be the corresponding ideals in the Hodge bundle algebras
$\B^\cdot_l$, $\B^\cdot_r$.

By \ref{left.right} we have bundle isomorphisms $\Lambda^\cdot_{tot} =
\Lambda^\cdot_l + \Lambda^\cdot_r$ and $\Lambda_o^\cdot =
\Lambda^\cdot_l \cap \Lambda^\cdot_r$, and the direct sum
decomposition $\Lambda_{tot}^\cdot \cong \Lambda_{ll}^\cdot \oplus
\Lambda_o^\cdot \oplus \Lambda_{rr}^\cdot$. Therefore we also have
\begin{align}\label{drct}
\begin{split}
\B^\cdot_{tot} &= \B_l^\cdot + \B_r^\cdot = \B^\cdot_{ll} \oplus
\B^\cdot_o \oplus \B^\cdot_{rr}\\
\B_o^\cdot &= \B_l^\cdot \bigcap \B_r^\cdot \subset \B^\cdot_{tot}
\end{split}
\end{align}
Moreover, the algebra $\Lambda^\cdot_o$ is isomorphic to the usual
de Rham complex $\Lambda^\cdot$, therefore the subalgebra
$\B_o^\cdot \subset \B^\cdot_{tot}$ is isomorphic to the usual Weil
algebra $\B^\cdot$. These isomorphisms are {\em not} weakly Hodge.

\punkt The total Weil algebra carries a canonical weight $0$ Hodge
bundle structure, and we will denote the corresponding Hodge type
grading by upper indices: $\B_{tot}^\cdot =
\oplus_p\left(\B_{tot}^\cdot\right)^{p,-p}$. The augmentation bigrading
on the Weil algebra introduced in \ref{aug} extends to a bigrading
of the total Weil algebra, which we will denote by lower indices. In
general, both these grading and the direct sum decomposition
\eqref{drct} are independent, so that, in general, for every $i \geq
0$ we have a decomposition
$$
\B_{tot}^i = \bigoplus_{n,p,q} \left(\B_{ll}^i\right)_{p,q}^{n,-n}
\oplus \left(\B_o^i\right)_{p,q}^{n,-n} \oplus
\left(\B_{rr}^i\right)_{p,q}^{n,-n}. 
$$ 
We would like to note, however, that some terms in this
decomposition vanish when $i = 0,1$. Namely, we have the following
fact.

\begin{lemma}\label{total.aug}
Let $n,k$ be arbitrary integers such that $k \geq 0$.  
\begin{enumerate}
\item If $n+k$ is odd, then $\left(\B^0_{tot}\right)_k^{n,-n} = 0$. 
\item If $n+k$ is even, then $\left(\B^1_{ll}\right)_k^{n,-n} =
\left(\B^1_{rr}\right)_k^{n,-n} = 0$, while if $n+k$ is odd, then
$\left(\B^1_o\right)_k^{n,-n} = 0$.
\end{enumerate}
\end{lemma}

\proof
\begin{enumerate}
\item The bundle $\B^0_{tot}$ by definition coincides with $\B^0$, and
it is generated by the subbundles $S^{1,-1},S^{-1,1} \subset
\B^0$. Both these subbundles have augmentation degree $1$ and Hodge
degree $\pm 1$, so that the sum $n+k$ of the Hodge degree with the
augmentation degree is even. Since both gradings are multiplicative,
for all non-zero components $\B^{n,-n}_k \subset \B^0$ the sum $n+k$
must also be even. 

\item By definition we have $\B^1_{tot} = \B^0 \otimes \Lambda^1_{tot}$. The
subbundle $\Lambda^1_{tot} \subset \B^1_{tot}$ has augmentation degree $1$,
and it decomposes 
$$
\Lambda^1_{tot} = \Lambda^1_o \oplus \Lambda^1_{ll} \oplus
\Lambda^1_{rr}.
$$
By \ref{S} we have $\Lambda^1_o \cong S^1 = S^{1,-1} \oplus
S^{-1,1}$ as Hodge bundles, so that the Hodge degrees on
$\Lambda^1_o \subset \Lambda^1_{tot}$ are odd. On the other hand, the
subbundles $\Lambda^1_{ll},\Lambda^1_{rr} \subset \Lambda^1_{tot}$ are
by \ref{S.Hodge.type} of Hodge bidegree $(0,0)$. Therefore the sum
$n+k$ of the Hodge and the augmentation degrees is even for
$\Lambda^1_o$ and odd for $\Lambda^1_{ll}$ and
$\Lambda^1_{rr}$. Together with \thetag{i} this proves the claim.
\endproof
\end{enumerate}

\subsection{Derivations of the Weil algebra}

\punkt We will now introduce certain canonical derivations of the
Weil algebra $\B^\cdot(M,\C)$ which will play an important part in
the rest of the paper. First of all, to simplify notation, for any
two linear maps $a,b$ let
$$
\{ a,b \} = a \circ b + b \circ a 
$$
be their anticommutator, and for any linear map $a:\B^\cdot \to
\B^{\cdot+i}$ let $a = \sum_{p+q=i} a^{p,q}$ be the Hodge type
decomposition. The following fact is well-known, but we have
included a proof for the sake of completeness.

\begin{lemma}
For every two odd derivations $P,Q$ of a
graded-commutative algebra $\A$, their anticommutator $\{P,Q\}$ is
an even derivation of the algebra $\A$. 
\end{lemma}

\proof Indeed, for every $a,b \in \A$ we have
\begin{align*}
\{P,Q\}(ab) &= P(Q(ab)) + Q(P(ab)) \\
& = P(Q(a)b + (-1)^{\deg a}aQ(b)) +
Q(P(a)b + (-1)^{\deg a}aP(b)) \\
& = P(Q(a))b + (-1)^{\deg Q(a)}Q(a)P(b)
+ (-1)^{\deg a}P(a)Q(b) \\
&\quad + aP(Q(b)) + Q(P(a))b + (-1)^{\deg P(a)}P(a)Q(b) \\
&\quad + (-1)^{\deg a}Q(a)P(b) + aQ(P(b)) \\
& = P(Q(a))b - (-1)^{\deg a}Q(a)P(b)
+ (-1)^{\deg a}P(a)Q(b) \\
&\quad + aP(Q(b)) + Q(P(a))b - (-1)^{\deg a}P(a)Q(b) \\
&\quad + (-1)^{\deg a}Q(a)P(b) + aQ(P(b)) \\
& = P(Q(a))b + aP(Q(b)) + Q(P(a))b + a Q(P(b)) \\
& = \{P,Q\}(a)b + a\{P,Q\}(b).
\end{align*} 
\endproof 

\punkt \label{C.and.sigma} 
Let $C:S^1 \to \Lambda^1$ be the canonical weakly Hodge map
introduced in \ref{C.Weil}. Extend $C$ to an algebra derivation
$C:\B^\cdot \to \B^{\cdot+1}$ by setting $C = 0$ on $\Lambda^1
\subset \B^1$. By \ref{C.Weil} the derivation $C$ is weakly Hodge.
The composition
$$
C \circ C = \half\{C,C\}:\B^\cdot \to \B^{\cdot+2} 
$$
is also an algebra derivation, and it obviously vanishes on
generators $S^1,\Lambda^1 \subset \B^\cdot$. Therefore $C \circ C =
0$ everywhere. 
 
Let also $\sigma:\B^{\cdot+1} \to \B^\cdot$ be the derivation
introduced in \ref{sigma}.  The derivation $\sigma$ is not weakly
Hodge; however, it is real and admits a decomposition $\sigma =
\sigma^{-1,0} + \sigma^{0,-1}$ into components of Hodge types
$(-1,0)$ and $(0,-1)$. Both these components are algebra derivations
of the Weil algebra $\B^\cdot$. We obviously have $\sigma \circ
\sigma = \sigma^{-1,0} \circ \sigma^{-1,0} = \sigma^{0,-1} \circ
\sigma^{0,-1} = 0$ on generators $S^1,\Lambda^1 \subset \B^\cdot$,
and, therefore, on the whole Weil algebra. 

\rem Up to a sign the derivations $C,\sigma$ and their Hodge
bidegree components coincide with the so-called {\em Koszul
differentials} on the Weil algebra $\B^\cdot = S^\cdot \otimes
\Lambda^\cdot$. 

\punkt \label{total.C} 
The derivation $C:\B^\cdot \to \B^{\cdot+1}$ is by definition weakly
Hodge. Applying the functor $\Gamma$ to it, we obtain a derivation
$C:\B_{tot}^\cdot \to \B_{tot}^{\cdot+1}$ of the total Weil algebra
$\B_{tot}^\cdot$ preserving the weight $0$ Hodge bundle structure on
$\B_{tot}^\cdot$. The canonical identification $\B^\cdot \cong
\B^\cdot_o \subset \B_{tot}^\cdot$ is compatible with the derivation
$C:\B^\cdot_{tot} \to \B^{\cdot+1}_{tot}$. Moreover, by \ref{C.and.sigma}
this derivation satisfies $C \circ C = 0:\B^\cdot_{tot} \to
\B^{\cdot+2}_{tot}$. Therefore the total Weil algebra $\B^\cdot_{tot}$
equipped with the derivation $C$ is a complex of Hodge bundles of
weight $0$. 

The crucial linear algebraic property of the total Weil algebra
$\B^\cdot_{tot}$ of the manifold $M$ which will allow us to classify
flat extended connections on $M$ is the following.

\begin{prop}\label{ac}
Consider the subbundle  
\begin{equation}\label{sbcmp}
\bigoplus_{p,q \geq 1}\left(\B^\cdot_{tot}\right)_{p,q} \subset
\B^\cdot_{tot}
\end{equation}
of the total Weil algebra $\B^\cdot_{tot}$ consisting of the components
of augmentation bidegrees $(p,q)$ with $p,q \geq 1$. This subbundle
equipped with the differential
$C:\left(\B^\cdot_{tot}\right)_{\cdot,\cdot} \to
\left(\B^{\cdot+1}_{tot}\right)_{\cdot,\cdot}$ is an acyclic complex of
Hodge bundles of weight $0$ on $M$.
\end{prop}

\punkt We sketch a more or less simple and conceptual proof of
Proposition~\ref{ac} in the Appendix. However, in order to be able
to study in Section~\ref{convergence} the analytic properties of our
formal constructions, we will need an explicit contracting homotopy
for the complex \eqref{sbcmp}, which we now introduce.

The restriction of the derivation $C:\B^0_{tot} \to \B^1_{tot}$ to the
subbundle $S^1 \subset \B_0 \cong \B^0_{tot}$ induces a Hodge bundle
isomorphism
$$
C:S^1 \to \Lambda^1_o \subset \Lambda^1_{tot} \subset \B^1_{tot}.
$$
Define a map $\sigma_{tot}:\Lambda^1_{tot} \to S^1$ by 
\begin{equation}\label{sigma.c.eq}
\sigma_{tot} = 
\begin{cases} 
0 \quad &\text{ on }\quad\Lambda^1_{ll},\Lambda^1_{rr} \subset \Lambda^1_{tot},\\
C^{-1} \quad &\text{ on }\quad\Lambda^1_o \subset \Lambda^1_{tot}.
\end{cases}
\end{equation}
The map $\sigma_{tot}:\Lambda^1_{tot} \to S^1$ preserves the Hodge bundle
structures of weight $0$ on both sides. Moreover, its restriction to
the subbundle $\Lambda^1 \cong \Lambda^1_o \subset \Lambda^1_{tot}$
coincides with the canonical map $\sigma:\Lambda^1 \to S^1$
introduced in \ref{C.and.sigma}. 

\punkt \label{sigma.l} 
Unfortunately, unlike $\sigma:\Lambda^1 \to S^1$, the map
$\sigma_{tot}:\Lambda^1_{tot} \to S^1$ does {\em not} admit an extension to
a derivation $\B^{\cdot+1}_{tot} \to \B_{tot}^\cdot$ of the total Weil
algebra $\B^\cdot_{tot}$. We will extend it to a {\em bundle map}
$\sigma_{tot}:\B^{\cdot+1}_{tot} \to \B^\cdot_{tot}$ in a somewhat roundabout
way. To do this, define a map $\sigma_l:\Lambda^1_l \to S^0 \subset
\B^0$ by
$$
\sigma_l =  
\begin{cases}
0 \quad &\text{ on } \quad\Lambda_{ll}^1 \subset \Lambda_l^1 \quad
\text{ and on }\quad \left(\Lambda^1_o\right)^{-1,1} \subset
\Lambda^1_{tot}, \\ 
C^{-1} \quad &\text{ on } \quad\left(\Lambda^1_o\right)^{1,-1} \subset
\Lambda^1_{tot}.
\end{cases}
$$
and set $\sigma_l = 0$ on $S^1$. By \ref{gamma.use} we have
$$
\Lambda^1_l = \Lambda^{1,0} \oplus \left( \Lambda^{0,1} \otimes
\W_1^*\right). 
$$
The map $\sigma_l:\Lambda^1_l \to S^1$ vanishes on the second
summand in this direct sum, and it equals $C^{-1}:\Lambda^{1,0} \to
S^{1,-1} \subset S^1$ on the first summand. The restriction of the map
$\sigma_l$ to the subbundle 
$$
\Lambda^{1,0} \oplus \Lambda^{0,1} = \Lambda^1 \cong \Lambda^1_o \subset
\Lambda^1_l
$$ 
vanishes on $\Lambda^{0,1}$ and equals $C^{-1}$ on
$\Lambda^{1,0}$. Thus it is equal to the Hodge type-$(0,-1)$
component $\sigma^{0,-1}:\Lambda^1 \to S^1$ of the canonical map
$\sigma:\Lambda^1 \to S^1$.

\punkt By \ref{l.r.rel} the algebra $\B^\cdot_l$ is generated by the
bundles $S^1$ and $\Lambda^1_l$, and the ideal of relations is
generated by the subbundle
\begin{equation}\label{rell}
S^2\left(\Lambda^{0,1}\right) \otimes \Lambda^2\left(\W_1^*\right) 
\subset \Lambda^2\left(\Lambda^1_l\right).
\end{equation}
Since the map $\sigma_l:\Lambda^1_l \to S^1$ vanishes on
$\Lambda^{0,1} \otimes \W_1^* \subset \Lambda^1_l$, it extends to an
algebra derivation $\sigma_l:\B_l^{\cdot+1} \to \B^\cdot_l$. The
restriction of the derivation $\sigma_l$ to the subalgebra $\B^\cdot
\cong \B^\cdot_0 \subset \B_{tot}^\cdot$ coincides with the
$(0,-1)$-component $\sigma^{0,-1}$ of the derivation
$\sigma:\B^{\cdot+1} \to \B^\cdot$.

Analogously, the $(-1,0)$-component $\sigma^{-1,0}$ of the
derivation $\sigma:\B^{\cdot+1} \to \B^\cdot$ extends to an algebra
derivation $\sigma_r:\B_r^{\cdot+1} \to \B_r^\cdot$ of the
subalgebra $\B^\cdot_r \subset \B_{tot}^\cdot$. By definition, the
derivation $\sigma_l$ preserves the decomposition $\B_l^\cdot =
\B_{ll} \oplus \B_o^\cdot$, while the derivation $\sigma_r$
preserves the decomposition $\B_r^\cdot = \B_{rr}^\cdot \oplus
\B_o^\cdot$. Both these derivations vanish on $\B_{tot}^0$, therefore
both are maps of $\B_{tot}^0$-modules. In addition, the compositions
$\sigma_l \circ \sigma_l$ and $\sigma_r \circ \sigma_r$ vanish on
generator and, therefore, vanish identically.

\punkt \label{sigma.c}
Extend both $\sigma_l$ and $\sigma_r$ to the whole $\B_{tot}^\cdot$ by
setting
\begin{equation}\label{nol}
\sigma_l = 0 \text{ on }\B^\cdot_r \qquad \sigma_r = 0 \text{ on
}\B^\cdot_l, 
\end{equation}
and let
$$
\sigma_{tot} = \sigma_l + \sigma_r:\B^{\cdot+1}_{tot} \to \B_{tot}^\cdot.
$$
On $\Lambda^1_{tot} \subset \B^1_{tot}$ this is the same map as in
\eqref{sigma.c.eq}. The bundle map $\sigma_{tot}:\B_{tot}^{\cdot+1} \to
\B_{tot}^\cdot$ preserves the direct sum decomposition \eqref{drct}, and
its restriction to $\B_o^\cdot \subset \B^\cdot_{tot}$ coincides with
the derivation $\sigma$.  Note that neither of the maps $\sigma_l$,
$\sigma_r$, $\sigma_{tot}$ is a derivation of the total Weil algebra
$\B_{tot}^\cdot$. However, all these maps are linear with respect to the
$\B_{tot}^0$-module structure on $\B_{tot}^\cdot$ and preserve the
decomposition \eqref{drct}. The map $\sigma_{tot}:\B^{\cdot+1}_{tot} \to
\B^\cdot_{tot}$ is equal to $\sigma_l$ on $\B^\cdot_{ll} \subset
\B^\cdot_{tot}$, to $\sigma_r$ on $\B^\cdot_{rr}$ and to
$\sigma:\B^{\cdot+1} \to \B^\cdot$ on $\B^\cdot \cong \B^\cdot_o
\subset \B^\cdot_{tot}$. Since $\sigma_l \circ \sigma_l = \sigma_r \circ
\sigma_r = \sigma \circ \sigma = 0$, we have $\sigma_{tot} \circ
\sigma_{tot} = 0$. 

\punkt The commutator 
$$
h = \{C, \sigma_{tot}\}:\B_{tot}^\cdot \to \B_{tot}^\cdot
$$
of the maps $C$ and $\sigma_{tot}$ also preserves the decomposition
\eqref{drct}, and we have the following. 

\begin{lemma}\label{h.acts}
The map $h$ acts as multiplication by $p$ on
$\left(\B_{ll}^\cdot\right)_{p,q}$, as multiplication by $q$ on
$\left(\B_{rr}^\cdot\right)_{p,q}$ and as multiplication by $(p+q)$
on $\left(\B_o^\cdot\right)_{p,q}$.
\end{lemma}

\proof It suffices to prove the claim separately on each term in the
decomposition \eqref{drct}. By definition $\sigma_{tot} = \sigma_l +
\sigma_r$, and $h = h_l + h_r$, where $h_l = \{\sigma_l,C\}$ and
$h_r = \{\sigma_r,C\}$. Moreover, $h_l$ vanishes on $\B_{rr}^\cdot$
and $h_r$ vanishes on $\B_{ll}^\cdot$. Therefore it suffices to
prove that $h_l = p\id$ on $\left(\B^\cdot_l\right)_{p,q}$ and that
$h_r = q\id$ on $\left(\B^\cdot_r\right)_{p,q}$. The proofs of these
two identities are completely symmetrical, and we will only give a
proof for $h_l$.

The algebra $\B^\cdot_l$ is generated by the subbundles $S^1 \subset
\B^0_l$ and $\Lambda^1_l \subset \B^1_l$. The augmentation bidegree
decomposition of $S^1$ is by definition given by
$$
S^1_{1,0} = S^{1,-1} \qquad\qquad S^1_{0,1} = S^{-1,1},
$$
while the augmentation bidegree decomposition of $\Lambda^1_l$ is
given by
$$
\left(\Lambda^1_l\right)_{1,0} = \Lambda^{1,0} \qquad
\left(\Lambda^1_l\right)_{0,1} = \Lambda^{0,1} \otimes \W_1^*. 
$$
By the definition of the map $\sigma_l:\Lambda^1_l \to S^1$ (see
\ref{sigma.l}) we have $h_l = \{C,\sigma_l\} = \id$ on
$\Lambda^{1,0}$ and $S^{1,-1}$, and $h_l = 0$ on $\Lambda^{0,1}
\otimes \W_1^*$ and on $S^{-1,1}$. Therefore for every $p,q \geq 0$
we have $h_l=p\id$ on the generator subbundles $S^1_{p,q}$ and on
$\left(\Lambda^1_l\right)_{p,q}$. Since the map $h_l$ is a
derivation and the augmentation bidegree is multiplicative, the same
holds on the whole algebra $\B^\cdot_l =
\oplus\left(\B^\cdot_l\right)_{p,q}$. 
\endproof

Lemma~\ref{h.acts} shows that the map $\sigma_{tot}$ is a homotopy,
contracting the subcomplex \eqref{sbcmp} in the total Weil algebra
$\B^\cdot_{tot}$, which immediately implies Proposition~\ref{ac}.

\rem In fact, in our classification of flat extended connections
given in Section~\ref{main.section} it will be more convenient for
us to use Lemma~\ref{h.acts} directly rather than refer to
Proposition~\ref{ac}.

\punkt We finish this section with the following corollary of
Lemma~\ref{total.aug} and Lemma~\ref{h.acts}, which we will need in
Section~\ref{convergence}.

\begin{lemma}\label{h.on.b1}
Let $n=\pm 1$. If the integer $k \geq 1$ is odd, then the map
$h:\B^\cdot_{tot} \to \B^\cdot_{tot}$ acts on $\left(\B_{tot}^1\right)^{n,-n}_k$
by multiplication by $k$. If $k = 2m \geq 1$ is even, then the
endomorphism $h:\left(\B_{tot}^1\right)^{n,-n}_k \to
\left(\B_{tot}^1\right)^{n,-n}_k$ is diagonalizable, and its only
eigenvalues are $m$ and $m-1$.
\end{lemma}

\proof If $k$ is odd, then $\left(\B_{tot}^1\right)^{n,-n}_k =
\left(\B_o^1\right)^{n,-n}_k$ by Lemma~\ref{total.aug}, and
Lemma~\ref{h.acts} immediately implies the claim.  Assume that the
integer $k = 2m$ is even. By Lemma~\ref{total.aug} we have
$$
\left(\B_{tot}^\cdot\right)^{n,-n}_k =
\left(\B_{ll}^\cdot\right)^{n,-n}_k \oplus
\left(\B_{rr}^\cdot\right)^{n,-n}_k = \B^{n,-n}_{k-1} \otimes
\left(\Lambda^1_{ll} \oplus \Lambda^1_{rr}\right).
$$
The bundle $\B^0$ is generated by subbundles $S^{1,-1}$ and
$S^{-1,1}$. The first of these subbundles has augmentation bidegree
$(1,0)$, while the second one has augmentation bidegree
$(0,1)$. Therefore for every augmentation bidegree component
$\B^{n,-n}_{p,q} \subset \B^{n,-n}_{k-1}$ we have $p - q = n$ and
$p+q=k-1$. This implies that $\B^{n,-n}_{k-1} = \B^{n,-n}_{p,q}$
with $p=m-(1-n)/2$ and $q=m-(1+n)/2$.
  
By definition the augmentation bidegrees of the bundles
$\Lambda^1_{ll}$ and $\Lambda^1_{rr}$ are, respectively, $(0,1)$ and
$(1,0)$. Lemma~\ref{h.acts} shows that the only eigenvalue of the
map $h$ on $\left(\B^1_{ll}\right)_{k+1}$ is $p=(m-(1-n)/2))$, while
its only eigenvalue on $\left(\B^1_{rr}\right)_{k+1}$ is
$q=(m-(1+n)/2)$. Since $n = \pm 1$, one of these numbers equals $m$
and the other one equals $m-1$.  
\endproof

\section{Classification of flat extended connections}\label{main.section} 

\subsection{K\"ah\-le\-ri\-an connections}

\punkt Let $M$ be a complex manifold. In
Section~\ref{formal.section} we have shown that formal Hodge
manifold structures on the tangent bundle $\tm$ are in one-to-one
correspondence with linear flat extended connections on the manifold
$M$ (see \ref{ext.con}--\ref{flat.ext.con} for the definitions). It
turns out that flat linear extended connections on $M$ are, in turn,
in natural one-to-one correspondence with differential operators of
a much simpler type, namely, connections on the cotangent bundle
$\Lambda^{1,0}(M)$ satisfying certain vanishing conditions
(Theorem~\ref{kal=ext}). We call such connections {\em
K\"ah\-le\-ri\-an}. In this section we use the results of
Section~\ref{Weil.section} establish the correspondence between
extended connections on $M$ and K\"ah\-le\-ri\-an connections on
$\Lambda^{1,0}(M)$.

\punkt We first give the definition of K\"ah\-le\-ri\-an
connections. Assume that the manifold $M$ is equipped with a
connection
$$
\nabla:\Lambda^1(M) \to \Lambda^1(M) \otimes \Lambda^1(M) 
$$
on its cotangent bundle $\Lambda^1(M)$. Let 
\begin{align*}
T &= \Alt \circ \nabla - d_M:\Lambda^1(M) \to \Lambda^2(M) \\
R &= \Alt \nabla \circ \nabla:\Lambda^1(M) \to \Lambda^1(M) \otimes
\Lambda^2(M) 
\end{align*}
be its torsion and curvature, and let $R = R^{2,0} + R^{1,1} +
R^{0,2}$ be the decomposition of the curvature according to the
Hodge type.  

\defn The connection $\nabla$ is called {\em K\"ah\-le\-ri\-an} if 
\begin{align*}
T &= 0 \tag{i}\\
R^{2,0} &= R^{0,2} = 0 \tag{ii}
\end{align*}

\ex The Levi-Civita connection on a K\"ahler manifold is K\"ah\-le\-ri\-an. 

\rem The condition $T=0$ implies, in particular, that the component 
$$
\nabla^{0,1}:\Lambda^{1,0}(M) \to \Lambda^{1,1}(M)
$$
of the connection $\nabla$ coincides with the Dolbeault
differential. Therefore a K\"ah\-le\-ri\-an connection is always
holomorphic.

\punkt Recall that in \ref{red} we have associated to any extended
connection $D$ on $M$ a connection $\nabla$ on the cotangent bundle
$\Lambda^1(M,\C)$ called the reduction of $D$. We can now formulate
the main result of this section.

\begin{theorem}\label{kal=ext}
\begin{enumerate}
\item If an extended connection $D$ on $M$ is flat and linear, then
its reduction $\nabla$ is K\"ah\-le\-ri\-an. 
\item Every K\"ah\-le\-ri\-an connection $\nabla$ on $\Lambda^1(M,\C)$ is
the reduction of a unique linear flat extended connection $D$ on
$M$. 
\end{enumerate}
\end{theorem} 

The rest of this section is taken up with the proof of
Theorem~\ref{kal=ext}. To make it more accessible, we first give an
informal outline. The actual proof starts with
Subsection~\ref{pf.first.sub}, and it is independent from the rest
of this subsection.

\punkt Assume given a K\"ah\-le\-ri\-an connection $\nabla$ on the
manifold $M$. To prove Theorem~\ref{kal=ext}, we have to construct a
flat linear extended connection $D$ on $M$ with reduction $\nabla$.
Every extended connection decomposes into a series $D = \sum_{k \geq
0}D_k$ as in \eqref{aug.con}, and, since $\nabla$ is the reduction
of $D$, we must have $D_1 = \nabla$. We begin by checking in
Lemma~\ref{lin.aug} that if $D$ is linear, then $D_0 = C$, where
$C:S^1(M,\C) \to \Lambda^1(M,\C)$ is as in \ref{C.Weil}. The sum $C
+ \nabla$ is already a linear extended connection on $M$. By
\ref{deriv} it extends to a derivation $D_{\leq 1}$ of the Weil
algebra $\B^\cdot(M,\C)$ of the manifold $M$, but this derivation
does not necessarily satisfy $D_{\leq 1} \circ D_{\leq 1} = 0$, thus
the extended connection $D_{\leq 1}$ is not necessarily flat.

We have to show that one can add the ``correction terms'' $D_k, k
\geq 2$ to $D_{\leq 1}$ so that $D = \sum_k D_k$ satisfies all the
conditions of Theorem~\ref{kal=ext}. To do this, we introduce in
\ref{red.Weil} a certain quotient $\wB^\cdot(M,\C)$ of the Weil
algebra $\B^\cdot(M,\C)$, called the reduced Weil algebra. The
reduced Weil algebra ia defined in such a way that for every
extended connection $D$ the associated derivation $D:\B^\cdot(M,\C)
\to \B^{\cdot+1}(M,\C)$ preserves the kernel of the surjection
$\B^\cdot \to \wB^\cdot$, thus inducing a derivation
$\wD:\wB^\cdot(M,\C) \to \wB^{\cdot+1}(M,\C)$. Moreover, the algebra
$\wB^\cdot(M,\C)$ has the following two properties:
\begin{enumerate}
\item \label{first.prop} The derivation $\wD:\wB^\cdot(M,\C) \to
\wB^{\cdot+1}(M,\C)$ satisfies $\wD \circ \wD = 0$ if and only if
the connection $D_1$ is K\"ah\-le\-ri\-an.
\item \label{second.prop} Let $\wD$ be the weakly Hodge derivation
of the quotient algebra $\wB^\cdot(M,\C)$ induced by an arbitrary
linear extended connection $D_{\leq 1}$ and such that $\wD \circ \wD
= 0$. Then the derivation $\wD$ lifts uniquely to a weakly Hodge
derivation $D$ of the Weil algebra $\B^\cdot(M,\C)$ such that $D
\circ D = 0$, and the derivation also $D$ comes from a linear
extended connection on $M$ (see Proposition~\ref{main} for a precise
formulation of this statement).
\end{enumerate}

\punkt The property \ref{first.prop} is relatively easy to check,
and we do it in the end of the proof, in
Subsection~\ref{pf.last.sub}. The rest is taken up with establishing
the property \ref{second.prop}. The actual proof of this statement
is contained in Proposition~\ref{main}, and
Subsection~\ref{pf.first.sub} contains the necessary preliminaries.

Recall that we have introduced in \ref{aug} a new grading on the
Weil algebra $\B^\cdot(M,\C)$, called the augmentation grading, so
that the component $D_k$ in the decomposition $D = \sum_kD_k$ is of
augmentation degree $k$. In order to lift $\wD$ to a derivation $D$
so that $D \circ D = 0$, we begin with the given lifting $D_{\leq
1}$ and then add components $D_k, k \geq 2$, one by one, so that on
each step for $D_{\leq k} = D_{\leq 1} + \sum_{2 \leq p \leq k}$ the
composition $D_{\leq k} \circ D_{\leq k}$ is zero in augmentation
degrees from $0$ to $k$. In order to do it, we must find for each
$k$ a solution to the equation
\begin{equation}\label{tslv}
D_0 \circ D_k  = - R_k,
\end{equation}
where $R_k$ is the component of augmentation degree $k$ in the
composition $D_{\leq k-1} \circ D_{\leq k-1}$. This solution must be
weakly Hodge, and the extended connection $D_{\leq k} = D_{\leq k-1}
+ D_k$ must be linear.

We prove in Lemma~\ref{lin.aug} that since $D_{\leq 0}$ is linear,
we may assume that $D_0 = C$. In addition, since $\wD \circ \wD =
0$, we may assume by induction that the image of $R_k$ lies in the
kernel $\I^\cdot$ of the quotient map $\B^\cdot(M,\C) \to
\wB^\cdot(M,\C)$.

\punkt In order to analyze weakly Hodge maps from $S^1(M,\C)$ to the
Weil algebra $\B^\cdot(M,\C)$, we apply the functor
$\Gamma:\prehodge_{\geq 0}(M) \to \prehodge_0(M)$ constructed in
\ref{gamma.m} to the bundle $\B^\cdot(M,\C)$ to obtain the total
Weil algebra $\B_{tot}^\cdot(M,\C) = \Gamma(\B^\cdot(M,\C))$ of weight
$0$, which we studied in Subsection~\ref{t.W.sub}. The Hodge bundle
$S^1(M,\C)$ on the manifold $M$ is of weight $0$, and, by the
universal property of the functor $\Gamma$, weakly Hodge maps from
$S^1(M,\C)$ to $\B^\cdot(M,\C)$ are in one-to-one correspondence
with Hodge bundle maps from $S^1(M,\C)$ to the total Weil algebra
$\B^\cdot_{tot}(M,\C)$. The canonical map $C:S^1(,\C) \to \B^1(M,\C)$
extends to a derivation $C:\B^\cdot_{tot}(M,\C) \to
\B^{\cdot+1}_{tot}(M,\C)$. Moreover, the weakly Hodge map $R_k:S^1(M,\C)
\to \B^2(M,\C)$ defines a Hodge bundle map $R_k^{tot}:S^1(M,\C) \to
\B^2_{tot}(M,\C)$, and solving \eqref{tslv} is equivalent to finding a
Hodge bundle map $D_k:S^1(M,\C) \to \B^1_{tot}(M,\C)$ such that
\begin{equation}\label{tslvc}
C \circ D_k = - R_k.
\end{equation}

\punkt Recall that by \ref{C.and.sigma} the derivation
$$
C:\B^\cdot_{tot}(M,\C) \to \B^{\cdot+1}_{tot}(M,\C)
$$ 
satisfies $C \circ C = 0$, so that the total Weil algebra
$\B^\cdot_{tot}(M,\C)$ becomes a complex with differential $C$. The
crucial part of the proof of Theorem~\ref{kal=ext} consists in
noticing that the subcomplex $\I_{tot}^\cdot(M,\C) =
\Gamma(\I^\cdot(M,\C)) \subset \B^\cdot_{tot}(M,\C)$ of the total Weil
algebra $\B^\cdot_{tot}(M,\C)$ corresponding to the kernel
$\I^\cdot(M,\C) \subset \B^\cdot(M,\C)$ of the quotient map
$\B^\cdot(M,\C) \to \wB^\cdot(M,\C)$ is canonically contractible.
This statement is analogous to Proposition~\ref{ac}, and we prove it
in the same way. Namely, we check that the subcomplex $\I^\cdot_{tot}
\subset \B^\cdot_{tot}$ is preserved by the bundle map
$\sigma_{tot}:\B^{\cdot+1}_{tot} \to \B^\cdot_{tot}$ constructed in
\ref{sigma.c}, and that the anticommutator $h =
\{\sigma_{tot},C\}:\B^\cdot_{tot}(M,\C) \to \B^\cdot_{tot}(M,\C)$ is invertible
on the subcomplex $\I_{tot}^\cdot(M,\C) \subset \B^\cdot_{tot}(M,\C)$
(Corollary~\ref{h.inv} of Lemma~\ref{h.acts}). We also check that $C
\circ R_k^{tot} = 0$, which implies that the Hodge bundle map
\begin{equation}\label{sltn}
D_k = -h^{-1} \circ \sigma_{tot} \circ R_k^{tot}:S^1(M,\C) \to \B^1_{tot}(M,\C)
\end{equation}
provides a solution to the equation \eqref{tslvc}. 

\punkt To establish the property \ref{second.prop}, we have to
insure additionally that the extended connection $D = D_{\leq k}$ is
linear, and and we have to show that the solution $D_k$ of
\eqref{tslvc} with this property is unique. This turns out to be
pretty straightforward. We show in Lemma~\ref{lin.aug} that $D_{\leq
k}$ is linear if and only if
\begin{equation}\label{lnr}
\sigma_{tot} \circ D_k = 0.
\end{equation}
Moreover, we show that the homotopy $\sigma_{tot}:\B^{\cdot+1}_{tot}(M,\C)
\to \B^\cdot_{tot}(M,\C)$ satisfies $\sigma_{tot} \circ \sigma_{tot} =
0$. Therefore the solution $D_k$ to \eqref{tslvc} given by
\eqref{sltn} satisfies \eqref{lnr} automatically.

The uniqueness of such a solution $D_k$ follows from the
invertibility of $h = C \circ \sigma_{tot} + \sigma_{tot} \circ C$. Indeed,
for every two solutions $D_k,D_k'$ to \eqref{tslv}, both satisfying
\eqref{lnr}, their difference $P = D_k - D'_k$ satisfies $C \circ P
= 0$. If, in addition, both $D_k$ and $D'_k$ were to satisfy
\eqref{lnr}, we would have had $\sigma_{tot} \circ P = 0$. Therefore $h
\circ P = 0$, and $P$ has to vanish.

\punkt These are the main ideas of the proof of
Theorem~\ref{kal=ext}. The proof itself begins in the next
subsection, and it is organized as follows. In
Subsection~\ref{pf.first.sub} we express the linearity condition on
an extended connection $D$ in terms of the associated derivation
$D_{tot}:\B^\cdot_C \to \B^{\cdot+1}_{tot}$ of the total Weil algebra
$\B^\cdot_{tot}$ of the manifold $M$. After that, we introduce in
Subsection~\ref{pf.third.sub} the reduced Weil algebra
$\wB^\cdot(M,\C)$ and prove Proposition~\ref{main}, thus reducing
Theorem~\ref{kal=ext} to a statement about derivations of the
reduced Weil algebra. Finally, in Subsection~\ref{pf.last.sub} we
prove this statement.

\rem In the Appendix we give, following Deligne and Simpson, a
more geomteric description of the functor $\Gamma:\prehodge_{\geq 0}
\to \prehodge$ and of the total Weil algebra $\B^\cdot_{tot}(M,\C)$,
which allows to give a simpler and more conceptual proof for the key
parts of Theorem~\ref{kal=ext}.

\subsection{Linearity and the total Weil algebra}\label{pf.first.sub}

\punkt Assume given an extended connection $D:S^1 \to \B^1$ on the
manifold $M$, and extend it to a derivation $D:\B^\cdot \to
\B^{\cdot+1}$ of the Weil algebra as in \ref{deriv}.  Let $D =
\sum_{k\geq 0} D_k$ be the augmentation degree decomposition.  The
derivation $D$ is weakly Hodge and defines therefore a derivation $D
= \sum_{k\geq 0} D_k:\B_{tot}^\cdot \to \B_{tot}^{\cdot+1}$ of the total
Weil algebra $\B_{tot}^\cdot$.

Before we begin the proof of Theorem~\ref{kal=ext}, we give the
following rewriting of the linearity condition \ref{lin.ext.con} on
the extended connection $D$ in terms of the total Weil algebra.

\begin{lemma}\label{lin.aug}
The extended connection $D$ is linear if and only if $D_0=C$ and
$\sigma_{tot} \circ D_k = 0$ on $S^1 \subset \B_{tot}^0$ for every $k \geq
0$. 
\end{lemma}

\proof Indeed, by Lemma~\ref{total.aug} for odd integers $k$ and $n
= \pm 1$ the subbundle $\left(\B_o^1\right)^{n,-n}_{k+1} \subset
\B_{tot}^1$ vanishes. Therefore the map $D_k:S^1 \to \B_{tot}^1$ factors
through $\B_{ll}^1 \oplus \B_{rr}^1$.  Since by definition
(\ref{sigma.c}) we have $\sigma_{tot}=0$ on both $\B_{ll}^1$ and
$\B_{rr}^1$, for odd $k$ we have $\sigma_{tot} \circ D_k = 0$ on $S^1$
regardless of the extended connection $D$. On the other hand, for
even $k$ we have $\left(\B_{tot}^1\right)^{n,-n}_{k+1} =
\left(\B_o^1\right)^{n,-n}_{k+1}$. Therefore on $S^1$ we have
$\sigma_{tot} \circ D_k = \sigma \circ D_k$ (where $\sigma:\B^{\cdot+1}
\to \B^\cdot$ is as in \ref{C.and.sigma}). Moreover, since
$\sigma:\Lambda^1 \to S^1$ is an isomorphism, $D_0 = C$ is
equivalent to $\sigma \circ D_0 = \id:S^1 \to \Lambda^1 \to
S^1$. Therefore the condition of the lemma is equivalent to the
following
\begin{equation}\label{tprv}
\sigma \circ D_k = 
\begin{cases} 
\id, \quad &\text{ for } k=0\\
0,   \quad &\text{ for even integers } k > 0.
\end{cases}
\end{equation}
Let now $\iota^*:\B^\cdot \to \B^\cdot$ be the operator given by the
action of the canonical involution $\iota:\tm \to \tm$, as in
\ref{iota.Weil}, and let $D^\iota = \sum_{k \geq 0} D_k^\iota =
\iota^* \circ D \circ(\iota^*)^{-1}$ be the operator
$\iota^*$-conjugate to the derivation $D$. The operator $\iota^*$
acts as $-\id$ on $S^1 \subset \B^0$ and as $\id$ on $\Lambda^1
\subset \B^1$. Since it is an algebra automorphism, it acts as
$(-1)^{i+k}$ on $\B_k^i \subset \B^\cdot$. Therefore $D_k^\iota =
(-1)^{k+1}D_k$, and \eqref{tprv} is equivalent to
$$
\sigma \circ \half(D - D^\iota) = \id:S^1 \to \B^1 \to S^1, 
$$
which is precisely the definition of a linear extended connection. 
\endproof 

\subsection{The reduced Weil algebra}\label{pf.third.sub}

\punkt We now begin the proof of Theorem~\ref{kal=ext}. Our first
step is to reduce the classification of linear flat extended
connections $D:S^1 \to \B^1$ on the manifold $M$ to the study of
derivations of a certain quotient $\wB^\cdot$ of the Weil algebra
$\B^\cdot$. We introduce this quotient in this subsection under the
name of reduced Weil algebra. We then show that every extended
connection $D$ on $M$ induces a derivation $\wD:\wB^\cdot \to
\wB^{\cdot+1}$ of the reduced Weil algebra, and that a linear flat
extended connection $D$ on $M$ is completely defined by the
derivation $\wD$.

\punkt By Lemma~\ref{h.acts} the anticommutator $h =
\{C,\sigma_{tot}\}:\B^\cdot_{tot} \to \B^\cdot_{tot}$ of the canonical bundle
endomorphisms $C,\sigma_{tot}$ of the total Weil algebra $\B_{tot}^\cdot$ is
invertible on every component $\left(\B^\cdot_{tot}\right)_{p,q}$ of
augmentation bidegree $(p,q)$ with $p,q \geq 1$.

The direct sum $\oplus_{p,q \geq 1}\left(\B_{tot}^\cdot\right)_{p,q}$ is
an ideal in the total Weil algebra $\B^\cdot_{tot}$, and it is obtained
by applying the functor $\Gamma$ to the ideal $\oplus_{p,q \geq
1}\B^\cdot_{p,q}$ in the Weil algebra $\B^\cdot$. For technical
reasons, it will be more convenient for us to consider the smaller
subbundle
$$
\I^\cdot = \bigoplus_{p \geq 2,q \geq 1}\B^\cdot_{p,q} +
\bigoplus_{p \geq 1, q \geq 2} \B^\cdot_{p,q} \subset \B^\cdot 
$$
of the Weil algebra $\B^\cdot$. The subbundle $\I^\cdot$ is a Hodge
subbundle in $\B^\cdot$, and it is an ideal with respect to the
multiplication in $\B^\cdot$. 

\punkt \defn \label{red.Weil}
The reduced Weil algebra $\wB^\cdot = \wB^\cdot(M,\C)$ of the
manifold $M$ is the quotient
$$
\wB^\cdot = \B^\cdot / \I^\cdot
$$
of the full Weil algebra $\B^\cdot$ by the ideal $\I^\cdot$. 

The reduced Weil algebra decomposes as 
$$
\wB^\cdot = \B^\cdot_{1,1} \oplus \bigoplus_{p \geq 0}\B^\cdot_{p,0}
\oplus \bigoplus_{q \geq 0}\B^\cdot_{0,q} 
$$
with respect to the augmentation bigrading on the Weil algebra
$\B^\cdot$. The two summands on the right are equal to  
\begin{align*}
\bigoplus_{p \geq 0}\B_{p,0}^\cdot &= \bigoplus_{p geq
0}S^p\left(S^{1,-1}\right) \otimes \Lambda^{\cdot,0}, \\ 
\bigoplus_{q \geq 0}\B_{0,q}^\cdot &= \bigoplus_{q geq
0}S^q\left(S^{-1,1}\right) \otimes \Lambda^{0,\cdot},
\end{align*} 

\punkt Since $\I^\cdot$ is a Hodge subbundle, the reduced Weil
algebra carries a canonical Hodge bundle structure compatible with
the multiplication. It also obviously inherits the augmentation
bigrading, and defines an ideal $\I_{tot}^\cdot = \Gamma(\I^\cdot)
\subset \B_{tot}^\cdot$ in the total Weil algebra
$\B_{tot}^\cdot$. Lemma~\ref{h.acts} immediately implies the following
fact.

\begin{corr}\label{h.inv}
The map $h=\{C,\sigma_{tot}\}:\B_{tot}^\cdot \to \B_{tot}^\cdot$ is invertible
on $\I_{tot}^\cdot \subset \B_{tot}^\cdot$. 
\end{corr} 

\punkt \label{wD} 
Let now $D:\B^\cdot \to \B^\cdot$ be the derivation associated to
the extended connection $D$ as in \ref{deriv}. The derivation $D$
does not increase the augmentation bidegree, it preserves the ideal
$\I^\cdot \subset \B^\cdot$ and defines therefore a weakly Hodge
derivation of the reduced Weil algebra $\wB^\cdot$, which we denote
by $\wD$. If the extended connection $D$ is flat, then the
derivation $\wD$ satisfies $\wD \circ \wD = 0$.

We now prove that every derivation $\wD:\wB^\cdot \to \wB^{\cdot+1}$
of this type comes from a linear flat extended connection $D$, and
that the connection $D$ is completely defined by $\wD$.  More
precisely, we have the following.

\begin{prop}\label{main}
Let $D:S^1 \to \B^1$ be a linear but not necessarily flat extended
connection on $M$, and let $\wD:\wB^\cdot \to \wB^{\cdot+1}$ be the
associated weakly Hodge derivation of the reduced Weil algebra
$\wB^\cdot$. Assume that $\wD \circ \wD = 0$. 

There exists a unique weakly Hodge bundle map $P:S^1 \to \I^1$ such
that the extended connection $D' = D + P:S^1 \to \B^1$ is linear
and flat.
\end{prop}

\punkt 
\proof Assume given a linear extended connection $D$ satisfying the
condition of Proposition~\ref{main}. To prove the proposition, we
have to construct a weakly Hodge map $P:S^1 \to \I^1$ such that the
extended connection $D+P$ is linear and flat. We do it by induction
on the augmentation degree, that is, we construct one-by-one the
terms $P_k$ in the augmentation degree decomposition $P = \sum_k
P_k$. The identity $\wD \circ \wD = 0$ is the base of the induction,
and the induction step is given by applying the following lemma to
$D + \sum_{i=2}^k P_k$, for each $k \geq 1$ in turn.

\begin{lemma}\label{main.ind}
Assume given a linear extended connection $D:S^1 \to \B^{\cdot+1}$
on $M$ and let $D:\B^\cdot \to \B^{\cdot+1}$ also denote the
associated derivation.  Assume also that the composition $D \circ
D:\B^\cdot \to \B^{\cdot+2}$ maps $S^1$ into $\I^2_{>k} =
\oplus_{p>k}\I^2_p$.

There exists a unique weakly Hodge bundle map $P_k:S^1 \to
\I^1_{k+1}$ such that the extended connection $D' = D + P_k:S^1 \to
\B^1$ is linear, and for the associated derivation $D':\B^\cdot \to
\B^{\cdot+1}$ the composition $D' \circ D'$ maps $S^1$ into
$\I^2_{>k+1} = \oplus_{p>k+1}\I^2_p$.
\end{lemma}

\proof Let $D:\B_{tot}^\cdot \to \B_{tot}^\cdot$ be the derivation of the
total Weil algebra associated to the extended connection $D$, and
let 
$$
R:\left(\B_{tot}^\cdot\right)_\cdot \to
\left(\B_{tot}^{\cdot+2}\right)_{\cdot+k}
$$ 
be the component of augmentation degree $k$ of the composition $D
\circ D:S^1 \to \B^2_{tot}$. Note that by \ref{deriv} the map $R$
vanishes on the subbundle $\Lambda^1_{tot} \subset
\left(\B^1_{tot}\right)_1$. Moreover, the composition $C \circ
R:\B^\cdot_{tot} \to \B^{\cdot+3}$ of the map $R$ with the canonical
derivation $C:\B^\cdot_{tot} \to \B^{\cdot+1}_{tot}$ vanishes on the
subbundle $S^1 \subset \B^0_{tot}$. Indeed, since $C$ maps $S^1$ into
$\Lambda^1_{tot}$, the composition $C \circ R$ is equal to the
commutator $[C,R]:S^1 \to \B^3_{tot}$.  Since the extended connection
$D$ is by assumption linear, we have $D_0=C$, and
\begin{align*}
C \circ R = [C,R] &= \sum_{0 \leq p \leq k}[C,D_p \circ D_{k-p}] =\\
&= [C,\{C,D_k\}] + \sum_{1 \leq p \leq k-1}[C,D_p \circ D_{k-p}].
\end{align*}
Since $C \circ C = 0$, the first term in the right hand side
vanishes. Let $\Theta = \sum_{1 \leq p \leq k-1}D_p:\B^\cdot_{tot} \to
\B^{\cdot+1}_{tot}$. Then the second term is the component of
augmentation degree $k$ in the commutator $[C,\Theta \circ
\Theta]:S^1 \to \B^3_{tot}$. By assumption $\{D,D\}=0$ in augmentation
degrees $<k$. Therefore we have $\{C,\Theta\} = - \{\Theta,\Theta\}$
in augmentation degrees $< k$. Since $\Theta$ increases the
augmentation degree, this implies that in augmentation degree $k$
$$
[C,\Theta \circ \Theta] = \{C,\Theta\} \circ \Theta - \Theta \circ
\{C,\Theta\} = [\Theta,\{\Theta,\Theta\}],
$$
which vanishes tautologically. 

The set of all weakly Hodge maps $P:S^1 \to \I_k^1$ coincides with
the set of all maps $P:S^1 \to \left(\I_{tot}^1\right)_k$ preserving the
Hodge bundle structures. Let $P$ be such a map, and let
$D':\B_{tot}^\cdot \to \B_{tot}^{\cdot+1}$ be the derivation associated to
the extended connection $D' = D + P$.

Since the extended connection $D$ is by assumption linear, by
Lemma~\ref{lin.aug} the extended connection $D'$ is linear if and
only if $\sigma_{tot} \circ P = 0$.  Moreover, since the augmentation
degree-$0$ component of the derivation $D$ equals $C$, the
augmentation degree-$k$ component $Q:S^1 \to \B^2_{tot}$ in the
composition $D' \circ D'$ is equal to
$$
Q = R + \{C,(D'-D)\}.
$$
By definition $D'-D:\B_{tot}^\cdot \to \B_{tot}^{\cdot+1}$ equals $P$ on
$S^1 \subset \B_{tot}^0$ and vanishes on $\Lambda^1_{tot} \subset
\B_{tot}^1$. Since $C$ maps $S^1$ into $\Lambda^1_{tot} \subset \B_{tot}^1$, we
have $Q = R + C \circ P$. Thus, a map $P$ satisfies the condition of
the lemma if and only if
$$
\begin{cases}
C \circ P = -R\\
\sigma_{tot} \circ P = 0
\end{cases}
$$
To prove that such a map $P$ is unique, note that these equations imply 
$$
h \circ P = (\sigma_{tot} \circ C + C \circ \sigma_{tot}) \circ P = \sigma_{tot}
\circ R, 
$$
and $h$ is invertible by Corollary~\ref{h.inv}. To prove that such a
map $P$ exists, define $P$ by 
$$
P = -h^{-1} \circ \sigma_{tot} \circ R:S^1 \to \I^1_{k+1}.
$$
The map $h = \{C,\sigma_{tot}\}$ and its inverse $h^{-1}$ commute with
$C$ and with $\sigma_{tot}$.  Since $\sigma_{tot} \circ \sigma_{tot} = C \circ
C0$, we have $\sigma_{tot} \circ P = h^{-1} \circ \sigma_{tot} \circ
\sigma_{tot} \circ R = 0$. On the other hand, $C \circ R = 0$. Therefore
\begin{align*}
C \circ P &= - C \circ h^{-1} \circ \sigma_C \circ R = - h^{-1} \circ
C \circ \sigma_{tot} \circ R = \\
&= h^{-1} \circ \sigma_{tot} \circ C \circ R - h^{-1} \circ h \circ R = -R.
\end{align*}
This finishes the proof of the lemma and of Proposition~\ref{main}.
\endproof

\subsection{Reduction of extended connections}\label{pf.last.sub}

\punkt We now complete the proof of Theorem~\ref{kal=ext}.  First we
will need to identify explicitly the low Hodge bidegree components
of the reduced Weil algebra $\wB^\cdot$. The following is easily
checked by direct inspection.

\begin{lemma}
We have 
\begin{align*}
&\wB^{2,-1} \oplus \wB^{1,0} \oplus \wB^{0,1} \oplus \wB^{-1,2} =
\Lambda^1 \oplus \left(S^1 \otimes \Lambda^1\right) \subset \wB^1\\ 
&\wB^{3,-2} \oplus \wB^{2,-1} \oplus \wB^{1,1} \oplus \wB^{-1,2}
\oplus \wB^{-2,3} = \\
&\qquad\qquad\qquad\qquad\qquad\qquad = \left(S^{1,-1} \otimes
\Lambda^{2,0}\right) \oplus \left(S^{-1,1} \otimes
\Lambda^{0,2}\right) \oplus \Lambda^2 \subset \wB^2   
\end{align*}
\end{lemma}

\punkt Let now $\nabla:S^1 \to S^1 \otimes \Lambda^1$ be an
arbitrary real connection on the bundle $S^1$. The operator 
$$
D = C + \nabla:S^1 \to \Lambda^1 \oplus \left(\Lambda^1 \otimes S^1\right) 
\subset \B^1 
$$
is then automatically weakly Hodge and defines therefore an extended
connection on $M$. This connection is linear by
Lemma~\ref{lin.aug}. Extend $D$ to a derivation $D:\B^\cdot \to
\B^{\cdot+1}$ as in \ref{deriv}, and let $\wD:\wB^\cdot \to
\wB^{\cdot+1}$ be the associated derivation of the reduced Weil
algebra. 

\begin{lemma}\label{kal=red}
The derivation $\wD$ satisfies $\wD \circ \wD = 0$ if and only if
the connection $\nabla$ is K\"ah\-le\-ri\-an. 
\end{lemma}

\proof Indeed, the operator $\wD \circ \wD$ is weakly Hodge, hence
factors through a bundle map
\begin{multline*}
\wD \circ \wD: S^1 \to \wB^{3,-1} \oplus \wB^{2,0} \oplus \wB^{1,1}
\oplus \wB^{0,2} \oplus \wB^{-1,3} = \\
= \left(S^{1,-1} \otimes \Lambda^{2,0}\right) \oplus \left(S^{-1,1}
\otimes \Lambda^{0,2}\right) \oplus \Lambda^2 \subset \B^2.  
\end{multline*}
By definition we have 
$$
\wD \circ \wD = (C + \nabla) \circ (C + \nabla) = \{C, \nabla\} + \{
\nabla, \nabla \}.
$$
An easy inspection shows that the sum is direct, and the first
summand equals 
$$
\{C, \nabla\} = T \circ C:S^1 \to \Lambda^2, 
$$
where $T$ is the torsion of the connection $\nabla$, while the
second summand equals
$$
\{\nabla, \nabla\} = R^{2,0} \oplus R^{0,2}: S^{1,-1} \oplus
S^{-1,1} \to \left(S^{1,-1} \otimes \Lambda^{2,0}\right) \oplus
\left(S^{-1,1} \otimes \Lambda^{0,2}\right), 
$$
where $R^{2,0}$, $R^{0,2}$ are the Hodge type components of the
curvature of the connection $\nabla$. Hence $\wD \circ \wD = 0$ if
and only if $R^{2,0} = R^{0,2} = T = 0$, which proves the lemma and
finishes the proof of Theorem~\ref{kal=ext}. 
\endproof

\punkt We finish this section with the following corollary of
Theorem~\ref{kal=ext} which gives an explicit expression for
the augmentation degree-$2$ component $D_2$ of a flat linear
extended connection $D$ on the manifold $M$. We will need this
expression in Section~\ref{metrics.section}.

\begin{corr}\label{D.2}
Let $D = \sum_{k \geq 0}D_k:S^1 \to \B^1$ be a flat linear extended
connection on $M$, so that $D_0 = C$ and $D_1$ is a K\"ah\-le\-ri\-an
connection on $M$. We have 
$$
D_2 = \frac{1}{3}\sigma \circ R, 
$$
where $\sigma:\B^{\cdot+1} \to \B^\cdot$ is the canonical derivation
introduced in \ref{sigma}, and $R = D_1 \circ D_1:S^1 \to S^1
\otimes \Lambda^{1,1} \subset \B^2$ is the curvature of the
K\"ah\-le\-ri\-an connection $D_1$.
\end{corr}

\proof Extend the connection $D$ to a derivation $D_C = \sum_{k \geq
0}D^{tot}_k:\B^\cdot_{tot} \to \B^{\cdot+1}_{tot}$ of the total Weil algebra. By
the construction used in the proof of Lemma~\ref{main.ind} we have
$D^{tot}_2 = h^{-1} \circ \sigma_{tot} \circ R_{tot}:S^1 \to
\left(\B^1_{tot}\right)_3$, where $h:\B^\cdot_{tot} \to \B^\cdot_{tot}$ is as in
Lemma~\ref{h.acts}, the map $\sigma_{tot}:\B^{\cdot+1}_{tot} \to \B^\cdot_{tot}$
is the canonical map constructed in \ref{sigma.c}, and $R^{tot}:S^1 \to
\left(\B^2_{tot}\right)_3$ is the square $R_{tot} = D^{tot}_1 \circ D^{tot}_1$ of
the derivation $D^{tot}_1:\B^\cdot_{tot} \to \B^{\cdot+1}_{tot}$. By
Lemma~\ref{h.on.b1} the map $h$ acts on $\left(\B_{tot}^1\right)_3$ by
multiplication by $3$. Moreover, it is easy to check that
$$
\left(B^2_{tot}\right)_3 = \left(S^1 \otimes \Lambda^2\right) \oplus
\left(S^{-1,1} \otimes \Lambda^{2,0}\right)\oplus
\left(S^{1,-1} \otimes \Lambda^{0,2}\right), 
$$
and the map $R^{tot}:S^1 \to \left(B^2_{tot}\right)_3$ sends $S^1$ into the
first summand in this decomposition and coincides with the curvature
$R:S^1 \to S^1 \otimes \Lambda^2 \subset \left(B^2_{tot}\right)_3$ of
the K\"ah\-le\-ri\-an connection $D_1$. Therefore $\sigma_{tot} \circ R^{tot} =
\sigma \circ R$, which proves the claim. 
\endproof

\section{Metrics}\label{metrics.section}

\subsection{Hyperk\"ahler metrics on Hodge manifolds}

\punkt Let $M$ be a complex manifold equipped with a K\"ah\-le\-ri\-an
connection $\nabla$, and consider the associated linear formal Hodge
manifold structure on the tangent bundle $\tm$. In this section we
construct a natural bijection between the set of all polarizations
on the Hodge manifold $\tm$ in the sense of
Subsection~\ref{polarization} and the set of all K\"ahler metrics on
$M$ compatible with the given connection $\nabla$.

\punkt \label{restr}
Let $h$ be a hyperk\"ahler metric on $\tm$, or, more
generally, a formal germ of such a metric in the neighborhood of the
zero section $M \subset \tm$. Assume that the metric $h$ is
compatible with the given hypercomplex structure and Hermitian-Hodge
in the sense of \ref{hermhodge}, and let $\omega_I$ be the K\"ahler
form associated to $h$ in the preferred complex structure $\tm_I$ on
$\tm$.

Let $h_M$ be the restriction of the metric $h$ to the zero section
$M \subset \tm$, and let $\omega \in C^\infty(M,\Lambda^{1,1}(M))$
be the associated real $(1,1)$-form on the complex manifold
$M$. Since the embedding $M \subset \tm_I$ is holomorphic, 
the form $\omega$ is the restriction
onto $M$ of the form $\omega_I$. In particular, it is closed, and
the metric $h_M$ is therefore K\"ahler.  

\punkt The main result of this section is the following. 

\begin{theorem}\label{metrics}
Restriction onto the zero section $M \subset \tm$ defines a
one-to-one correspondence between 
\begin{enumerate}
\item K\"ahler metrics on $M$ compatible with the K\"ah\-le\-ri\-an
connection $\nabla$, and 
\item formal germs in the neighborhood on $M \subset \tm$ of
Hermitian-Hodge hyperk\"ahler metrics on $\tm$ compatible with the
given formal Hodge manifold structure. 
\end{enumerate}
\end{theorem}

\noindent 
The rest of this section is devoted to the proof of this theorem. 

\punkt In order to prove Theorem~\ref{metrics}, we reformulate it in
terms of polarizations rather than metrics. Recall (see
\ref{positive}) that a {\em polarization} of the formal Hodge
manifold $\tm$ is by definition a $(2,0)$-form $\Omega \in
C^\infty_M(\tm,\Lambda^{2,0}(\tm_J))$ for the complementary complex
structure $\tm_J$ which is holomorphic, real and of $H$-type $(1,1)$
with respect to the canonical Hodge bundle structure on
$\Lambda^{2,0}(\tm_J)$, and satisfies a certain positivity condition
\eqref{P}.

\punkt By Lemma~\ref{pol.hm} Hermitian-Hodge hyperk\"ahler metrics
on $\tm$ are in one-to-one correspondence with polarizations. Let
$h$ be a metric on $M$, and let $\omega_I$ and $\omega$ be the
K\"ahler forms for $h$ on $\tm_I$ and on $M \subset \tm_I$.
The corresponding polarization $\Omega \in
C^\infty(\tm,\Lambda^{2,0}(\tm_J))$ satisfies by
\eqref{omega.and.Omega} 
\begin{equation}\label{pol2kal}
\omega_I = \half\left(\Omega + \nu(\Omega)\right) \in \Lambda^2(\tm,\C),
\end{equation}
where $\nu:\Lambda^2(\tm,\C) \to \overline{\Lambda^2(\tm,\C)}$ is
the usual complex conjugation. 

\punkt Let $\rho:\tm \to M$ be the natural projection, and let 
$$
\Res:\rho_*\Lambda^{\cdot,0}(\tm_J) \to \Lambda^\cdot(M,\C)
$$
be the map given by the restriction onto the zero section $M \subset
\tm$. Both bundles are naturally Hodge bundles of the same weight on
$M$ in the sense of \ref{hodge.bundles}, and the bundle map $\Res$
preserves the Hodge bundle structures. Since $\Omega$ is of $H$-type
$(1,1)$, the form $\Res\Omega \in C^\infty(M,\Lambda^2(M,\C))$ is
real and of Hodge type $(1,1)$. By \eqref{pol2kal} 
\begin{multline*}
\Res\Omega = \half\left(\Res\Omega + \overline{\Res\Omega}\right) =
\half(\Omega+\nu(\Omega))|_{M \subset \tm} = \\
= \omega_I|_{M \subset
\tm} = \omega \in \Lambda^{1,1}(M).
\end{multline*}
Therefore to prove Theorem~\ref{metrics}, it suffices to prove the
following.  
\begin{itemize}
\item \label{to.prove} For every polarization $\Omega$ of the formal
Hodge manifold $\tm$ the restriction $\omega = \Res\Omega \in
C^\infty(\Lambda^{1,1}(M))$ is compatible with the connection
$\nabla$, that is, $\nabla\omega=0$. Vice versa, every real positive
$(1,1)$-form $\omega \in C^\infty(\Lambda^{1,1}(M))$ satisfying
$\nabla\omega=0$ extends to a polarization $\Omega$ of $\tm$, and
such an extension is unique.
\end{itemize}
This is what we will actually prove. 

\subsection{Preliminaries}\label{drm.mod}

\punkt\label{descr} We begin with introducing a convenient model for
the holomorphic de Rham algebra $\Lambda^{\cdot,0}(\tm_J)$ of the
complex manifold $\tm_J$, which would be independent of the Hodge
manifold structure on $\tm$.  To construct such a model, consider
the relative de Rham complex $\Lambda^\cdot(\tm/M,\C)$ of $\tm$ over
$M$ (see \ref{relative.de.rham.sub} for a reminder of its definition
and main properties). Let $\pi:\Lambda^\cdot(\tm,\C) \to
\Lambda^\cdot(\tm/M,\C)$ be the canonical projection. Recall that
the bundle $\Lambda^i(\tm/M,\C)$ of relative $i$-forms on $\tm$ over
$M$ carries a natural structure of a Hodge bundle of weight
$i$. Moreover, we have introduced in \eqref{eta} a Hodge bundle
isomorphism
$$
\eta:\rho^*\Lambda^\cdot(M,\C) \to \Lambda^\cdot(\tm/M,\C) 
$$
between $\Lambda^i(\tm/M,\C)$ and the pullback
$\rho^*\Lambda^\cdot(M,\C)$ of the bundle of $\C$-valued $i$-forms
on $M$.

\begin{lemma}\label{ident.bis}
\begin{enumerate}
\item The projection $\pi$ induces an algebra isomorphism 
$$
\pi:\Lambda^{\cdot,0}(\tm_J) \to \Lambda^\cdot(\tm/M,\C) 
$$
compatible with the natural Hodge bundle structures. 
\item Let $\alpha \in C^\infty_M(\tm,\Lambda^{i,0}(\tm_J))$ be a
smooth $(i,0)$-form on $\tm_J$, and consider the smooth $i$-form 
$$
\beta = \eta^{-1}\pi(\alpha) \in C^\infty_M(\tm,\Lambda^i(\tm,\C)) 
$$
on $\tm$. The forms $\alpha$ and $\beta$ have the same restriction
to the zero section $M \subset \tm$. 
\end{enumerate}
\end{lemma}

\proof Since $\eta$, $\pi$ and the restriction map are compatible
with the algebra structure on $\Lambda^\cdot(M,\C)$, it suffices to
prove both claims for $\Lambda^1(M,\C)$. For every bundle $\E$ on
$\tm$ denote by $\E|_{M \subset \tm}$ the restriction of $\E$ to the
zero section $M \subset \tm$. Consider the bundle map
$$
\chi = \eta \circ \Res:\Lambda^{1,0}(\tm_J)|_{M \subset \tm} \to
\Lambda^1(M,\C) \to \Lambda^1(\tm/M,\C)|_{M \subset \tm}. 
$$
The second claim of the lemma is then equivalent to the identity
$\chi = \pi$.\! Moreover, note that the contraction with
the canonical vector field $\phi$ on $\tm$ defines an injective map
$i_\phi:C^\infty(M,\Lambda^1(\tm/M,\C)|_M) \to
C^\infty(\tm,\C)$. Therefore it suffices to prove that $i_\phi \circ
\chi = i_\phi \circ \pi$.

Every smooth section $s$ of the bundle $\Lambda^{1,0}(\tm_J)|_{M \subset \tm}$
is of the form 
$$
s = (\rho^*\alpha + \sqrt{-1}j\rho^*\alpha)|_{M \subset \tm}, 
$$
where $\alpha \in C^\infty(M,\Lambda^1(M,\C))$ is a smooth $1$-form
on $M$, and $j:\Lambda^1(\tm,\C) \to \overline{\Lambda^1(\tm,\C)}$
is the map induced by the quaternionic structure on $\tm$. For such
a section $s$ we have $\Res(s) = \alpha$, and by \eqref{phi.and.tau}
$i_\phi(\chi(s)) = \sqrt{-1}\tau(\alpha)$, where
$\tau:C^\infty(M,\Lambda^1(M,\C)) \to C^\infty(\tm,\C)$ is the
tautological map introduced in \ref{tau}. On the other hand, since $\pi
\circ \rho^* = 0$, we have 
$$
i_\phi(\pi(s)) = i_\phi(\pi(\sqrt{-1}j\rho^*\alpha)) =
i_\phi(\sqrt{-1}j\rho^*\alpha).
$$ 
Since the Hodge manifold structure on $\tm$ is linear, this equals 
$$
i_\phi(\pi(s)) = \sqrt{-1}i_\phi(j(\rho^*\alpha)) =
\sqrt{-1}\tau(\alpha) = \chi(s), 
$$
which proves the second claim of the lemma. Moreover, it shows that
the restriction of the map $\pi$ to the zero section $M \subset \tm$
is an isomorphism. As in the proof of Lemma~\ref{ident}, this
implies that the map $\pi$ is an isomorphism on the whole $\tm$,
which proves the first claim and finishes the proof of the lemma. 
\endproof 

\punkt Lemma~\ref{ident.bis}~\thetag{i} allows to define the bundle
isomorphism
$$
\pi^{-1} \circ \eta:\rho^*\Lambda^\cdot(M,\C) \to
\Lambda^\cdot(\tm/M,\C) \to \Lambda^{\cdot,0}(\tm_J), 
$$
between $\rho^*\Lambda^\cdot(M,\C)$ and $\Lambda^{\cdot,0}(\tm_J)$,
and it induces an isomorphism
$$
\rho_*(\eta \circ \pi^{-1}):\rho_*\rho^*\Lambda^\cdot(M,\C) \cong
\rho_*\Lambda^{\cdot,0}(\tm_J)
$$
between the direct images of these bundles under the canonical
projection $\rho:\tm \to M$. 

On the other hand, by adjunction we have the canonical embedding
$$
\Lambda^\cdot(M,\C) \hookrightarrow \rho_*\rho^*\Lambda^\cdot(M,\C),
$$
and by the projection formula it extends to an isomorphism
$$
\rho_*\rho^*\Lambda^\cdot(M,\C) \cong \Lambda^\cdot(M,\C) \otimes
\B^0, 
$$
where $\B^0 = \rho_*\Lambda^0(\tm,\C)$ is the $0$-th component of
the Weil algebra $\B^\cdot$ of $M$. All these isomorphisms are
compatible with the Hodge bundle structures and with the
multiplication.

\punkt \label{ident.punkt}
It will be convenient to denote the image $\rho_*(\eta \circ
\pi^{-1})\left(\Lambda^\cdot(M,\C)\right) \subset
\rho_*\Lambda^{\cdot,0}(\tm_J)$ by $L^\cdot(M,\C)$ or, to simplify
the notation, by $L^\cdot$. (The algebra $L^\cdot(M,\C)$ is, of
course, canonically isomorphic to $\Lambda^\cdot(M,\C)$.) We then
have the identification
\begin{equation}\label{ident.zero}
L^\cdot \otimes \B^0 \cong \rho_*\rho^*\Lambda^\cdot(M,\C) \cong 
\rho_*\Lambda^{\cdot,0}(\tm_J). 
\end{equation}
This identification is independent of the Hodge manifold structure
on $\tm$. Moreover, by Lemma~\ref{ident.bis}~\thetag{ii} the
restriction map $\Res:\rho_*\Lambda^{\cdot,0}(\tm_J) \to
\Lambda^\cdot(M,\C)$ is identified under \eqref{ident.zero} with the
canonical projection $L^\cdot \otimes \B^0 \to L^\cdot \otimes
\B^0_0 \cong L^\cdot$.

By Lemma~\ref{ident} we also have the identification
$\rho_*\Lambda^{0,\cdot}(\tm_J) \cong \B^\cdot$. Therefore
\eqref{ident.zero} extends to an algebra isomorphism
\begin{equation}\label{ident.formula}
\rho_*\Lambda^{\cdot,\cdot}(\tm_J) \cong
\rho_*\Lambda^\cdot(\tm/M,\C) \otimes \Lambda^\cdot(M,\C) \cong
L^\cdot \otimes \B^\cdot. 
\end{equation}
This isomorphism is also compatible with the Hodge bundle structures
on both sides. 

\subsection{The Dolbeault differential on $\protect\tm_J$}

\punkt Our next goal is to express the Dolbeault differential
$\bar\6_J$ of the complex manifold $\tm_J$ in terms of the model for
the de Rham complex $\Lambda^{\cdot,\cdot}(\tm_J)$ given by
\eqref{ident.formula}. For every $k \geq 0$ denote by
$$
D:L^k \otimes \B^\cdot \to L^k \otimes \B^{\cdot+1}. 
$$
the differential operator induced by
$\bar\6_J:\Lambda^{\cdot,\cdot}(\tm_J) \to
\Lambda^{\cdot,\cdot+1}(\tm_J)$ under \eqref{ident.formula}. The
operator $D$ is weakly Hodge. It satisfies the Leibnitz rule with
respect to the algebra structure on $L^\cdot \otimes \B^\cdot$, and
we have $D \circ D = 0$. By definition for $k=0$ it coincides with
the derivation $D:\B^\cdot \to \B^{\cdot+1}$ defined by the Hodge
manifold structure on $\tm$. For $k > 0$ the complex $\langle L^k
\otimes \B^\cdot, D \rangle$ is a free differential graded module
over the Weil algebra $\langle\B^\cdot,D\rangle$.

\punkt \label{dr.L}
The relative de Rham differential $d^r$ (see
Subsection~\ref{relative.de.rham.sub}) induces under the isomorphism
\eqref{ident.formula} an algebra derivation
$$
d^r:L^\cdot \otimes \B^\cdot \to L^{\cdot+1} \otimes \B^\cdot. 
$$
The derivation $d^r$ also is weakly Hodge, and we have the
following. 

\begin{lemma}\label{D.dr}
The derivations $D$ and $d^r$ commute, that is,  
$$
\{D,d^r\} = 0:L^\cdot \otimes \B^\cdot \to L^{\cdot+1} \otimes
\B^{\cdot+1}. 
$$
\end{lemma}

\proof The operator $\{D,d^r\}$ satisfies the Leibnitz rule, so it
suffices to prove that it vanishes on $\B^0$, $\B^1$ and $L^1
\otimes \B^0$.  Moreover, the $\B^0$-modules $\B^1$ and $L^1 \otimes
\B^0$ are generated, respectively, by local sections of the form $Df$ and
$d^rf$, $f \in \B^0$. Since $\{D,d^r\}$ commutes with $D$ and $d^r$,
it suffices to prove that it vanishes on $\B^0$. Finally,
$\{D,d^r\}$ is continuous in the adic topology on $\B^0$. Since the
subspace
$$
\{fg|f,g \in \B^0, Df = d^rg = 0\} \subset \B^0
$$
is dense in this topology, it suffices to prove that for a local
section $f \in \B^0$ we have $\{D,d^r\}f=0$ if either $d^rf=0$ of
$Df=0$.

It is easy to see that for every local section $\alpha \in \B^\cdot$
we have $d^r\alpha = 0$ if and only if $\alpha \in \B^\cdot_0$ is of
augmentation degree $0$. By definition the derivation $D$ preserves
the component $\B^\cdot_0 \subset \B^\cdot$ of augmentation degree
$0$ in $\B^0$. Therefore $d^rf=0$ implies $d^rDf=0$ and consequently
$\{D,d^r\}f=0$. This handles the case $d^rf=0$. To finish the proof,
assume given a local section $f \in \B^0$ such that $Df=0$. Such a
section by definition comes from a germ at $M \subset \tm$ of a
holomorphic function $\wh{f}$ on $\tm_J$. Since $\wh{f}$ is
holomorphic, we have $\bar\6_J\wh{f}=0$ and $d\wh{f} =
\6_J\wh{f}$. Therefore $d^rf = \pi(d\wh{f}) = \pi(\6_Jf)$, and
$$
Dd^rf = \pi(\bar\6_J\6_J\wh{f}) = -\pi(\6_J\bar\6_J\wh{f}) = 0, 
$$
which, again, implies $\{D,d^r\}f = 0$. 
\endproof

\punkt Let now $\nabla=D_1:S^1 \to S^1 \otimes \Lambda^1$ be the
reduction of the extended connection $D$. It induces a connection on
the bundle $L^1 \cong S^1$, and this connection extends by the
Leibnitz rule to a connection on the exterior algebra $L^\cdot$ of
the bundle $L^1$, which we will also denote by $\nabla$.

Denote by $R = \nabla \circ \nabla:L^\cdot \to L^\cdot \otimes
\Lambda^2$ the curvature of the connection $\nabla$. Since $\nabla
\circ \nabla = \half\{\nabla,\nabla\}$, the operator $R$ also
satisfies the Leibnitz rule with respect to the multiplication in
$L^\cdot$. 

\punkt Introduce the augmentation grading on the bundle $L^\cdot
\otimes \B^\cdot$ by setting $\deg L^\cdot = 0$. The derivation
$D:L^\cdot \otimes \B^\cdot \to L^\cdot \otimes \B^{\cdot+1}$
obviously does not increase the augmentation degree, and we have the
decomposition $D = \sum_{k \geq 0}D_k$. On the other hand, the
derivation $d^r$ preserves the augmentation degree. Therefore
Lemma~\ref{D.dr} implies that for every $k \geq 0$ we have
$\{D_k,d^r\} = 0$. This in turn implies that $D_0 = 0$ on $L^p$ for
$p > 0$, and therefore $D_0 = \id \otimes C:L^p \otimes \B^\cdot \to
L^p \otimes \B^{\cdot+1}$. Moreover, this allows to identify
explicitly the components $D_1$ and $D_2$ of the derivation
$D:L^\cdot \to L^\cdot \otimes \B^1$. Namely, we have the following.

\begin{lemma}\label{D.1.2}
We have 
\begin{align*}
D_1 &= \nabla:L^\cdot \to L^\cdot \otimes \B^1_1 = L^\cdot \otimes
\Lambda^1\\ 
D_2 &= \frac{1}{3}\sigma \circ R:L^\cdot \to L^\cdot \otimes \B^1_2,
\end{align*}
where $\sigma = \id \otimes \sigma:L^\cdot \otimes \B^{\cdot+1} \to
L^\cdot \otimes \B^\cdot$ is as in \ref{C.and.sigma}. 
\end{lemma}

\proof Since both sides of these identities satisfy the Leibnitz
rule with respect to the multiplication in $L^\cdot$, it suffices to
prove them for $L^1$. But $d^r:\B^0 \to L^1 \otimes \B^0$ restricted
to $S^1 \subset \B^0$ becomes an isomorphism $d^r:S^1 \to
L^1$. Since $\{D_1,d^r\} = \{D_2,d^r\} = 0$, it suffices to prove
the identities with $L^1$ replaced with $S^1$. The first one then
becomes the definition of $\nabla$, and the second one is
Corollary~\ref{D.2}.  \endproof

\subsection{The proof of Theorem~\ref{metrics}}

\punkt We can now prove Theorem~\ref{metrics} in the form
\ref{to.prove}. We begin with the following corollary of
Lemma~\ref{D.1.2}.

\begin{corr}\label{red.corr}
Let $\I^\cdot \subset \B^\cdot$ be the ideal introduced in
\ref{red.Weil}. An arbitrary smooth section $\alpha \in
C^\infty(M,L^\cdot)$ satisfies
\begin{equation}\label{rrd}
D\alpha \in C^\infty(M, L^\cdot \otimes \I^1) \subset C^\infty(M,
L^\cdot \otimes \B^1)
\end{equation}
if and only if $\nabla \alpha = 0$. 
\end{corr}

\proof Again, both the identity \eqref{rrd} and the equality
$\nabla\alpha$ are compatible with the Leibnitz rule with respect to
the multiplication in $\alpha$. Therefore it suffices to prove that
they are equivalent for every $\alpha \in L^1$. By definition of the
ideal $\I^\cdot$ the equality \eqref{rrd} holds if and only if
$D_1\alpha = D_2\alpha = 0$. By Lemma~\ref{D.1.2} this is equivalent
to $\nabla\alpha = \sigma \circ R(\alpha) = 0$. But since $R =
\nabla \circ \nabla$, $\nabla\alpha = 0$ implies $\sigma \circ
R(\alpha) = 0$, which proves the claim.  
\endproof

\punkt Let now $\Omega \in C^\infty(M, \rho_*\Lambda^{2,0}(\tm_J)
\cong L^2 \otimes \B^0$ be a polarization of the Hodge manifold
$\tm_J$, so that $\Omega$ is of Hodge type $(1,1)$ and
$D\Omega=0$. Let $\omega = \Res\Omega \in
C^\infty(M,\Lambda^{1,1}(M,\C))$ be its restriction, and let $\Omega
= \sum_{k \geq 0}\Omega_k$ be its augmentation degree decomposition.

As noted in \ref{ident.punkt}, the restriction map
$\Res:\rho_*\Lambda^{\cdot,0}(\tm_J) \to \Lambda^\cdot(M,\C)$ is
identified under the isomorphism \eqref{ident.formula} with the
projection $\LL^\cdot \otimes \B^0 \to L^\cdot$ onto the component of
augmentation degree $0$. Therefore $\omega = \Omega_0$. Since the
augmentation degree-$1$ component $\left(L^2 \otimes
B^0_{tot}\right)^{1,1}_1=0$ and $D\Omega = 0$, we have $\nabla\omega =
D_1\Omega_0 = 0$, which proves the first claim of
Theorem~\ref{metrics}.

\punkt To prove the second claim of the theorem, let $\omega$ be a
K\"ahler form on $M$ compatible with the connection $\nabla$, so
that $\nabla\omega = 0$. We have to show that there exists a unique
section $\Omega = \sum_{k \geq 0}\Omega_k \in C^\infty(M, L^2
\otimes \B^0)$ which is of Hodge type $(1,1)$ and satisfies
$D\Omega=0$ and $\Omega_0 = \omega$. As in the proof of
Theorem~\ref{kal=ext}, we will use induction on $k$. Since $\Omega$
is of Hodge type $(1,1)$, we must have $\Omega_1 = 0$, and by
Corollary~\ref{red.corr} we have $D(\Omega_0 + \Omega_1) \in
C^\infty(M, L^2 \otimes \I^1)$, which gives the base of our
induction. The induction step is given by applying the following
proposition to $\sum_{0 \geq p \geq k}\Omega_k$ for each $k \geq 1$
in turn. 

\begin{prop}\label{metrics.ind}
Assume given integers $p,q,k; p,q \geq 0, k \geq 1$ and assume given
a section $\alpha \in C^\infty(M, L^{p+q} \otimes \B^0)$ of Hodge
type $(p,q)$ such that
$$
D\alpha \in C^\infty(M, L^{p+q} \otimes \I^1_{\geq k}), 
$$
where $\I^\cdot_{\geq k} = \oplus_{m \geq k}\I^\cdot_m$. Then there
exists a unique section $\beta \in C^\infty(M, L^{p+q} \otimes
\B_k)$ of the same Hodge type $(p,q)$ and such that $D(\alpha +
\beta) \in C^\infty(M, L^{p+q} \otimes \I^1_{\geq k+1})$. 
\end{prop}

\proof Let $\B_{tot}^\cdot$ be the total Weil algebra introduced in
\ref{total.Weil}, and consider the free module $L^{p+q} \otimes
\B_{tot}^\cdot$ over $\B_{tot}^\cdot$ generated by the Hodge bundle
$L^{p+q}$. This module carries a canonical Hodge bundle structure of
weight $p+q$.  Consider the maps $C:\B_{tot}^\cdot \to \B_{tot}^{\cdot+1}$,
$\sigma_{tot}:\B^{\cdot+1} \to \B_{tot}^\cdot$ introduced in \ref{total.C}
and \ref{sigma.c}, and let $C = \id \otimes C,\sigma_{tot} = \id \otimes
\sigma_{tot}:L^{p+q} \otimes \B_{tot} \to L^{p+q} \otimes \B_{tot}$ be the
associated endomorphisms of the free module $L^{p+q} \otimes
\B_{tot}^\cdot$. 

The maps $C$ and $\sigma_{tot}$ preserve the Hodge bundle structure.
The commutator $h = \{C, \sigma_{tot}\}: \B^\cdot_{tot} \to \B^\cdot_{tot}$ is
invertible on $\I_{tot}^\cdot \subset \B_{tot}^\cdot$ by
Corollary~\ref{h.inv} and acts as $k\id$ on $\B^0_k \subset
\B^\cdot_{tot}$. Therefore the endomorphism
$$
h = \id \otimes h = \{C,\sigma_{tot}\}:L^{p+q} \otimes \B^\cdot_{tot} \to
L^{p+q} \otimes \B^\cdot_{tot}
$$
is invertible on $L^{p+q} \otimes \I^\cdot_{tot}$ and acts as $k\id$ on
$L^{p+q} \otimes \B^0_k$.

Since the derivation $D:L^{p+q} \otimes \B^\cdot \to L^{p+q} \otimes
\B^{\cdot+1}$ is weakly Hodge, it induces a map $D^{tot}:L^{p+q} \otimes
\B_{tot}^\cdot \to L^{p+q} \otimes \B_{tot}^{\cdot+1}$, and $D\alpha \in
L^{p+q} \otimes \I^1_{\geq k}$ if and only if the same holds for
$D^{tot}\alpha$. To prove uniqueness, note that $D_0 = C$ is injective
on $L^{p+q} \otimes \B^0_k$. If there are two sections
$\beta,\beta'$ satisfying the conditions of the proposition, then
$D_0(\beta-\beta')=0$, hence $\beta = \beta'$. 

To prove existence, let $\gamma = (D^{tot}\alpha)_k$ be the component of
the section $D^{tot}\alpha$ of augmentation degree $k$. Since $D^{tot} \circ
D^{tot}=0$, we have $C\gamma = D^{tot}_0\gamma = 0$ and $C \circ \sigma_{tot}
\gamma = h \gamma$.  Let $\beta = - \frac{1}{k} \circ
\sigma_{tot}(\gamma)$. The section $\beta$ is of Hodge type $(p,q)$ and
of augmentation degree $k$. Moreover, it satisfies
$$
D^{tot}_0\beta = C\beta = -Ch^{-1}\sigma_{tot}(\gamma) = -h^{-1} \circ C
\circ \sigma_{tot} \gamma = -h^{-1} \circ h\gamma = -\gamma.
$$
Therefore $D_{tot}(\alpha+\beta)$ is indeed a section of $L^{p+q}
\otimes (I^1_{tot})_{\geq k+1}$, which proves the proposition and
finishes the proof of Theorem~\ref{metrics}. 
\endproof 

\subsection{The cotangent bundle}

\punkt For every K\"ahler manifold $M$ Theorem~\ref{metrics}
provides a canonical formal hyperk\"ahler structure on the total
space $\tm$ of the complex-conjugate to the tangent bundle to
$M$. In particular, we have a closed holomorphic $2$-form $\Omega_I$
for the preferred complex structure $\tm_I$ on $M$. 

Let $T^*M$ be the total space to the cotangent bundle to $M$
equipped with the canonical holomorphic symplectic form $\Omega$. To
obtain a hyperk\"ahler metric of the formal neighborhood of the zero
section $M \subset T^*M$, one can apply an appropriate version of
the Darboux Theorem, which gives a local symplectic isomorphism
$\kappa:\tm \to T^*M$ in a neighborhood of the zero
section. However, this theorem is not quite standard in the
holomorphic and formal situations. For the sake of completeness, we
finish this section with a sketch of a construction of such an
isomorphism $\kappa:\tm \to T^*M$ which can be used to obtain a
hyperk\"ahler metric on $T^*M$ rather than on $\tm$.

\punkt We begin with the following preliminary fact on the
holomorphic de Rham complex of the manifold $\tm_I$. 

\begin{lemma}\label{de.rham.exact}
Assume given either a formal Hodge manifold structure on the
$U(1)$-manifold $\tm$ along the zero section $M \subset \tm$, or an
actual Hodge manifold structure on an open neighborhood $U \subset
\tm$ of the zero section. 
\begin{enumerate}
\item For every point $m \in M$ there exists an open neighborhood $U
\subset \tm_I$ such that the spaces $\Omega^\cdot(U)$ of holomorphic
forms on the complex manifold $\tm_I$ (formally completed along $M
\subset \tm$ if necessary) equipped with the holomorphic de Rham
differential $\6_I:\Omega^\cdot(U) \to \Omega^{\cdot+1}(U)$ form an
exact complex. 
\item If the subset $U \subset \tm$ is invariant under the
$U(1)$-action on $\tm$, then the same is true for the subspaces
$\Omega^\cdot_k(U) \subset \Omega^\cdot(U)$ of forms of weight $k$
with respect to the $U(1)$-action. 
\item Assume further that the canonical projection $\rho:\tm_I \to
M$ is holomorphic for the preferred complex structure $\tm_I$ on
$\tm$. Then both these claims hold for the spaces
$\Omega^\cdot(U/M)$ of relative holomorphic forms on $U$ over $M$.
\end{enumerate}
\end{lemma}

\proof The claim \thetag{i} is standard. To prove \thetag{ii}, note
that, both in the formal and in the analytic situation, the spaces
$\Omega^\cdot(U)$ are equipped with a natural topology. Both this
topology and the $U(1)$-action are preserved by the holomorphic
Dolbeault differential $\6_I$.

The subspaces $\Omega^\cdot_{fin}(U) \subset \Omega^\cdot(U)$ of
$U(1)$-finite vectors are dense in the natural topology. Therefore
the complex $\langle \Omega^\cdot_{fin}(U),\6_I\rangle$ is also
exact. Since the group $U(1)$ is compact, we have
$$
\Omega^\cdot_{fin}(U) = \bigoplus_k\Omega^\cdot_k(U),
$$
which proves \thetag{ii}. The claim \thetag{iii} is, again,
standard. 
\endproof 

\punkt We can now formulate and prove the main result of this
subsection.

\begin{prop}
Assume given a formal polarized Hodge manifold structure on the
manifold $\tm$ along the zero section $M \subset \tm$ such that the
canonical projection $\rho:\tm_I \to M$ is holomorphic for the
preferred complex structure $\tm_I$ on $\tm$. Let $\Omega_I$ be the
associated formal holomorphic $2$-form on $\tm_I$. Let $T^*M$ be the
total space of the cotangent bundle to the manifold $M$ equipped
with a canonical holomorphic symplectic form $\Omega$. There exists
a unique $U(1)$-equivariant holomorphic map $\kappa:\tm_I \to T^*M$,
defined in a formal neighborhood of the zero section, which commutes
with the canonical projections onto $M$ and satisfies $\Omega_I =
\kappa^*\Omega$. Moreover, if the polarized Hodge manifold structure
on $\tm_I$ is defined in an open neighborhood $U \subset \tm$ of the
zero section $M \subset \tm$, then the map $\kappa$ is also defined
in a (possibly smaller) open neighborhood of the zero section.
\end{prop}

\proof By virtue of the uniqueness, the claim is local on $M$, so
that we can assume that the whole $M$ is contained in a
$U(1)$-invariant neighborhood $U \subset \tm_I$ satisfying the
conditions of Lemma~\ref{de.rham.exact}. Holomorphic maps $\kappa:U
\to T^*M$ which commute with the canonical projections onto $M$ are
in a natural one-to-one correspondence with holomorphic sections
$\alpha$ of the bundle $\rho^*\Lambda^{1,0}(M)$ on $\tm_I$. Such a
map $\kappa$ is $U(1)$-equivariant if and only if the corresponding
$1$-form $\alpha \in \Omega^1(U)$ is of weight $1$ with respect to
the $U(1)$-action. Moreover, it satisfies $\kappa^*\Omega=\Omega_I$
if and only if $\6_I\alpha = \Omega_I$. Therefore to prove the
formal resp. analytic parts of the proposition it suffices to prove
that there exists a unique holomorphic formal resp. analytic section
$\alpha \in C^\infty(U,\rho^*\Lambda^{1,0}(M)$ which is of weight
$1$ with respect to the $U(1)$-action and satisfies
$\6_I\alpha=\Omega_I$.

The proof of this fact is the same in the formal and in the analytic
situations. By definition of the polarized Hodge manifold the
$2$-form $\Omega_I \in \Omega^2(U)$ is of weight $1$ with respect to
the $U(1)$-action. Therefore by
Lemma~\ref{de.rham.exact}~\thetag{ii} there exists a holomorphic
$1$-form $\alpha \in \Omega^1(U)$ of weight $1$ with respect to the
$U(1)$-action and such that $\6_I\alpha = \Omega_I$. Moreover, the
image of the form $\Omega_I$ under the canonical projection
$\Omega^2(U) \to \Omega^2(U/M)$ is zero. Therefore by
Lemma~\ref{de.rham.exact}~\thetag{iii} we can arrange so that the
image of the form $\alpha$ under the projection $\Omega^1(U) \to
\Omega^1(U/M)$ is also zero, so that $\alpha$ is in fact a section
of the bundle $\rho^*\Lambda^{1,0}(M)$. This proves the existence
part. To prove uniqueness, note that every two such $1$-forms must
differ by a form of the type $\6_if$ for a certain holomorphic
function $f \in \Omega^0(U)$. Moreover, by
Lemma~\ref{de.rham.exact}~\thetag{ii} we can assume that the
function $f$ is of weight $1$ with respect to the $U(1)$-action.  On
the other hand, by Lemma~\ref{de.rham.exact}~\thetag{iii} we can
assume that the function $f$ is constant along the fibers of the
canonical projection $\rho:\tm_I \to M$. Therefore we have $f=0$
identically on the whole $U$.  
\endproof

\section{Convergence}\label{convergence}

\subsection{Preliminaries}

\punkt Let $M$ be a complex manifold. By Theorem~\ref{kal=ext} every
K\"ah\-le\-ri\-an connection $\nabla:\Lambda^1(M,\C) \to \Lambda^1(M,\C)
\otimes \Lambda^1(M,\C)$ on the cotangent bundle $\Lambda^1(M,\C)$
to the manifold $M$ defines a flat linear extended connection
$D:S^1(M,\C) \to \B^1(M,\C)$ on $M$ and therefore a formal Hodge
connection $D$ on the total space $\tm$ of the complex-conjugate to
the tangent bundle to $M$. By Proposition~\ref{equiv.bis}, this
formal Hodge connection defines, in turn, a formal Hodge manifold
structure on $\tm$ in the formal neighborhood of the zero section $M
\subset \tm$.

In this section we show that if the K\"ah\-le\-ri\-an connection $\nabla$
is real-analytic, then the corresponding formal Hodge manifold
structure on $\tm$ is the completion of an actual Hodge manifold
structure on an open neighborhood $U \subset \tm$ of the zero
section $M \subset \tm$. We also show that if the connection
$\nabla$ comes from a K\"ahler metric $\omega$ on $M$, then the
corresponding polarization $\Omega$ of the formal Hodge manifold
$\tm$ defined in Theorem~\ref{metrics} converges in a neighborhood
$U' \subset U \subset \tm$ of the zero section $M \subset \tm$ to a
polarization of the Hodge manifold structure on $U'$. Here is the
precise formulation of these results. 

\begin{theorem}\label{converge}
Let $M$ be a complex manifold equipped with a real-analytic
K\"ah\-le\-ri\-an connection $\nabla:\Lambda^1(M,\C) \to \Lambda^1(M,\C)
\otimes \Lambda^1(M,\C)$ on its cotangent bundle
$\Lambda^1(M,\C)$. There exists an open neighborhood $U \subset \tm$
of the zero section $M \subset \tm$ in the total space $\tm$ of the
complex-conjugate to the tangent bundle to $M$ and a Hodge manifold
structure on $U \subset \tm$ such that its completion along the zero
section $M \subset \tm$ defines a linear flat extended connection
$D$ on $M$ with reduction $\nabla$. 

Moreover, assume that $M$ is equipped with a K\"ahler metric
$\omega$ such that $\nabla\omega=0$, and let $\Omega \in
C^\infty_M(\tm,\Lambda^2(\tm,\C))$ be the formal polarization of the 
Hodge manifold structure on $U \subset \tm$ along $M \subset
\tm$. Then there exists an open neighborhood $U' \subset U$ of $M
\subset U$ such that $\Omega \in C^\infty(U',\Lambda^2(\tm,\C))
\subset C^\infty_M(\tm,\Lambda^2(\tm,\C))$. 
\end{theorem}

\punkt \label{formal.Weil}
We begin with some preliminary observations. First of all, the
question is local on $M$, therefore we may assume that $M$ is an
open neighborhood of $0$ in the complex vector space $V = \C^n$. Fix
once and for all a real structure and an Hermitian metrics on the
vector space $V$, so that it is isomorphic to its dual $V \cong
V^*$.

The subspace $\J \subset C^\infty(M,\C)$ of functions vanishing at
$0 \in M$ is an ideal in the algebra $C^\infty(M,\C)$, and $\J$-adic
topology on $C^\infty(M,\C)$ extends canonically to the de Rham
algebra $\Lambda^\cdot(M,\C)$ of the manifold $M$ and, further, to
the Weil algebra $\B^\cdot(M,\C)$ of $M$ introduced in
\ref{Weil.defn}. Instead of working with bundle algebra
$\B^\cdot(M,\C)$ on $M$, it will be convenient for us to consider
the vector space
$$
\B^\cdot = C^\infty_\J(M,B^\cdot(M,\C)),
$$ 
which is by definition the $\J$-adic completion of the space
$C^\infty(M,\B^\cdot(M,\C))$ of global sections of the Weil
algebra. This vector space is canonically a (pro-)algebra over
$\C$. Moreover, the Hodge bundle structure on $\B^\cdot(M,\C)$
induces an $\R$-Hodge structure on the algebra $\B^\cdot$. 

\punkt The $\J$-adic completion $C^\infty_\J(M,\C)$ of the space of
smooth functions on $M$ is canonically isomorphic to the completion
$\wh{S}^\cdot(V)$ of the symmetric algebra of the vector space $V
\cong V^*$. The cotangent bundle $\Lambda^1(M,\C)$ is isomorphic to 
the trivial bundle $\V$ with fiber $V$ over $M$, and the completed de Rham
algebra $C^\infty_\J(M,\Lambda^\cdot(M,\C))$ is isomorphic to the
product 
$$
C^\infty_\J(M,\Lambda^\cdot(M,\C)) \cong S^\cdot(V) \otimes
\Lambda^\cdot(V). 
$$
This is a free graded-commutative algebra generated by two copies of
the vector space $V$, which we denote by $V_1 = V \subset
\Lambda^1(V)$ and by $V_2 = V \subset S^1(V)$. It is convenient to
choose the trivialization $\Lambda^1(M,\C) \cong \V$ in such a way
that the de Rham differential $d_M:\Lambda^\cdot(M,\C) \to
\Lambda^{\cdot+1}(M,\C)$ induces an identity map $d_M:V_2 \to V_1
\subset C^\infty_\J(M,\Lambda^1(M,\C))$.

\punkt The complex vector bundle $S^1(M,\C)$ on $M$ is also
isomorphic to the trivial bundle $\V$. Choose a trivialization
$S^1(M,\C) \cong \V$ in such a way that the canonical map
$C:S^1(M,\C) \to \Lambda^1(M,\C)$ is the identity map. Denote by
$$
V_3 = V \subset C^\infty(M,S^1(M,\C)) \subset \B^0 
$$
the subset of constant sections in $S^1(M,\C) \cong \V$. Then the
Weil algebra $\B^\cdot$ becomes isomorphic to the product
$$
\B^i \cong \wh{S}^\cdot(V_2 \oplus V_3) \otimes \Lambda^i(V_1)
$$
of the completed symmetric algebra $\wh{S}^\cdot(V_2 \oplus V_3)$ of
the sum $V_2 \oplus V_3$ of two copies of the vector space $V$ and
the exterior algebra $\Lambda^\cdot(V_1)$ of the third copy of the
vector space $V$.

\punkt Recall that we have introduced in \ref{aug} a grading on the
Weil algebra $\B^\cdot(M,\C)$ which we call {\em the augmentation
grading}. It induces a grading on the the Weil algebra
$\B^\cdot$. The augmentation grading on $\B^\cdot$ is
multiplicative, and it is obtained by assigning degree $1$ to the
generator subspaces $V_1,V_3 \subset \B^\cdot$ and degree $0$ to the
generator subspace $V_2 \subset \B^\cdot$. As in \ref{aug}, we will
denote the augmentation grading on $\B^\cdot$ by lower indices.

We will now introduce yet another grading on the algebra $\B^\cdot$
which we will call {\em the total grading}. It is by definition the
multiplicative grading obtained by assigning degree $1$ to {\em all}
the generators $V_1,V_2,V_3 \subset \B^\cdot$ of the Weil algebra
$\B^\cdot$ We will denote by $\B^\cdot_{k,n} \subset \B^\cdot$ the
component of augmentation degree $k$ and total degree $n$. Note that
by definition $n,k \geq 0$ and, moreover, $n \geq k$.

\rem In \ref{aug} we have also defined a finer {\em augmentation
bigrading} on the Weil algebra $\B^\cdot(M,\C)$ \, and it this
bigrading that was denoted by double lower indices throughout
Section~\ref{main.section}. We will now longer need the augmentation
bigrading, so there is no danger of confusion.

\punkt The trivialization of the cotangent bundle to $M$ defines an
isomorphism $\tm \cong M \times V$ and a constant Hodge connection
on the pair $\langle \tm,M \rangle$. The corresponding extended
connection $D^{const}:S^1(M,\C) \to \B^1(M,\C)$ is the sum of the
trivial connection
$$
\nabla^{const}_1:S^1(M,\C) \to S^1(M,\C) \otimes
\Lambda^1(M,\C) \subset \B^1(M,\C)
$$ 
on $S^1(M,\C) \cong \V$ and the canonical isomorphism 
$$
C = \id:S^1(M,\C) \to \Lambda^1(M,\C) \subset \B^1(M,\C).
$$ 
The derivation $D^{const}:\B^\cdot \to \B^{\cdot+1}$ of the Weil
algebra associated to the extended connection $D^{const}$ by
\ref{deriv} is equal to
\begin{align*}
D^{const} &= C = \id:V_3 \to V_1\\
D^{const} &= d_M = \id:V_2 \to V_1\\
D^{const} &= d_M = 0 \text{ on }V_1
\end{align*}
on the generator spaces $V_1,V_2,V_3 \subset \B^\cdot$. In
particular, the derivation $D^{const}$ preserves the total degree.

\punkt Let now $D:S^1(M,\C) \to \B^1(M,\C)$ be the an arbitrary
linear extended connection on the manifold $M$, and let
$$
D = \sum_{k \geq 0}D_k:\B^\cdot \to \B^{\cdot+1}
$$
be the derivation of the Weil algebra $\B^\cdot$ associated to the
extended connection $D$ by \ref{deriv}. The derivation $D$ admits a
finer decomposition 
$$
D = \sum_{k,n \geq 0} D_{k,n}:\B^\cdot \to \B^{\cdot+1}
$$
according to both the augmentation and the total degree on
$\B^\cdot$. The summand $D_{k,n}$ by definition raises the
augmentation degree by $k$ and the total degree by $n$. 

\punkt Since the extended connection $D:S^1(M,\C) \to \B^1(M,\C)$ is
linear, its component $D_0:S^1(M,\C) \to \Lambda^1(M,\C)$ of
augmentation degree $0$ coincides with the canonical isomorphism
$C:S^1(M,\C) \to \Lambda^1(M,\C)$. Therefore the restriction of the
derivation $D:\B^\cdot \to \B^{\cdot+1}$ to the generator subspace
$V_3 \subset S^1(M,\C) \subset \B^0$ satisfies
$$
D_0 = C = D^{const}_0 = D^{const}_{0,0}:V_3 \to V_1 \subset \B^1.
$$
In particular, all the components $D_{0,n}$ except for $D_{0,0}$
vanish on the subspace $V_3 \subset \B^\cdot$. 

The restriction of the derivation $D$ to the subspace
$\Lambda^\cdot(M,\C) \subset \B^\cdot(M,\C)$ by definition coincides
with the de Rham differential $d_M:\Lambda^\cdot(M,\C) \to
\Lambda^{\cdot+1}(M,\C)$. Therefore on the generator subspaces
$V_1,V_2 \subset \B^\cdot$ we have $D = d_M = D^{const}$. In
particular, all the components $D_{k,n}$ except for $D_{1,0}$ vanish
on the subspaces $V_1,V_2 \subset \B^\cdot$.

\punkt The fixed Hermitian metric on the generator spaces $V_1 = V_2
= V_3 = V$ extends uniquely to a metric on the whole Weil algebra
such that the multiplication map $\B^\cdot \otimes \B^\cdot \to
\B^\cdot$ is an isometry. We call this metric {\em the standard
metric} on $\B^\cdot$. We finish our preliminary observations with
the following fact which we will use to deduce
Theorem~\ref{converge} from estimates on the components $D_{n,k}$ of
the derivation $D:\B^\cdot \to \B^{\cdot+1}$.

\begin{lemma}\label{est=conv}
Let $D = \sum_{n,k}D_{n,k}:\B^\cdot \to \B^{\cdot+1}$ be a
derivation associated to an extended connection $D$ on the manifold
$M$.  Consider the norms $\|D_{k,n}\|$ of the restrictions
$D_{k,n}:V_3 \to \B^1_{n+1,k+1}$ of the derivations
$D_{k,n}:\B^\cdot \to \B^{\cdot+1}$ to the generator subspace $V_3
\subset \B^0$ taken with respect to the standard metric on the Weil
algebra $\B^\cdot$. If for certain constants $C,\eps > 0$ and for
every natural $n \geq k \geq 0$ we have
\begin{equation}\label{est}
\|D_{k,n}\| < C\eps^n,
\end{equation}
then the formal Hodge connection on $\tm$ along $M$ associated to
$D$ converges to an actual real-analytic Hodge connection on the
open ball of radius $\eps$ in $\tm$ with center at $0 \in M \subset
\tm$. Conversely, if the extended connection $D$ comes from a
real-analytic Hodge connection on an open neighborhood $U \subset
\tm$, and if the Taylor series for this Hodge connection converge in
the closed ball of radius $\eps$ with center at $0 \in M \subset
\tm$, then there exists a constant $C>0$ such that \eqref{est} holds
for every $n,k \geq 0$.
\end{lemma}

\proof The constant Hodge connection $D^{const}$ is obviously defined
on the whole $\tm$, and every other formal Hodge connection on $\tm$
is of the form
$$
D = D^{const} + d^r \circ \Theta:\Lambda^0(\tm,\C) \to
\rho^*\Lambda^1(M,\C), 
$$
where $d^r:\Lambda^0(\tm,\C) \to \Lambda^1(\tm/M,\C)$ is the
relative de Rham differential, and $\Theta \in
C^\infty_M(\tm,\Lambda^1(\tm/M,\C) \otimes \rho^*\Lambda^1(M,\C))$
is a certain relative $1$-form on the formal neighborhood of $M
\subset \tm$ with values in the bundle $\rho^*\Lambda^1(M,\C)$. Both
bundles $\Lambda^1(\tm/M,\C)$ and $\rho^*\Lambda^1(M,\C)$ are
canonically isomorphic to the trivial bundle $\V$ with fiber $V$ on
$\tm$. Therefore we can treat the $1$-form $\Theta$ as a formal germ
of a $\End(V)$-valued function on $\tm$ along $M$. The Hodge
connection $D$ converges on a subset $U \subset \tm$ if and only if
this formal germ comes from a real-analytic $\End(V)$-valued
function $\Theta$ on $U$.

The space of all formal Taylor series for $\End(V)$-valued functions
on $\tm$ at $0 \in M \subset \tm$ is by definition equal to $\End(V)
\otimes \B^0$. Moreover, for every $n \geq 0$ the component
$\Theta_n \in \B^0_n \otimes \End(V) = \Hom(V,\B^0_n \otimes V)$ of
total degree $n$ of the formal power series for the function
$\Theta$ at $0 \in \tm$ is equal to the derivation
$$
\sum_{0 \leq k \leq n}D_{k,n}:V=V_3 \to V=V_1 \otimes \bigoplus_{0 \leq k
\leq n}\B^0_{k,n}.
$$
Every point $x \in \tm$ defines the ``evaluation at $x$'' map
$$
\ev_x:C^\infty(\tm,\End(V)) \to \C, 
$$
and the formal Taylor series for $\Theta \in \End(V) \otimes \B^0$
converges at the point $x \in \tm$ if and only if the series
$$
\Theta(x) = \sum_{n \geq 0}\ev_x(\Theta_n) \in \End(V)
$$
converges. But we have
$$
\|\ev_x(\Theta_n)\| = \left\|\sum_{0 \leq k \leq n}D_{k,n}\right\||x|^n,
$$
where $|x|$ is the distance from the point $x$ to $0 \in \tm$. Now
the application of standard criteria of convergence finishes the
proof of the lemma.  
\endproof

\subsection{Combinatorics}

\punkt We now derive some purely combinatorial facts needed to
obtain estimates for the components $D_{k,n}$ of the extended
connection $D$. First, let $a_n$ be the Catalan numbers, that is,
the numbers defined by the recurrence relation 
$$
a_n = \sum_{1 \leq k \leq n-1} a_k a_{n-k}
$$
and the initial conditions $a_1 = 1$, $a_n = 0$ for $n \leq 0$. As
is well-known, the generating function $f(z) = \sum_{k \geq
0}a_kz^k$ for the Catalan numbers satisfies the equation $f(z) =
f(z)^2 + z$ and equals therefore
$$
f(z) = \half - \sqrt{\frac{1}{4} - z}.
$$
The Taylor series for this function at $z=0$ converges for $4|z|
< 1$, which implies that 
$$
a_k < C(4+\eps)^k
$$
for some positive constant $C > 0$ and every $\eps > 0$. 

\punkt \label{b.n.k}
We will need a more complicated sequence of integers, numbered by
two natural indices, which we denote by $b_{k,n}$. The sequence
$b_{k,n}$ is defined by the recurrence relation
$$
b_{k,n} = \sum_{p,q;1 \leq p \leq k-1} \frac{q+1}{k} b_{p,q}b_{k-p,n-q}
$$
and the initial conditions 
$$
\begin{cases}
b_{k,n} &= 0 \quad \text{ for } \quad k \leq 0,\\
b_{k,n} &= 0 \quad \text{ for } \quad k = 1, n < 0,\\
b_{k,n} &= 1 \quad \text{ for } \quad k = 1, n \geq 0,
\end{cases}
$$
which imply, in particular, that if $n < 0$, then $b_{k,n} = 0$ for
every $k$. For every $k \geq 1$ let $g_k(z) = \sum_{n \geq
0}b_{k,n}z^n$ be the generating function for the numbers
$b_{k,n}$. The recurrence relations on $b_{k,n}$ give 
\begin{align*}
g_k(z) &= \frac{1}{k}\sum_{1 \leq p \leq
k-1}g_{k-p}(z)\left(1+z\frac{\6}{\6 z}\right)(g_p(z)) \\ 
&=\frac{1}{2k}\sum_{1 \leq p \leq k-1}\left(2+z\frac{\6}{\6
z}\right)(g_p(z)g_{k-p}(z)), 
\end{align*}
and the initial conditions give 
$$
g_1(z) = \frac{1}{1-z}.
$$

\punkt Say that a formal series $f(z)$ in the variable $z$ is {\em
non-negative} if all the terms in the series are non-negative real
numbers. The sum and product of two non-negative series and the
derivative of a non-negative series is also obviously
non-negative. For two formal series $s(z)$, $t(z)$ write $s(z) \ll
t(z)$ if the difference $t(z) - s(z)$ is a non-negative power
series.

Our main estimate for the generating functions $g_k(z)$ is the
following. 

\begin{lemma}
For every $k \geq 1$ we have 
$$
g_k(z) \ll a_k\frac{1}{(1-z)^{2k-1}}, 
$$
where $a_k$ are the Catalan numbers. 
\end{lemma}

\proof Use induction on $k$. For $k=1$ we have $g_1(z) =
\frac{1}{1-z}$ and $a_1 = 1$, which gives the base for
induction. Assume that the claim is proved for all $p < k$. Since
all the $g_n(z)$ are non-negative power series, this implies that for
every $p$, $1 \leq p \leq k-1$ we have
$$
g_p(z)g_{k-p}(z) \ll
a_pa_{k-p}\frac{1}{(1-z)^{2p-1}}\frac{1}{(1-z)^{2k-2p-1}} =
a_pa_{k-p}\frac{1}{(1-z)^{2k-2}}. 
$$
Therefore 
\begin{align*}
\left(2+z\frac{\6}{\6 z}\right)(g_p(z)g_{k-p}(z)) &\ll
a_pa_{k-p}\left(2+z\frac{\6}{\6 z}\right)\frac{1}{(1-z)^{2k-2}} \\
&= a_pa_{k-p}\left(\frac{2}{(1-z)^{2k-2}} +
\frac{(2k-2)z}{(1-z)^{2k-1}} \right) \\  
&= a_pa_{k-p}\left(\frac{2k-2}{(1-z)^{2k-1}} -
\frac{2k-4}{(1-z)^{2k-2}}\right) \\ 
&\ll (2k-2)a_pa_{k-p}\frac{1}{(1-z)^{2k-1}}.
\end{align*}
Hence
\begin{align*}
g_k(z) &= \frac{1}{2k}\sum_{1 \leq p \leq k-1}\left(2+z\frac{\6}{\6
z}\right)(g_p(z)g_{k-p}(z)) \\ 
&\ll \frac{2k-2}{2k}\sum_{1 \leq p \leq
k-1}a_pa_{k-p}\frac{1}{(1-z)^{2k-1}} \\ 
&\ll \frac{1}{(1-z)^{2k-1}}\sum_{1 \leq p \leq k-1}a_pa_{k-p} =
a_k\frac{1}{(1-z)^{2k-1}}, 
\end{align*}
which proves the lemma. 
\endproof 

\punkt This estimate yields the following estimate for the numbers
$b_{k,n}$. 

\begin{corr}\label{combin}
The power series 
$$
g(z) = \sum_{k,n} b_{k,n}z^{n+k} = \sum_{k \geq 1}g_k(z)z^k
$$
converges for $z < 3-\sqrt{8}$. Consequently, for every $C_2$ such
that $(3-\sqrt{8})C_2 > 1$ there exists a positive constant $C > 0$
such that
$$
b_{n,k} < CC_2^{n+k}
$$
for every $n$ and $k$. (One can take, for example, $C_2 = 6$.) 
\end{corr}

\proof Indeed, we have 
\begin{equation}\label{g(z)}
g(z) \ll \sum_{k \geq 1}a_kz^k\frac{1}{(1-z)^{2k-1}} =
(1-z)f\left(\frac{z}{(1-z)^2}\right), 
\end{equation}
where $f(z) = \frac{1}{2} - \sqrt{\frac{1}{4}-z}$ is the generating
function for the Catalan numbers. Therefore
$$
g(z) \ll (1-z)\left(\frac{1}{2} - \sqrt{\frac{1}{4} -
\frac{z}{(1-z)^2}}\right), 
$$
and the right hand side converges absolutely when 
$$
\frac{|z|}{(1-z)^2} < \frac{1}{4}. 
$$
Since $3 - \sqrt{8}$ is the root of the quadratic equation $(1-z)^2
= 4z$, this inequality holds for every $z$ such that $|z| < 3 -
\sqrt{8}$.  
\endproof

\punkt \label{b.n.k.m}
To study polarizations of Hodge manifold structures on $\tm$, we
will need yet another recursive sequence of integers, which we
denote by $b_{k,n}^m$. This sequence is defined by the recurrence
relation
$$
b^m_{k,n} = \sum_{p,q;1 \leq p \leq k-1} \frac{q+m(k-p)}{k}
b^m_{p,q}b_{k-p,n-q}
$$
and the initial conditions 
$$
\begin{cases}
b^m_{k,n} &= 0 \quad \text{ unless } \quad k,m \leq 0,\\
b^m_{k,n} &= 0 \quad \text{ for } \quad k = 1, n < 0,\\
b^m_{k,n} &= 1 \quad \text{ for } \quad k = 1, n \geq 0. 
\end{cases}
$$

\punkt \label{c.k.n}
To estimate the numbers $b^m_{k,n}$, consider the auxiliary sequence
$c_{k,n}$ defined by setting
$$
c_{k,n} = \sum_{p,q;1 \leq p \leq k-1}c_{p.q}b_{k-p,n-q} \qquad k
\geq 2, 
$$
and $c_{k,n} = b_{k,n}$ for $k \leq 1$. The generating series $c(z)
= \sum_{k,n \geq 0}c_{k,n}z^{n+k}$ satisfies
$$
c(z) = c(z)g(z) + \frac{z}{1-z},
$$
so that we have $c(z) = \frac{z}{(1-z)(1-g(z))}$, which is
non-singular when $|z|, |g(z)| < 1$ and $g(z)$ is non-singular. 
By \eqref{g(z)} the latter inequality holds if 
$$
\left|(1-z)\left(\frac{1}{2} - \sqrt{\frac{1}{4} -
\frac{z}{(1-z)^2}}\right)\right| < 1,
$$
which holds in the whole disc where $g(z)$ converges, that is, for
$|z| < 3-\sqrt{8}$. Therefore, as in Corollary~\ref{combin}, we have
\begin{equation}\label{cnk}
c_{k,n} < C6^{n+k}
\end{equation}
for some positive constant $C$. 

\punkt We can now estimate the numbers $b^m_{k,n}$. 

\begin{lemma}\label{combin.2}
For every $m,k,n$ we have 
\begin{equation}\label{bnkm.indu}
b^m_{k,n} \leq (2m)^{k-1}c_{k,n}b_{k,n},
\end{equation}
where $c_{k,n}$ are as in \ref{c.k.n} and $b_{k,n}$ are the
numbers introduced in \ref{b.n.k}. Consequently, we have 
$$
b^m_{k,n} < C(72m)^{n+k+m}
$$
for some positive constant $C > 0$. 
\end{lemma}

\proof Use induction on $k$. The case $k=1$ follows from the initial
conditions. Assume the estimate \eqref{bnkm.indu} proved for all
$b^m_{p,q}$ with $p < k$. Note that by the recurrence relations we
have $b_{p,q} \leq b_{k,n}$ and $c_{p,q} \leq c_{k,n}$ whenever $p <
k$. Therefore
\begin{align*}
b^m_{k,n} &= \sum_{p,q;1 \leq p \leq k-1} \frac{q+m(k-p)}{k}
b^m_{p,q}b_{k-p,n-q} \\
&\leq \sum_{p,q;1 \leq p \leq k-1} \frac{q}{k}b^m_{p,q}b_{k-p,n-q} +
\sum_{p,q;1 \leq p \leq k-1} b^m_{p,q}b_{k-p,n-q} \\
&\leq \sum_{p,q;1 \leq p \leq k-1}2^{p-1} c_{p,q} \frac{q}{k}
b_{p,q}b_{k-p,n-q} + \sum_{p,q;1 \leq p \leq k-1} 2^{p-1} b_{p,q}
c_{p,q} b_{k-p,n-q} \\
&\leq 2^{k-2} c_{k,n}\sum_{p,q;1 \leq p \leq k-1} \frac{q+1}{k}
b_{p,q} b_{k-p,n-q} \\
&\quad + 2^{k-2}b_{k,n}\sum_{p,q;1 \leq p \leq k-1} 
c_{p,q} b_{k-p,n-q} \\
&= 2^{k-2}c_{k,n}b_{k,n} + 2^{k-2}c_{k,n}b_{k,n} = 2^{k-1}c_{k,n}b_{k,n},
\end{align*}
which proves \eqref{bnkm.indu} for $b^m_{k,n}$. The second estimate of
the lemma now follows from \eqref{cnk} and Corollary~\ref{combin}. 
\endproof 

\subsection{The main estimate}

\punkt \label{norms}
Let now $D = \sum_{k,n}D_{n,k}:\B^\cdot \to \B^{\cdot+1}$ be a
derivation of the Weil algebra $\B^\cdot$ associated to a flat
linear extended connection on $M$.  Consider the restriction
$D_{k,n}:\B^\cdot_{p,q} \to \B^\cdot{p+k,q+n}$ of the derivation
$D_{k,n}$ to the component $\B^\cdot_{p,q}\subset\B^\cdot$ of
augmentation degree $p$ and total degree $q$. Since both
$\B^\cdot_{p,q}$ and $\B^\cdot{p+k,q+n}$ are finite-dimensional
vector spaces, the norm of this restriction with respect to the
standard metric on $\B^\cdot$ is well-defined. Denote this norm by
$\|D_{k,n}\|_{p,q}$.

By Lemma~\ref{est=conv} the convergence of the Hodge manifold
structure on $\tm$ corresponding to the extended connection $D$ is
related to the growth of the norms $\|D_{k,n}\|_{1,1}$.  Our main
estimate on the norms $\|D_{k,n}\|_{1,1}$ is the following.

\begin{prop}\label{estimate}
Assume that there exist a positive constant $C_0$ such that for
every $n$ the norms $\|D_{1,n}\|_{1,1}$ and $\|D_{1,n}\|_{0,1}$ 
satisfy 
$$
\|D_{1,n}\|_{1,1}, \|D_{1,n}\|_{1,0} < C_0^n.
$$
Then there exists a positive constant $C_1$ such that for every
$n,k$ the norm $\|D_{k,n}\|_{1,1}$ satisfies
$$
\|D_{k,n}\|_{1,1} < C_1^n.
$$
\end{prop}

\punkt \label{metr}
In order to prove Proposition~\ref{estimate}, we need some
preliminary facts. Recall that we have introduced in
\ref{total.Weil} the total Weil algebra $\B_{tot}^\cdot(M,\C)$ of the
manifold $M$, and let
$$
\B^\cdot_{tot} = C^\infty_\J(M,\B_{tot}^\cdot(M,\C))
$$
be the algebra of its smooth sections completed at $0 \subset M$. By
definition for every $k \geq 0$ we have $\B_{tot}^k = \B^k \otimes
\W_k^*$, where $\W_k$ is the $\R$-Hodge structure of weight $k$
universal for weakly Hodge maps, as in \ref{w.k.uni}. There exists a
unique Hermitian metric on $\W_k$ such that all the Hodge components
$\W^{p,q} \subset \W_k$ are orthogonal and all the Hodge degree
components $w_k^{p,q}$ of the universal weakly Hodge map $w_k:\R(0)
\to \W_k$ are isometries. This metric defines a canonical Hermitian
metric on $\W_k^*$.

\defn The {\em standard metric} on the total Weil algebra
$\B_{tot}^\cdot$ is the product of the canonical metric and the standard
metric on $\B^\cdot$.

\punkt By Lemma~\ref{total.rel} the total Weil algebra $\B^\cdot_{tot}$
is generated by the subspaces $V_2,V_3 \subset \B^0 = \B^0_{tot}$ and
the subspace $V_1 \otimes \W_1^* \subset \B^1_{tot}$, which we denote by
$V_1^{tot}$. The ideal of relations for the algebra $\B^\cdot_{tot}$ is the
ideal in $S^\cdot(V_2 \oplus V_3)\otimes\Lambda^\cdot(V_1^{tot})$
generated by $S^2(V_1)\otimes\Lambda^2(\W_1^*) \subset
\Lambda^2(V_1^{tot})$. 

The direct sum decomposition \eqref{drct} induces a direct sum
decomposition 
$$
V^{tot}_1 = V^{ll}_1 \oplus \V^o_1 \oplus V^{rr}_1 
$$
of the generator subspace $V^{tot}_1 \subset \B^1_{tot}$. The subspaces
$V^o_1 \subset V^{tot}_1$ and $V^{ll}_1 \oplus V^{rr}_1 \subset V^{tot}_1$
are both isomorphic to the vector space $V_1$. More precisely, the
universal weakly Hodge map $w_1:\R(0) \to \W_1$ defines a projection
$$
P:V^{tot}_1 = V_1 \otimes \W_1^* \to V_1,
$$
and the restriction of the projection $P$ to either of the subspaces
$V^o_1,V^{ll}_1 \oplus V^{rr}_1 \subset V^{tot}_1$ is an
isomorphism. Moreover, either of these restrictions is an isometry
with respect to the standard metrics. 

\punkt \label{P2}
The multiplication in $\B_{tot}^\cdot$ is not an isometry with respect
to this metric. However, for every $b_1,b_2 \subset \B_{tot}^\cdot$ we
have the inequality
$$
\|b_1b_2\| \leq \|b_1\|\cdot\|b_2\|. 
$$
Moreover, this inequality becomes an equality when $b_1 \subset
\B^0_{tot}$. In particular, if we extend the map $P:V^{tot}_1 \to V_1$ to a
$\B^0$-module map 
$$
P:\B^1_{tot} \to \B^1,
$$
then the restriction of the map $P$ to either of the subspaces
$\B^1_o, \B^1_{ll}\oplus\B^1_{rr} \subset \B^1_{tot}$ is an isometry
with respect to the standard metric. Therefore the norm of the
projection $P:\B^1_{tot} \to \B^1$ is at most $2$. 

\punkt The total and augmentation gradings on the Weil algebra
$\B^\cdot$ extend to gradings on the total Weil algebra
$\B_{tot}^\cdot$, also denoted by lower indices. The extended connection
$D$ on $M$ induces a derivation $D^{tot} =
\sum_{n,k}D_{k,n}^{tot}:\B^\cdot_{tot} \to \B^{\cdot+1}_{tot}$. As in
\ref{norms}, denote by $\|D^{tot}_{k,n}\|_{p,q}$ the norm of the map
$D^{tot}_{k,n}:\left(\B^\cdot_{tot}\right)_{p,q} \to
\left(\B^\cdot_{tot}\right)_{p+k,q+n}$ with respect to the standard
metric on $\B^\cdot_{tot}$. The derivations $D^{tot}_{k,n}$ are related to
$D_{k,n}$ by
$$
D_{k,n} = P \circ D^{tot}_{k,n}:\B^0 \to \B^1, 
$$
and we have the following.

\begin{lemma}\label{total.ne.total}
For every $k$, $n$ and $p=0,1$ we have 
$$
\|D_{k,n}\|_{p,1} = \|D^{tot}_{k,n}\|_{p,1}.
$$
\end{lemma}

\proof By definition we have $\B^\cdot_{0,1} \oplus \B^\cdot_{1,1} =
V_1 \oplus V_2 \oplus V_3$. Moreover, the derivation $D_{k,n}$
vanishes on $V_1$, hence $D^{tot}_{k,n}$ vanishes on $V_1^{tot}$. Therefore
it suffices to compare their norms on $V_2 \oplus V_3 \subset \B^0 =
\B^0_{tot}$.

Since on $S^\cdot(V_2) \subset \B^0$ the derivation $D_{k,n}$
coincides with the de Rham differential, the derivation $D^{tot}_{k,n}$
maps the subspace $V_2$ into $\B^1_{ll}\oplus\B^1_{rr}$. Moreover,
by Lemma~\ref{total.aug} the derivation $D_{k,n}^{tot}$ maps $V_3$
either into $\B_o^1$ or into $\B^1_{ll} \oplus \B^1_{rr}$, depending
on the parity of the number $k$. Since $D_{k,n} = P \circ D^{tot}_{k,n}$
and the map $P:\B^1_{tot} \to \B^1$ is an isometry on both $\B^1_o
\subset \B^1_{tot}$ and $\B^1_{ll}\oplus\B^1_{rr} \subset \B^1_{tot}$, we
have $\|D_{k,n}\| = \|D^{tot}_{k,n}\|$ on both $V_2 \subset \B^0$ and
$V_3 \subset \B^0$, which proves the lemma.
\endproof

This lemma allows to replace the derivations $D_{k,n}$ in
Proposition~\ref{estimate} with associated derivations $D_{k,n}^{tot}$ of
the total Weil algebra $\B_{tot}^\cdot$.

\punkt Since the extended connection $D$ is linear and flat, the
construction used in the proof of Lemma~\ref{main.ind} shows that 
\begin{equation}\label{indu}
D^{tot}_k = h^{-1} \circ \sigma_{tot} \circ \sum_{1 \leq p \leq k-1}D_p^{tot}
\circ D^{tot}_{k-p}:V_3 \to \left(\B^1_{tot}\right)_{k+1},
\end{equation}
where $\sigma_{tot}:\B^{\cdot+1}_{tot} \to \B^\cdot_{tot}$ is the canonical map
constructed in \ref{sigma.c}, and $h:\B^\cdot_{tot} \to \B^\cdot_{tot}$ is
as in Lemma~\ref{h.acts}. Both $\sigma_{tot}$ and $h$ preserve the
augmentation degree. In order to obtain estimates on
$\|D_{k,n}\|_{1,1} = \|D_{k,n}^{tot}\|_{1,1}$, we have to estimate the
norms $\|h^{-1}\|$ and $\|\sigma_{tot}\|$ of the restrictions of maps
$h^{-1}$ and $\sigma_{tot}$ on the subspace
$\left(\B_{tot}^\cdot\right)_{k+1} \subset \B^\cdot_{tot}$.

By Lemma~\ref{h.on.b1} the map $h:\left(\B^\cdot\right)_{k+1} \to
\left(\B^\cdot\right)_{k+1}$ is diagonalizable, with eigenvalues
$k+1$ if $k$ is even and $(k+1)/2,(k-1)/2$ if $k$ is odd. Since $m
\geq 2$, in any case on $\left(\B^\cdot\right)_{k+1}$ we have
\begin{equation}\label{h.est}
\|h^{-1}\| < \frac{3}{k}.
\end{equation}

\punkt \label{sigma.est.punkt}
To estimate $\sigma_{tot}:\B^2_{tot} \to \B^2_{tot}$, recall that, as noted in
\ref{sigma.c}, the map $\sigma_{tot}$ is a map of $\B^0$-modules. The
space $\B^2_{tot} = \B^0 \otimes \left(\B^2_{tot}\right)_2$ is a free
$\B^0$-module generated by a finite-dimensional vector space
$\left(\B^2_{tot}\right)_2$. The map $\sigma_{tot}$ preserves the
augmentation degree, hence it maps $\left(\B^2_{tot}\right)_2$ into the
finite-dimensional vector space $\left(\B^1_{tot}\right)_2$. Therefore
there exists a constant $K$ such that
\begin{equation}\label{sigma.est}
\|\sigma_{tot}\| \leq K
\end{equation}
on $\left(\B^2_{tot}\right)_2$. Since we have
$\|b_1b_2\|=\|b_1\|\cdot\|b_2\|$ for every $b_1 \in \B^0$, $b_2 \in
\B^2_{tot}$, and the map $\sigma_{tot}$ is $\B^0$-linear, the estimate
\eqref{sigma.est} holds on the whole $\B^2_{tot} = \B^0 \otimes
\left(\B^2_{tot}\right)_2$. We can assume, in addition, that $K \geq 1$. 

\rem In fact $K = 2$, but we will not need this. 

\punkt Let now $b_{k,n}$ be the numbers defined recursively in
\ref{b.n.k}. Our estimate for $\|D^{tot}_{k,n}\|_{p,q}$ is the
following. 

\begin{lemma}\label{est.indu}
In the assumptions and notations of Proposition~\ref{estimate}, we
have 
$$
\|D^{tot}_{k,n}\|_{p,q} < q(3K)^{k-1} C_0^nb_{k,n}
$$
for every $k,n$. 
\end{lemma}

\proof Use induction on $k$. The base of induction is the case
$k=1$, when the inequality holds by assumption. Assume that for some
$k$ we have proved the inequality for all 
$\|D^{tot}_{m,n}\|_{p,q}$ with $m < k$, and fix a number $n \geq k$.  

Consider first the restriction of $D^{tot}_{k,n}$ onto the generator
subspace $V_3 \subset \B^0$.  Taking into account the total degree,
we can rewrite \eqref{indu} as
$$
D^{tot}_{k,n} = h^{-1} \circ \sigma_{tot} \circ \sum_{1 \leq p \leq k-1}\sum_q
D^{tot}_{k-p,n-q} \circ D^{tot}_{p,q}:V_3 \to \left(\B^1_{tot}\right)_{k+1}. 
$$
Therefore the norm of the map $D^{tot}_{k,n}:V_3 \to \B_{tot}^1$ satisfies 
$$
\|D^{tot}_{k,n}\} \leq \|h^{-1}\| \cdot \|\sigma_{tot}\| \cdot \sum_{1 \leq p
\leq k-1}\sum_q\|D^{tot}_{k-p,n-q}\|_{p+1,q+1} \cdot \|D^{tot}_{p,q}\|_{1,1}. 
$$
Substituting into this the estimates \eqref{h.est},
\eqref{sigma.est} and the inductive assumption, we get 
$$
\|D^{tot}_{k,n}\| < \frac{3}{k}K\sum_{1 \leq p \leq
k-1}\sum_q(q+1)(3K)^{k-2}C_0^nb_{k-p,n-q}b_{p,q} = (3K)^{k-1}C_0^nb_{k,n}.
$$
Since by definition $D^{tot}_{k,n}$ vanishes on $V_2 \subset \B^0$ and
on $V^{tot}_1 \subset \B^1_{tot}$, this proves that 
$$
\|D^{tot}_{k,n}\|_{p,1} < (3K)^{k-1}C_0^nb_{k,n}
$$
when $p=0,1$. Since the map $D^{tot}_{k,n}:\B^\cdot_{tot} \to
\B^{\cdot+1}_{tot}$ is a derivation, the Leibnitz rule and the triangle
inequality show that for every $p$, $q$
$$
\|D^{tot}_{k,n}\|_{p,q} < q(3K)^{k-1}C_0^nb_{k,n}, 
$$
which proves the lemma.
\endproof 

\punkt \proof[Proof of Proposition~\ref{estimate}]
By Lemma~\ref{total.ne.total} we have $\|D_{k,n}\|_{p,1} =
\|D^{tot}_{k,n}\|_{p,1}$, and by Lemma~\ref{est.indu} we have 
$$
\|D_{k,n}\|_{p,1} = \|D^{tot}_{k,n}\|_{p,1} < (3K)^{k-1}C_0^nb_{n,k}.
$$
Since $k \leq n$ and $K \geq 1$, this estimate together with
Corollary~\ref{combin} implies that
$$
\|D_{k,n}\|_{p,1} < C(3K)^{k-1}C_0^n6^{2n} < C(108KC_0)^n 
$$
for some positive constant $C > 0$, which proves the proposition. 
\endproof 

\punkt Proposition~\ref{estimate} gives estimates for the derivation
$D:\B^\cdot \to \B^{\cdot+1}$ or, equivalently, for the Dolbeault
differential 
$$
\bar\6_J:\Lambda^{0,\cdot}(\tm_J) \to \Lambda^{0,\cdot+1}(\tm_J)
$$ 
for the complementary complex structure $\tm_J$ on $\tm$ associated
to the extended connection $D$ on $M$. To prove the second part of
Theorem~\ref{converge}, we will need to obtain estimates on the
Dolbeault differential $\bar\6_J:\Lambda^{p,0}(\tm_J) \to
\Lambda^{p,1}(\tm_J)$ with $p > 0$. To do this, we use the model for
the de Rham complex $\Lambda^{\cdot,\cdot}(\tm_J)$ constructed in
Subsection~\ref{drm.mod}.

Recall that in \ref{ident.punkt} we have identified the direct image
$\rho_*\Lambda^{\cdot,\cdot}(\tm_J)$ of the de Rham algebra of the
manifold $\tm$ with the free module over the Weil algebra
$\B^\cdot(M,\C)$ generated by a graded algebra bundle
$L^\cdot(M,\C)$ on $M$. The Dolbeault differential $\bar\6_J$ for
the complementary complex structure $\tm_J$ induces an algebra
derivation $D:L^\cdot(M,\C) \otimes \B^\cdot(M,\C) \to
L^\cdot(M,\C) \otimes \B^{\cdot+1}(M,\C)$, so that the free module
$\rho_*\Lambda^{\cdot,\cdot}(\tm_J) \cong L^\cdot(M,\C) \otimes
\B^\cdot(M,\C)$ becomes a differential graded module over the Weil
algebra. 

\punkt The algebra bundle $L^\cdot(M,\C)$ is isomorphic
to the de Rham algebra $\Lambda^\cdot(M,\C)$. In particular, the
bundle $L^1(M,\C)$ is isomorphic to the trivial bundle $\V$ with
fiber $V$ over $M$.  By \ref{dr.L} the relative de Rham differential
$$
d^r:\Lambda^\cdot(\tm/M,\C) \to \Lambda^{\cdot+1}(\tm/M,\C)
$$
induces a derivation
$$
d^r:L^\cdot(M,\C) \otimes \B^\cdot(M,\C) \to
L^{\cdot+1}(M,\C) \otimes \B^\cdot(M,\C), 
$$
and we can choose the trivialization $L^1(M,\C) \cong \V$ in such a
way that $d^r$ identifies the generator subspace $V_3 \subset \B^0$
with the subspace of constant sections of $\V \cong L^1(M,\C)
\subset L^1(M,\C) \otimes \B^0(M,\C)$. 

\punkt \label{LLL}
Denote by
$$
\LL^{\cdot,\cdot} = C^\infty_\J(L^\cdot(M,\C) \otimes
\B_{tot}^\cdot(M,\C))
$$
the $\J$-adic completion of the space of smooth sections of the
algebra bundle $L^\cdot(M,\C) \otimes \B_{tot}^\cdot(M,\C)$ on $M$.  The
space $\LL^{\cdot,\cdot}$ is a bigraded algebra equipped with the
derivations $d^r:\LL^{\cdot,\cdot} \to \LL^{\cdot+1,\cdot}$,
$D^{tot}:\LL^{\cdot,\cdot} \to \LL^{\cdot,\cdot+1}$, which commute by
Lemma~\ref{D.dr}. The algebra $\LL^{\cdot,\cdot}$ is the free
graded-commutative algebra generated by the subspaces $V_1,V_2,V_3
\subset \LL^{0,\cdot} = \B^\cdot$ and the subspace $V = d^r(V_3)
\subset \LL^{1,0}$, which we denote by $V_4$.

\punkt As in \ref{metr}, the given metric on the generator subspaces
$V_1=V_2=V_3=V_4=V$ extends uniquely to a multiplicative metric on
the algebra $\LL^{\cdot,\cdot}$, which we call {\em the standard
metric}. For every $k > 0$, introduce the total and augmentation
gradings on the free $\B^\cdot$-module $\LL^{k,\cdot} =
\Lambda^k(V_4) \otimes \B^\cdot$ by setting $\deg \Lambda^k(V_4) =
(0,0)$.  Let $D = D_{k,n}$ be the decomposition of the derivation
$D:\LL^{\cdot,\cdot} \to \LL^{\cdot,\cdot+1}$ with respect to the
total and the augmentation degrees. Denote by $\|D_{k,n}\|^p_q$ the
norm with respect to the standard metric of the restriction of the
derivation $D_{k,n}:\LL^{p,0} \to \LL^{p,1}$ to the component in
$\LL^{p,0}$ of total degree $q$.

\punkt Let now $\LL_{tot}^{\cdot,\cdot} = \Lambda^\cdot(V_4) \otimes
\B^\cdot_{tot}$ be the product of the exterior algebra
$\Lambda^\cdot(V_4)$ with the total Weil algebra $\B^\cdot_{tot}$. We
have the canonical identification $\LL^{p,q}_{tot} = \LL^{p,q} \otimes
\W_q^*$, and the canonical projection $P:\LL^{\cdot,\cdot}_{tot} \to
\LL^{\cdot,\cdot}$, identical on $\LL^{\cdot,0}_{tot} = \LL^{\cdot,0}$. 
As in \ref{P2}, the norm of the projection $P$ on $\LL^{\cdot,1}_{tot}$
is at most $2$. 

The derivation $D:\LL^{\cdot,0} \to \LL^{\cdot,1}$ induces a
derivation $D^{tot}:\LL_{tot}^{\cdot,0} \to \LL_{tot}^{\cdot,1}$, related to $D$
by $D = P \circ D^{tot}$. The gradings and the metric on
$\LL^{\cdot,\cdot}$ extend to $\LL^{\cdot,\cdot}_{tot}$, in particular,
we have the decomposition $D^{tot} = \sum_{k,n}D^{tot}_{k,n}$ with respect
to the total and the augmentation degrees.  Denote by
$\|D_{k,n}^{tot}\|^p_q$ the norm with respect to the standard metric of
the restriction of the derivation $D^{tot}_{k,n}:\LL_{tot}^{p,0} \to
\LL_C^{p,1}$ to the component in $\LL^{p,0}$ of total degree
$q$. Since $\|P\|\leq 2$ on $\LL^{\cdot,1}_{tot}$, we have
\begin{equation}\label{DcL}
\|D_{k,n}^{tot}\|^p_q \leq 2\cdot\|D_{k,n}\|^p_q.
\end{equation}

\punkt The estimate on the norms $\|D^{tot}_{k,n}\|^\LL_{p,q}$ that we
will need is the following.

\begin{lemma}\label{D2est}
In the notation of Lemma~\ref{est.indu}, we have 
$$
\|D^{tot}_{k,n}\|^p_q < 2(q+pk)(3K)^{k-1} C_0^nb_{k,n},
$$
for every $p,q,k,n$. 
\end{lemma}

\proof Since $D_{k,n}$ satisfies the Leibnitz rule, it suffices to
prove the estimate for the restriction of the derivation $D^{tot}_{k,n}$
to the generator subspaces $V_2,V_3,V_4 \subset \LL^{\cdot,0}$. By
\eqref{DcL} it suffices to prove that on $V_2,V_3,V_4$ we have
$$
\|D_{k,n}\|^p_q < (q+pk)(3K)^{k-1} C_0^nb_{k,n}. 
$$
On the generator subspaces $V_2,V_3 \subset \B^0$ we have $p=0$, and
this equality is the claim of Lemma~\ref{est.indu}. Therefore it
suffices to consider the restriction of the derivation $D_{k,n}$ to
the subspace $V_4 = d^r(V_3) \subset \LL^{1,0}$. We have $d^r \circ
D_{k,n} = D_{k,n} \circ d^r:V_3 \to \LL^{1,1}$. Moreover, $d^r =
\id:V_3 \to V_4$ is an isometry. Since the operator $d^r:\B^\cdot
\to \LL^{1,\cdot}$ satisfies the Leibnitz rule and vanishes on the
generators $V_1,V_2 \subset \B^\cdot$, the norm $\|d^r\|_k$ of its
restriction to the subspace $\B^\cdot_k$ of augmentation degree $k$
does not exceed $k$. Therefore
\begin{align*}
\|D_{k,n}\|^1_0 &= \left\|D_{k,n}|_{V_4}\right\| =
\left\|D_{k,n}\circ d^r|_{V_3}\right\| = \left\|d^r \circ
D_{k,n}|_{V_3}\right\| \\
&\leq \left\|D_{k,n}|_{V_3}\right|\cdot\left\|d^r_{\B^1_k}\right\| < k
(3K)^{k-1} C_0^nb_{k,n},
\end{align*}
which proves the lemma.
\endproof

\subsection{The proof of Theorem~\ref{converge}} 

\punkt We can now prove Theorem~\ref{converge}. Let $D:S^1(M,\C) \to
\B^1(M,\C)$ be a flat linear extended connection on $M$. Assume that
its reduction $D_1=\nabla:S^1(M,\C) \to S^1(M,\C) \otimes
\Lambda^1(M,\C)$ is a real-analytic connection on the bundle
$S^1(M,\C)$.

The operator $D_1:S^1(M,\C) \to S^1(M,\C) \otimes \Lambda^1(M,\C)
\subset \B^1(M,\C)$ considered as an extended connection on $M$
defines a Hodge connection $D_1:\Lambda^0(\tm,\C) \to
\rho^*\Lambda^1(M,\C)$ on the pair $\langle \tm,M \rangle$, and this
Hodge connection is also real-analytic. Assume further that the
Taylor series at $0 \subset M \subset \tm$ for the Hodge connection
$D_1$ converge in the closed ball of radius $\eps>0$.

\punkt Let $D = \sum_{n,k}D_{k,n}:\B^\cdot \to \B^{\cdot+1}$ be the
derivation of the Weil algebra $\B^\cdot$ associated to the extended
connection $D$.  Applying Lemma~\ref{est=conv} to the Hodge
connection $D_1$ proves that there exists a constant $C>0$ such that
for every $n \geq 0$ the norm $\|D_{1,n}\|_{1,1}$ of the restriction
of the derivation $D_{1,n}$ to the generator subspace $V_3 =
\B^0_{1,1} \subset \B^0$ satisfies 
$$
\|D_{1,n}\|_{1,1} < CC_0^n,
$$
where $C_0 = 1/\eps$. 

By definition the derivation $D_{1,n}$ vanishes on the generator
subspace $V_1 \subset \B^1$. If $n > 0$, then it also vanishes on
the generator subspace $V_2 = \B^0_{0,1} \subset \B^0$. If $n=0$,
then its restriction to $V_2 = \B^0_{0,1} \subset \B^0$ is the
identity isomorphism $D_{1,0}=\id:V_2 \to V_1$. In any case, we have
$\|D_{1,n}\|_{0,1} \leq 1$. Increasing if necessary the constant
$C_0$, we can assume that for any $n$ and for $p=0,1$ we have
$$
\|D_{1,n}\|_{p,1} < C_0^n.
$$

\punkt We can now apply our main estimate,
Proposition~\ref{estimate}. It shows that there exists a constant
$C_1 > 0$ such that for every $k,n$ we have
$$
\|D_{k,n}\|_{1,1} < C_1^n.
$$
Together with Lemma~\ref{est=conv} this estimate implies that the
formal Hodge connection $D$ on $\tm$ along $M \subset \tm$
corresponding to the extended connection $D$ converges to a
real-analytic Hodge connection on an open neighborhood $U \subset
\tm$ of the zero section $M \subset \tm$. This in turn implies the
first claim of Theorem~\ref{converge}.

\punkt To prove the second claim of Theorem~\ref{converge}, assume
that the manifold $M$ is equipped with a K\"ahler form $\omega$
compatible with the K\"ah\-le\-ri\-an connection $\nabla$, so that
$\nabla\omega=0$. The differential operator $\nabla:\Lambda^{1,1}(M)
\to \Lambda^{1,1}(M) \otimes \Lambda^1(M,\C)$ is elliptic and
real-analytic. Since $\nabla\omega=0$, the K\"ahler form $\omega$ is
also real-analytic.

For every $p,q \geq 0$ we have introduced in \ref{LLL} the space
$\LL^{p,q}$, which coincides with the space of formal germs at $0
\in M \subset \tm$ of smooth forms on $\tm$ of type $(p,q)$ with
respect to the complementary complex structure $\tm_J$. The spaces
$\LL^{\cdot,\cdot}$ carry the total and the augmentation
gradings. Consider $\omega$ as an element of the vector space
$\LL^{2,0}_0$, and let $\omega = \sum_n\omega_n$ be the total degree
decomposition.  The decomposition $\omega = \sum_n \omega_n$ is the
Taylor series decomposition for the form $\omega$ at $0 \in M$. Since
the form $\omega$ is real-analytic, there exists a constant $C_2$
such that
\begin{equation}\label{omega.est}
\|\omega_n\| < C_2^n
\end{equation}
for every $n$. 

\punkt Let $\Omega = \sum_k \Omega_k = \sum_{k,n} \Omega_{k,n}
\subset \LL^{2,0}$ be the formal polarization of the Hodge manifold
$\tm$ at $M \subset \tm$ corresponding to the K\"ahler form $\omega$
by Theorem~\ref{metrics}. By definition we have $\omega =
\Omega_0$. Moreover, by construction used in the proof of
Proposition~\ref{metrics.ind} we have 
\begin{equation}\label{Omega.rec}
\Omega_k = -\frac{1}{k}\sum_{1 \leq p \leq
k-1}\sigma_{tot}(D^{tot}_{k-p}\Omega_p), 
\end{equation}
where $D^{tot}:\LL^{2,0} \to \LL^{2,1}$ is the derivation associated to
the extended connection $D$ on $M$ and $\sigma_{tot}:\LL^{2,\cdot+1} \to
\LL^{2,\cdot}$ is the extension to $\LL^{2,\cdot} = L^2_0 \otimes
\B^\cdot_{tot}$ of the canonical endomorphism of the total Weil algebra
$\B^\cdot_{tot}$ constructed in the proof of
Proposition~\ref{metrics.ind}.

\punkt The map $\sigma_{tot}:\LL^{2,1} \to \LL^{2,0}$ is a map of
$\B^0_{tot}$-modules. Therefore, as in \ref{sigma.est.punkt}, there
exists a constant $K_1 > 0$ such that 
\begin{equation}\label{K1}
\|\sigma_{tot}\| < K_1
\end{equation}
on $\LL^{2,1}$. We can assume that $K_1 > 3K$, where $K$ is as in
\eqref{sigma.est}.  Together with the recursive formula
\eqref{Omega.rec}, this estimate implies the following estimate on
the norms $\|\Omega_{k,n}\|$ of the components $\Omega_{k,n}$ of the
formal polarization $\Omega$ taken with respect to the standard
metric.

\begin{lemma}
For every $k,n$ we have 
$$
\|\Omega_{k,n}\| < (2K_1)^{k-1}C^nb^2_{k,n},
$$
where $C = \max(C_0,C_2)$ is the bigger of the constants $C_0$,
$C_2$, and $b^2_{k,n}$ are the numbers defined recursively in
\ref{b.n.k.m}. 
\end{lemma}

\proof Use induction on $k$. The base of the induction is the case
$k=1$, where the estimate holds by \eqref{omega.est}. Assume the
estimate proved for all $\Omega_{p,n}$ with $p < k$, and fix a
number $n$. By \eqref{Omega.rec} we have
$$
\Omega_{k,n} = -\frac{1}{k}\sum_{1 \leq p \leq
k-1}\sum_q\sigma_{tot}(D^{tot}_{k-p,n-q}\Omega_{p,q}).
$$
Substituting the estimate \eqref{K1} together with the inductive
assumption and the estimate on $\|D^{tot}_{k-p,n-q}\|^2_q$ obtained in
Lemma~\ref{D2est}, we get
\begin{align*}
\|\Omega_{k,n}\| &< \frac{1}{k}\sum_{1 \leq p \leq k-1}\sum_q
\|\sigma_{tot}\|\cdot\|D^{tot}_{k-p,n-q}\|^2_q\cdot\|\Omega_{p,q}\| \\
&< \frac{1}{k}\sum_{1 \leq p \leq k-1}\sum_qK_1 \cdot
2(q+2(k-p))(3K)^{k-p-1}C_0^{n-q}b_{k-p,n-q} \\
&\qquad\qquad\qquad\qquad\cdot (2K_1)^{p-1}C^qb_{p,q} \\
&< \frac{1}{k}(2K_1)^{k-1}C^n\sum_{1 \leq p \leq
k-1}\sum_q(q+2(k-p))b_{k-p,n-q}b^2_{p,q} \\
&= (2K_1)^{k-1}C^nb^2_{k,n},
\end{align*}
which proves the lemma. 
\endproof 

\punkt This estimate immediately implies the last claim of
Theorem~\ref{converge}. Indeed, together with Lemma~\ref{combin.2}
it implies that for every $k,n$ 
$$
\|\Omega_{k,n}\| < (C_3)^n
$$
for some constant $C_3>0$. But $\Omega = \sum_n\sum_{0 \leq k \leq
n}\Omega_{k,n}$ is the Taylor series decomposition for the formal
polarization $\Omega$ at $0 \subset M$. The standard convergence
criterion shows that this series converges in an open ball of radius
$1/C_3 > 0$. Therefore the polarization $\Omega$ is indeed
real-analytic in a neighborhood of $0 \subset M \subset \tm$.
\endproof

\section*{Appendix}
\refstepcounter{section}
\refstepcounter{subsection}
\renewcommand{\thesection}{A}
\addcontentsline{toc}{section}{Appendix}

\punkt In this appendix we describe a well-known Borel-Weyl type
localization construction for quaternionic vector spaces (see,
e.g. \cite{HKLR}) which provides a different and somewhat more
geometric approach to many facts in the theory of Hodge
manifolds. In particular, we establish, following Deligne and
Simpson (\cite{De2}, \cite{Simpson}), a relation between Hodge
manifolds and the theory of mixed $\R$-Hodge structures. For the
sake of simplicity, we consider only Hodge manifold structures on
the formal neighborhood of $0 \in \R^{4n}$ instead of actual Hodge
manifolds, as in Section~\ref{convergence}. To save the space all
proofs are either omitted or only sketched.

\punkt Let $\SB$ be the Severi-Brauer variety associated to the
algebra $\h$, that is, the real algebraic variety of minimal right 
ideals in $\h$. The variety $\SB$ is a twisted $\R$-form of the
complex projective line $\cp^1$. 

For every algebra map $I:\C \to \h$ let the algebra $\h \otimes_\R
\C$ act on the $2$-dimensional complex vector space $\h_I$ by left
multiplication, and let $\widehat{I} \subset \h \otimes_\R \C$ be
the annihilator of the subspace $I(\C) \subset \h_I$ with respect to
this action. The subspace $\widehat{I}$ is a minimal right ideal in
$\h \otimes_\R \C$. Therefore it defines a $\C$-valued point
$\widehat{I} \subset \SB(\C)$ of the real algebraic variety
$\SB$. This establishes a bijection between the set $\SB(\C)$ and
the set of algebra maps from $\C$ to $\h$.

\punkt Let $\Shv(\SB)$ be the category of flat coherent sheaves on
$\SB$. Say that a sheaf $\E \in \Ob\Shv(\SB)$ is {\em of weight $p$}
if the sheaf $\E \otimes \C$ on $\cp^1 = \SB \otimes \C$ is a sum of
several copies of the sheaf $\calo(p)$. 

Consider a quaternionic vector space $V$.  Let $\I \in \h
\otimes \calo_\SB$ be the tautological minimal left ideal in the
algebra sheaf $\h \otimes \calo_\SB$, and let $\loc{V} \in
\Ob\Shv(\SB)$ be the sheaf defined by
$$
\loc{V} = V \otimes \calo_\SB / \I \cdot V \otimes \calo_\SB. 
$$
The correspondence $V \mapsto \loc{V}$ defines a functor from
quaternionic vector spaces to $\Shv(\SB)$. It is easy to check that
this functor is a full embedding, and its essential image is the
subcategory of sheaves of weight $1$. Call $\loc{V}$ {\em the
localization} of the quaternionic vector space $V$. For every
algebra map $\I:\C \to \h$ the fiber $\loc{V}|_{\widehat{I}}$ of the
localization $\loc{V}$ over the point $\widehat{I} \subset \SB(\C)$
corresponding to the map $i:\C \to \h$ is canonically isomorphic to
the real vector space $V$ with the complex structure $V_I$.

\punkt The compact Lie group $U(1)$ carries a canonical structure of
a real algebraic group. Fix an algebra embedding $I:\C \to \h$ and
let the group $U(1)$ act on the algebra $\h$ as in
\ref{u.acts.on.h}. This action is algebraic and induces therefore an
algebraic action of the group $U(1)$ on the Severi-Brauer variety
$\SB$. The point $\widehat{I}:\Spec\C \subset \SB$ is preserved by
the $U(1)$-action.  The action of the group $U(1)$ on the complement
$\SB \setminus \widehat{I}(\Spec\C) \subset \SB$ is free, so that
the variety $\SB$ consists of two $U(1)$-orbits. The corresponding
orbits of the complexified group $\C^* = U(1) \times \Spec\C$ on the
complexification $\SB \times \Spec\C \cong \cp$ are the pair of
points $0,\infty \subset \cp$ and the open complement $\cp \setminus
\{0,\infty\} \cong \C^* \subset \cp$.

Let $\Shv^{U(1)}(\SB)$ be the category of $U(1)$-equivariant flat
coherent sheaves on the variety $\SB$. The localization construction
immediately extends to give the equivalence $V \mapsto \loc{V}$
between the category of equivariant quaternionic vector spaces and
the full subcategory in $\Shv^{U(1)}(\SB)$ consisting of sheaves of
weight $1$. For an equivariant quaternionic vector space $V$, the
fibers of the sheaf $\loc{V}$ over the point $\widehat{I} \subset
\SB(\C)$ and over the complement $\SB \setminus
\widehat{I}(\Spec\C)$ are isomorphic to the space $V$ equipped,
respectively, with the preferred and the complementary complex
structures $V_I$ and $V_J$.

\punkt The category of $U(1)$-equivariant flat coherent sheaves on
the variety $\SB$ admits the following beautiful description, due to
Deligne. 

\begin{lemma}[ (\cite{De2},\cite{Simpson})]\label{DS}
\begin{enumerate}
\item For every integer $n$ the full subcategory
$$
\Shv^{U(1)}_n(\SB) \subset \Shv^{U(1)}(\SB)
$$ 
of sheaves of weight $n$ is equivalent to the category of pure
$\R$-Hodge structures of weight $n$.
\item The category of pairs $\langle \E, W_\cdot\rangle$ of a flat
$U(1)$-equivariant sheaf 
$$
\E \in \Ob\Shv^{U(1)}(\SB)
$$ 
and an increasing filtration $W_\cdot$ on $\E$ such that for every
integer $n$
\begin{equation}\label{hg}
W_n\E/W_{n-1}\E \quad \text{ is a sheaf of weight } \quad n \quad
\text{ on } \quad \SB
\end{equation}
is equivalent to the category of mixed $\R$-Hodge structures. (In
particular, it is abelian.)
\end{enumerate}
\end{lemma}

\punkt For every pure $\R$-Hodge structure $V$ call the
corresponding $U(1)$-e\-qui\-va\-ri\-ant flat coherent sheaf on the
variety $\SB$ {\em the localization} of $V$ and denote it by
$\loc{V}$. For the trivial $\R$-Hodge structure $\R(0)$ of weight
$0$ the sheaf $\loc{\R(0)}$ coincides with the structure sheaf
$\calo$ on $\SB$. If $V$, $W$ are two pure $\R$-Hodge structures,
then the space $\Hom(\loc{V},\loc{W})$ of $U(1)$-equivariant maps
between the corresponding sheaves coincides with the space of weakly
Hodge maps from $V$ to $W$ in the sense of Subsection~\ref{w.H.sub}.

For every pure $\R$-Hodge structure $V$ the space
$\Gamma(\SB,\loc{V})$ of the global sections of the sheaf $\loc{V}$
is equipped with an action of the group $U(1)$ and carries therefore
a canonical $\R$-Hodge structure of weight $0$. This $\R$-Hodge
structure is the same as the universal $\R$-Hodge structure
$\Gamma(V)$ of weight $0$ constructed in Lemma~\ref{g.ex}. This
explains our notation for the functor $\Gamma:\prehodge_{\geq 0} \to
\prehodge_0$.

\punkt Assume given a complex vector space $V$ and let $M$ be the
formal neighborhood of $0 \in V$. Let $\B^\cdot$ be the Weil algebra
of the manifold $M$, as in \ref{formal.Weil}. For every $n \geq 0$
the vector space $\B^n$ is equipped with an $\R$-Hodge structure of
weight $n$, so that we can consider the localization
$\loc{\B^n}$. The sheaf $\oplus\loc{\B^\cdot}$ is a commutative
algebra in the tensor category $\Shv^{U(1)}(\SB)$. We will call {\em
the localized Weil algebra}.

The augmentation grading on $\B^\cdot$ defined in \ref{aug} is
compatible with the $\R$-Hodge structures.  Therefore it defines an
augmentation grading on the localized Weil algebra
$\loc{\B^\cdot}$. The finer augmentation bigrading on $\B^\cdot$
does not define a bigrading on $\loc{\B^\cdot}$. However, it does
define a bigrading on the complexified algebra $\loc{\B^\cdot}
\otimes \C$ of $\C^*$-equivariant sheaves on the manifold $\SB
\otimes \C \cong \cp$.

\punkt Assume now given a flat extended connection $D:\B^\cdot \to
\B^{\cdot+1}$ on $M$. Since the derivation $D$ is weakly Hodge, it
corresponds to a derivation $D:\loc{\B^\cdot} \to
\loc{\B^{\cdot+1}}$ of the localized Weil algebra $\loc{\B^\cdot}$.
It is easy to check that the complex $\langle \loc{\B^\cdot}, D
\rangle$ is acyclic in all degrees but $0$. Denote by $\HH^0$ the
$0$-th cohomology sheaf $H^0(\loc{\B^\cdot})$. The sheaf $\HH^0$
carries a canonical algebra structure. Moreover, while the
derivation $D$ does not preserve the augmentation grading on
$\loc{\B^\cdot}$, it preserves the decreasing {\em augmentation
filtration} $\left(\loc{\B^\cdot}\right)_{\geq \cdot}$.  Therefore
we have a canonical decreasing filtration on the algebra $\HH^0$,
which we also call the augmentation filtration.

\punkt It turns out that the associated graded quotient algebra $\gr \HH^0$
with respect to the augmentation filtration does not depend on the
extended connection $D$. To describe it, introduce the $\R$-Hodge
structure $W$ of weight $-1$ by setting $W = V$ as a real vector
space and 
\begin{align}\label{W.dfn}
\begin{split}
W^{-1,0} &= V \subset V \otimes_\R \C,\\ 
W^{0,-1} &= \overline{V} \subset V \otimes_\R \C.
\end{split}
\end{align} 
The $k$-th graded piece $\gr_k \HH^0$ with respect to the
augmentation filtration is then isomorphic to the symmetric power
$S^k(\loc{W})$ of the localization $\loc{W}$ of the $\R$-Hodge
structure $W$. In particular, it is a sheaf of weight $-n$, so that
up to a change of numbering the augmentation filtration on $\HH^0$
satisfies the condition \eqref{hg}. The extension data between these
graded pieces depend on the extended connection $D$. The whole
associated graded algebra $\gr \HH^0$ is isomorphic to the completed
symmetric algebra $\widehat{S}^\cdot(\loc{W})$.

\punkt Using standard deformation theory, one can show that the
algebra map $\HH^0 \to \loc{\B^0}$ is the universal map from the
algebra $\HH^0$ to a complete commutative pro-algebra in the tensor
category of $U(1)$-equivariant flat coherent sheaves of weight $0$
on $\SB$. Moreover, the localized Weil algebra $\loc{\B^\cdot}$
coincides with the relative de Rham complex of $\loc{\B^0}$ over
$\HH^0$. Therefore one can recover, up to an isomorphism, the whole
algebra $\loc{\B^\cdot}$ and, consequently, the extended connection
$D$, solely from the algebra $\HH^0$ in $\Shv^{U(1)}(\SB)$. Together
with Lemma~\ref{DS}~\thetag{ii} this gives the following, due also
to Deligne (in a different form).

\begin{prop}\label{DS.prop}
The correspondence $\langle \B^\cdot, D \rangle \mapsto \HH^0$ is a
bijection between the set of all isomorphism classes of flat
extended connections on $M$ and the set of all algebras $\HH^0$ in
the tensor category of mixed $\R$-Hodge structures equipped with an
isomorphism $\gr^W_{-1}\HH^0 \cong W$ between the $-1$-th associated
graded piece of the weight filtration on $\HH^0$ and the pure
$\R$-Hodge structure $W$ defined in \eqref{W.dfn} which induces for
every $n \geq 0$ an isomorphism $\gr^W_{-n}\HH^0 \cong S^nW$.
\end{prop}

\rem The scheme $\Spec\HH^0$ over $\SB$ coincides with the so-called
twistor space of the manifold $\tm$ with the hypercomplex structure
given by the extended connection $D$ (see \cite{HKLR} for the
definition). Deligne's and Simpson's (\cite{De2}, \cite{Simpson})
approach differs from ours in that they use the language of twistor
spaces to describe the relation between $U(1)$-equivariant
hypercomplex manifolds and mixed $\R$-Hodge structures. Since this
requires some additional machinery, we have avoided introducing
twistor spaces in this paper.

\punkt We will now try to use the localization construction to
eludicate some of the complicated linear algebra used in
Section~\ref{main.section} to prove our main theorem. As we have
already noted, the category $\prehodge$ of pure $\R$-Hodge
structures with weakly Hodge maps as morphisms is identified by
localization with the category $\Shv^{U(1)}(\SB)$ of
$U(1)$-equivariant flat coherent sheaves on $\SB$. Moreover, the
functor $\Gamma:\prehodge_{\geq 0} \to \prehodge_0$ introduced in
Lemma~\ref{g.ex} is simply the functor of global sections
$\Gamma(\SB,\cdot)$. 

\punkt Consider the localized Weil algebra $\loc{\B^\cdot}$ with the
derivation $C:\loc{\B^\cdot} \to \loc{\B^{\cdot+1}}$ associated to
the canonical weakly Hodge derivation introduced in
\ref{C.and.sigma}. The differential graded algebra $\langle
\loc{\B^\cdot},C \rangle$ is canonically an algebra over the
completed symmetric algebra $\widehat{S}^\cdot(V)$ generated by the
constant sheaf on $\SB$ with the fiber $V$. Moreover, it is a free
commutative algebra generated by the complex 
\begin{equation}\label{cmplx}
V \longrightarrow V(1)
\end{equation}
placed in degrees $0$ and $1$, where $V(1)$ is the
$U(1)$-equivariant sheaf of weight $1$ on $\SB$ corresponding to the
$\R$-Hodge structure given by $V(1)^{1,0} = V$ and $V(1)^{0,1} =
\overline{V}$. 

\punkt The homology sheaves of the complex \eqref{cmplx} are
non-trivial only in degree $1$. This non-trivial homology sheaf is a
skyscraper sheaf concentrated in the point $\widehat{I}(\Spec\C)
\subset \SB$ with fiber $V$. The associated sheaf on the
complexification $\SB \otimes \C \cong \cp$ splits into the sum of
skyscraper sheaf with fiber $V$ concentrated at $0 \in \cp$ and the
skyscraper sheaf with fiber $\overline{V}$ concentrated at $\infty
\in \cp$. This splitting cooresponds to the splitting of the complex
\eqref{cmplx} itself into the components of augmentation bidegrees
$(1,0)$ and $(0,1)$. 

\punkt Let now $\I^\cdot = \B^\cdot_{\geq 0,\geq 0}$ be the sum of
the components in the Weil algebra $\B^\cdot$ of augmentation
bidegree greater or equal than $(1,1)$. The subspace $\I^\cdot
\subset \B^\cdot$ is compatible with the $\R$-Hodge structure. The
crucial point in the proof of Theorem~\ref{kal=ext} is
Proposition~\ref{ac}, which claim the acyclycity of the complex
$\langle \Gamma(\I^\cdot), C \rangle$. It is this fact that becomes
almost obvious from the point of view of the localization
construction. To show it, we first prove the following.

\begin{lemma}
The complex $\langle\loc{\I^\cdot},C\rangle$ of $U(1)$-equivariant
sheaves on $\SB$ is acyclic.
\end{lemma}

\proof It suffices to prove that the complex $\loc{I^\cdot} \otimes
\C$ of sheaves on $\SB \otimes \C \cong \cp$ is acyclic. To prove
it, let $p,q \geq 1$ be arbitrary integer, and consider the
component $\left(\loc{B^\cdot}\right)_{p,q}$ 
of augmentation bidegree $(p,q)$ in the localized Weil algebra
$\loc{\B^\cdot}$. By definition we have 
$$
\left(\loc{B^\cdot}\right)_{p,q} = \left(\loc{B^\cdot}\right)_{p,0}
\otimes \left(\loc{B^\cdot}\right)_{0,q} 
$$
Since the complex $\left(\loc{B^\cdot}\right)_{p,0} =
S^p\left(\loc{B^\cdot}\right)_{1,0}$ has homology concentrated at $0
\in \cp$, while the complex $\left(\loc{B^\cdot}\right)_{0,q} =
S^q\left(\loc{B^\cdot}\right)_{0,1}$ has homology concentrated at
$\infty \in \cp$, their product is indeed acyclic.
\endproof 

Now, we have $\Gamma(\I^\cdot) = \Gamma(\SB,\loc{\I^\cdot})$, and
the functor $\Gamma(\SB,\cdot)$ is exact on the full subcategory in
$\Shv^{U(1)}(\SB)$ consisting of sheaves of positive
weight. Therefore the complex $\langle \Gamma(\I^\cdot),C \rangle$
is also acyclic, which gives an alternative proof of
Proposition~\ref{ac}.

\punkt We would like to finish the paper with the following
observation. Proposition~\ref{DS.prop} can be extended to the
following claim.

\begin{prop}
Let $M$ be a complex manifold.  There exists a naturla bijection
between the set of isomorphism classes of germs of Hodge manifold
structures on $\tm$ in the neighborhood of the zero section $M
\subset \tm$ and the set of multiplicative filtrations $F^\cdot$ on
the sheaf $\calo_\R(M) \otimes \C$ of $\C$-valued real-analytic
functions on $M$ satisfying the following condition:
\begin{itemize}
\item For every point $m \in M$ let $\widehat{\calo}_m$ be the
formal completion of the local ring $\calo_m$ of germs of
real-analytic functions on $M$ in a neighborhood of $m$ with respect
to the maximal ideal. Consider the filtration $F^\cdot$ on
$\widehat{\calo}_m \otimes \C$ induced by the filtration $f^\cdot$
on the sheaf $\calo_\R(M) \otimes \C$, and for every $k \geq 0$ let
$W_{-k}\widehat{\calo}_m \subset \widehat{\calo}_m$ be the $k$-th
power of the maximal ideal in $\widehat{\calo}_m$. Then the triple
$\langle \widehat{\calo}_m, F^\cdot, W^\cdot \rangle$ is a mixed
$\R$-Hodge structure. (In particular, $F^k = 0$ for $k > 0$.)
\end{itemize}
\end{prop}

If the Hodge manifold structure on $\tm$ is such that the projection
$\rho:\tm_I \to M$ is holomorphic for the preferred complex
structure $\tm_I$ on $\tm$, then it is easy to see that the first
non-trivial piece $F^0\calo_\R(M) \otimes \C$ of the filtration
$F^\cdot$ on the sheaf $\calo_\R(M) \otimes \C$ coincides with the
subsheaf $\calo(M) \subset \calo_\R(M) \otimes \C$ of holomorphic
functions on $M$. Moreover, since the filtration $F^\cdot$ is
multiplicative, it is completely defined by the subsheaf
$F^{-1}\calo_\R(M) \otimes \C \subset \calo_\R(M) \otimes \C$. It
would be very interesting to find an explicit description of this
subsheaf in terms of the K\"ah\-le\-ri\-an connection $\nabla$ on $M$
which corresponds to the Hodge manifold structure on $\tm$.

\medskip

\noindent
{\sc Independent University of Moscow\\
B. Vlassievsky, 11\\
Moscow, Russia, 121002\\
{\em E-mail:}
kaledin$@$balthi.dnttm.rssi.ru}

\end{document}